\documentclass[10pt]{article}

\usepackage[doi=false,isbn=false,url=false,
	hyperref=auto,
	sorting=nyt,
	maxnames=10,
	maxcitenames=3,
	backend=biber,
	texencoding=auto,
	giveninits=true,
	block=space,
]{biblatex}

\DefineBibliographyStrings{english}{%
	andothers = {et al}
} 
 
\DeclareFieldFormat{eprint:arxiv}{%
	\href{https://arxiv.org/abs/#1}{\emph{arXiv:\allowbreak #1}}}

\DeclareFieldFormat{eprint:mrnumber}{%
	\href{http://www.ams.org/mathscinet-getitem?mr=MR#1}{MR\allowbreak #1}}

\DeclareFieldFormat{eprint:jstor}{%
	\href{http://www.jstor.org/stable/#1}{JSTOR\allowbreak #1}}

\DeclareFieldFormat{eprint:onlineshown}{%
	\em Available at \allowbreak \upshape \ttfamily \href{https://#1}{#1}}

\DeclareFieldFormat{eprint:onlinehidden}{%
	\em Available \href{https://#1}{online}}

\DeclareFieldFormat{eprint:inprep}{%
	\emph{In preparation}}

\DeclareFieldFormat{eprint:manual}{%
	\emph{#1}}

\DeclareFieldFormat{eprint:onarxiv}{%
	\emph{#1}}

\DeclareFieldFormat{eprint:toappear}{%
	To appear in \emph{#1}}

\DeclareFieldFormat{eprint:accepted}{%
	\emph{Accepted at {#1}; to appear}}

\DeclareSourcemap{
	\maps[datatype=bibtex]{
		\map{
			\step[fieldsource=mrnumber,		fieldtarget=eprint, final]
			\step[fieldset=eprinttype,		fieldvalue=mrnumber]
		}
		\map{
			\step[fieldsource=arxiv,		fieldtarget=eprint, final]
			\step[fieldset=eprinttype,		fieldvalue=arxiv]
		}
		\map{
			\step[fieldsource=jstor,		fieldtarget=eprint, final]
			\step[fieldset=eprinttype,		fieldvalue=jstor]
		}
		\map{
			\step[fieldsource=onlineshown,	fieldtarget=eprint, final]
			\step[fieldset=eprinttype,		fieldvalue=onlineshown]
		}
		\map{
			\step[fieldsource=onlinehidden,		fieldtarget=eprint, final]
			\step[fieldset=eprinttype,		fieldvalue=onlinehidden]
		}
		\map{
			\step[fieldsource=inprep,		fieldtarget=eprint, final]
			\step[fieldset=eprinttype,		fieldvalue=inprep]
		}
		\map{
			\step[fieldsource=manual,		fieldtarget=eprint, final]
			\step[fieldset=eprinttype,		fieldvalue=manual]
		}
		\map{
			\step[fieldsource=onarxiv,		fieldtarget=eprint, final]
			\step[fieldset=eprinttype,		fieldvalue=onarxiv]
		}
		\map{
			\step[fieldsource=toappear,		fieldtarget=eprint, final]
			\step[fieldset=eprinttype,		fieldvalue=toappear]
		}
		\map{
			\step[fieldsource=accepted,		fieldtarget=eprint, final]
			\step[fieldset=eprinttype,		fieldvalue=accepted]
		}
	}
}

\DeclareFieldFormat{doi}{%
	\href{https://doi.org/#1}{DOI}%
}

\DeclareFieldFormat{url}{%
	\href{https://#1}{\texttt{#1}}%
}

\DeclareFieldFormat{comments}{%
	\emph{#1}%
} 
 
\DeclareFieldFormat
	[article,inbook,incollection,inproceedings,patent,thesis,unpublished]
	{title}{{#1\isdot}}

\DeclareFieldFormat
	[article,book,inbook,incollection,inproceedings,patent,thesis,unpublished, online]
	{date}{\mkbibparens{#1}}
\DeclareFieldFormat
	[article,book,inbook,incollection,inproceedings,patent,thesis,unpublished]
	{volume}{\mkbibbold{#1}}
\DeclareFieldFormat
	{pages}{\mkbibparens{#1}}

\DeclareFieldFormat{edition}{%
	\ifinteger{#1}
	{\mkbibordedition{#1}~\bibstring{edition}\addcomma}
	{#1~\bibstring{edition}\addcomma}}

\renewbibmacro*{volume+number+eid}{%
	\printfield{volume}%
\setunit*{\adddot}%
	\printfield{number}%
\setunit{\space}%
	\printfield{eid}%
}

\DeclareBibliographyDriver{article}{%
	\usebibmacro{bibindex}%
	\usebibmacro{begentry}%
	\usebibmacro{author}%
\setunit{\addspace}%
	\usebibmacro{date}%
\newunit\newblock
	\usebibmacro{title}%
\newunit\newblock
	\printfield{journaltitle}\addperiod%
\setunit*{\addspace}%
	\usebibmacro{volume+number+eid}%
\setunit*{\addspace}\newblock
	\printfield{pages}
\setunit*{\addspace}
	\printfield{comments}%
	\usebibmacro{eprint}%
\setunit*{\addspace}%
	\printfield{doi}%
	\usebibmacro{finentry}%
}

\DeclareBibliographyDriver{online}{%
	\usebibmacro{bibindex}%
	\usebibmacro{begentry}%
	\usebibmacro{author}%
\setunit{\addspace}%
	\usebibmacro{date}%
\newunit\newblock
	\usebibmacro{title}%
\newunit\newblock
	\usebibmacro{eprint}%
	\usebibmacro{finentry}%
}

\DeclareBibliographyDriver{book}{%
	\usebibmacro{bibindex}%
	\usebibmacro{begentry}%
	\usebibmacro{author}%
\setunit{\addspace}\newblock
	\usebibmacro{date}%
\newunit\newblock
	\usebibmacro{title}%
\setunit*{\addspace}\newblock
	\printfield{edition}%
\newunit\newblock
	\printlist{publisher}
\setunit*{\addspace}%
	\printfield{volume}%
\setunit*{\addspace}\newblock%
	\printfield{comments}%
	\usebibmacro{eprint}%
\setunit*{\addspace}
	\printfield{doi}
	\usebibmacro{finentry}%
}

\DeclareBibliographyDriver{thesis}{%
	\usebibmacro{bibindex}%
	\usebibmacro{begentry}%
	\usebibmacro{author}%
\setunit{\addspace}\newblock
	\usebibmacro{date}%
\newunit\newblock
	\usebibmacro{title}.%
\setunit*{\addspace}\newblock
	\space Thesis, \printlist{institution}
	\usebibmacro{eprint}%
	\printfield{comments}%
\setunit*{\addspace}
	\printfield{doi}%
	\usebibmacro{finentry}%
}

\DeclareBibliographyDriver{inproceedings}{%
	\usebibmacro{bibindex}%
	\usebibmacro{begentry}%
	\usebibmacro{author}%
\setunit{\addspace}%
	\usebibmacro{date}%
\newunit\newblock
	\usebibmacro{title}%
\newunit\newblock
	\printfield{booktitle}\addcomma%
\setunit*{\addspace}%
	\printfield{series}\addcomma%
\setunit*{\addspace}\newblock
	\usebibmacro{editor}\addcomma
\setunit*{\addspace}\newblock
	\printlist{publisher}
\setunit*{\addspace}%
	\printfield{volume}%
\setunit*{\addspace}%
	\printfield{pages}
\setunit*{\addspace}
	\usebibmacro{eprint}%
	\printfield{comments}%
\setunit*{\addnnspace}
	\printfield{doi}
	\usebibmacro{finentry}%
}

\DeclareBibliographyDriver{inbook}{%
	\usebibmacro{bibindex}%
	\usebibmacro{begentry}%
	\usebibmacro{author}%
\setunit{\addspace}%
	\usebibmacro{date}%
\newunit\newblock
	\usebibmacro{title}%
\newunit\newblock
	\printfield{booktitle}\addcomma%
\setunit*{\addspace}%
	\printfield{series}\addcomma%
\setunit*{\addspace}\newblock
	\usebibmacro{editor}\addcomma
\setunit*{\addspace}\newblock
	\printlist{publisher}
\setunit*{\addspace}%
	\printfield{volume}%
\setunit*{\addspace}%
	\printfield{pages}
\setunit*{\addspace}
	\usebibmacro{eprint}%
	\printfield{comments}%
\setunit*{\addspace}
	\printfield{doi}
	\usebibmacro{finentry}%
}

\DeclareBibliographyDriver{incollection}{%
	\usebibmacro{bibindex}%
	\usebibmacro{begentry}%
	\usebibmacro{author}%
\setunit{\addspace}%
	\usebibmacro{date}%
\newunit\newblock
	\usebibmacro{title}%
\newunit\newblock
	\printfield{booktitle}\addcomma%
\setunit*{\addspace}%
	\printfield{series}\addcomma%
\setunit*{\addspace}\newblock
	\printlist{publisher}
\setunit*{\addspace}%
	\printfield{volume}%
\setunit*{\addspace}%
	\printfield{pages}
\setunit*{\addspace}
	\usebibmacro{eprint}%
	\printfield{comments}%
\setunit*{\addspace}
	\printfield{doi}
	\usebibmacro{finentry}%
}

\AtEveryBibitem{%
	\clearfield{day}%
	\clearfield{month}%
} 
 
\DeclareNameAlias{sortname}{given-family}
\DeclareNameAlias{default}{given-family}

\newcommand{\acknofootnote}{%
	\blfootnote{%
		\hspace*{1em}%
		Jonathan Hermon%
	\hfill%
		Sam Olesker-Taylor%
		\hspace*{1em}%
	\\%
		\href{mailto:jhermon@math.ubc.ca}{jhermon@math.ubc.ca},
		\href{http://www.math.ubc.ca/~jhermon/}{math.ubc.ca/$\sim$jhermon/}%
	\hfill%
		\href{mailto:oleskertaylor.sam@gmail.com}{oleskertaylor.sam@gmail.com},
		\href{https://sites.google.com/view/sam-ot/}{sites.google.com/view/sam-ot/}%
	\\%
		University of British Columbia, Vancouver, Canada%
	\hfill%
		Department of Mathematical Sciences, University of Bath, UK%
	\\%
		Supported by EPSRC EP/L018896/1 and an NSERC Grant%
	\hfill%
		Supported by EPSRC Grants 1885554 and EP/N004566/1%
	\par\smallskip\par
	\centering%
		The vast majority of this work was undertaken whilst both authors were at the University of Cambridge%
	}
}

\newcommand*{\mm}{\ensuremath{L}}

\newcommand*{\qq}{\ensuremath{q}}

\newcommand*{\RR}{\ensuremath{\mathcal R}}
\newcommand*{\LL}{\ensuremath{L}}

\newcommand*{\WW}{\ensuremath{W}}
\newcommand*{\ww}{\ensuremath{w}}
\newcommand*{\VV}{\ensuremath{V}}
\newcommand*{\vv}{\ensuremath{v}}

\newcommand{\printtoc}[1]{%
	\ifthenelse%
		{\equal{#1}{1}}%
		{\sffamily\boldmath\tableofcontents\unboldmath\normalfont}%
		{\newpage\small\sffamily\boldmath\tableofcontents\unboldmath\normalfont\normalsize}%
	}

\newcommand{\nextresult}{%
	\setcounter{introthm}{\value{introthm}}
	\setcounter{introcor}{\value{introthm}}
	\setcounter{introconj}{\value{introthm}}
	\setcounter{introdefn}{\value{introthm}}
	\setcounter{intrormkT}{\value{introthm}}
}

\usepackage[utf8]{inputenc}

\usepackage{empheq}

\usepackage{manyfoot}
\SetFootnoteHook{\hspace*{-1.8em}}
\DeclareNewFootnote{bl}[gobble]
\newcommand{\blfootnote}[1]{\footnotebl{\sffamily#1}}

\usepackage{graphicx}
\graphicspath{{../Figures/}}
\usepackage[labelsep=period]{caption}

\usepackage{chngcntr}
\counterwithin{figure}{section}

\usepackage[UKenglish]{babel}
\usepackage{amsmath, amsthm, amssymb, mathrsfs, bm, mathtools}
\usepackage{wasysym}
\usepackage{mathtools}

\allowdisplaybreaks

\usepackage{cases}

\usepackage{titlesec}
\usepackage{titletoc}

\titleformat{\title}
{\sffamily \Huge \bfseries \scshape \boldmath}{}{}{}
\titleformat{\section}
{\sffamily \Large \bfseries \scshape \boldmath}{\thesection}{1em}{}
\titleformat{\subsection}
{\sffamily \large \bfseries \scshape \boldmath}{\thesubsection}{1em}{}
\titleformat{\subsubsection}
{\sffamily \normalsize \bfseries \scshape \boldmath}{\thesubsubsection}{1em}{}
\titleformat{\paragraph}
{\sffamily \normalsize \bfseries \scshape \boldmath}{\theparagraph}{1em}{}
\titleformat{\subparagraph}[runin]
{\sffamily \normalsize \bfseries \scshape \boldmath}{\thesubparagraph}{1em}{}

\def\IfAmpersandUseAlign#1#2&#3\EndIfAmpersandUseAlign
{%
	\if\relax\detokenize{#3}\relax
	\begin{equation*}%
	#1%
	\end{equation*}%
	\else
	\begin{align*}%
	#1%
	\end{align*}%
	\fi
}
\def\[#1\]%
{%
	\IfAmpersandUseAlign{#1}#1&\EndIfAmpersandUseAlign
}

\usepackage{eqparbox}

\usepackage{scrextend}

\usepackage{enumitem}

\newcommand{\numberingroman}{%
	\renewcommand{\labelenumi}{(\roman{enumi})}%
	\renewcommand{\theenumi}{(\roman{enumi})}%
}

\setlist[description]{%
	topsep		= 0pt,		
	noitemsep,				
	font		= {\mdseries\itshape},	
}

\usepackage[dvipsnames,hyperref]{xcolor}

\newcommand{\nt}{\addtocounter{equation}{1}\tag{\theequation}}

\newcommand{\bcdot}{\ensuremath{\bm{\cdot}}}
\let\mod\relax
\DeclareMathOperator{\mod}{\, mod}

\newcommand{\Quad}[1]{
	\mathchoice
	{\quad\text{#1}\quad}
	{\text{ #1 }}
	{\text{ #1 }}
	{\text{ #1 }}
}

\newcommand{\whp}{\text{whp}\xspace}

\newcommand{\Qforall}{\Quad{for all}}
\newcommand{\Qfor}{\Quad{for}}
\newcommand{\Qand}{\Quad{and}}
\newcommand{\Qwhere}{\Quad{where}}
\newcommand{\Qwhen}{\Quad{when}}

\newcommand{\id}{\mathsf{id}}

\newcommand{\typ}{\mathsf{typ}}

\newcommand{\cq}{\coloneqq}
\renewcommand{\epsilon}{\varepsilon}

\newcommand{\eps}{\epsilon}

\usepackage{xifthen}

\newcommand{\binomt}[2]{ \textstyle \binom{#1}{#2} \displaystyle }

\newcommand{\maxt}[1]{ \textstyle \max_{#1} \displaystyle }

\newcommand{\mint}[1]{ \textstyle \min_{#1} \displaystyle }

\newcommand{\limt}[1]{ \textstyle \lim_{#1} \displaystyle }
\newcommand{\limsupt}[1]{ \textstyle \limsup_{#1} \displaystyle }

\usepackage{calc}
\newlength{\halfplusheight}
\newcommand{\MAX}[1]{\mathop{\raisebox{\halfplusheight}{\(\displaystyle\max_{#1}\)}}}
\newcommand{\MIN}[1]{\mathop{\raisebox{\halfplusheight}{\(\displaystyle\min_{#1}\)}}}
\newcommand{\LIM}[1]{\mathop{\raisebox{\halfplusheight}{\(\displaystyle\lim_{#1}\)}}}
\newcommand{\LIMSUP}[1]{\mathop{\raisebox{\halfplusheight}{\(\displaystyle\limsup_{#1}\)}}}
\newcommand{\LIMINF}[1]{\mathop{\raisebox{\halfplusheight}{\(\displaystyle\liminf_{#1}\)}}}

\DeclareMathOperator*{\sumTT}{\textstyle\sum}
\newcommand{\sumT}[2][]{
	\ifthenelse{\isempty{#1}}
	{\sumTT_{#2}}
	{\sumTT_{#2}^{#1}}
}
\newcommand{\sumt}[2][]{
	\ifthenelse{\isempty{#1}}
	{\textstyle \sum_{#2}      \displaystyle}
	{\textstyle \sum_{#2}^{#1} \displaystyle}
}
\newcommand{\sumd}[2][]{
	\ifthenelse{\isempty{#1}}
	{\displaystyle \sum_{#2}}
	{\displaystyle \sum_{#2}^{#1}}
}

\newcommand{\intt}[2][]{
	\ifthenelse{\isempty{#1}}
	{\textstyle \int_{#2}      \displaystyle}
	{\textstyle \int_{#2}^{#1} \displaystyle}
}

\newcommand{\prodt}[2][]{
	\ifthenelse{\isempty{#1}}
	{\textstyle \prod_{#2}      \displaystyle}
	{\textstyle \prod_{#2}^{#1} \displaystyle}
}
\newcommand{\prodd}[2][]{
	\ifthenelse{\isempty{#1}}
	{\prod_{#2}}
	{\prod_{#2}^{#1}}
}

\usepackage{xspace}

\let\originalexp\exp
\let\exp\relax

\DeclareRobustCommand{\exp} [1]{\originalexp(#1)}
\newcommand{\expb} [1]{\originalexp\bigl( #1 \bigr)}

\newcommand{\abs}  [1]{| #1 |}
\newcommand{\absb} [1]{\bigl| #1 \bigr|}

\newcommand{\norm}  [1]{\lVert #1 \rVert}
\newcommand{\normb} [1]{\big\lVert #1 \bigr\rVert}

\newcommand{\rbr} [1]{ ( #1 ) }
\newcommand{\rbb} [1]{\bigl( #1 \bigr)}

\newcommand{\sbb} [1]{\bigl[ #1 \bigr]}

\newcommand{\bra} [1]{ \{ #1 \} }
\newcommand{\brb} [1]{\bigl\{ #1 \bigr\}}

\newcommand{\dist}{\mathrm{dist}}

\DeclareMathOperator{\diam}{diam}

\newcommand{\mix}{\mathrm{mix}}

\newcommand{\rel}{\mathrm{rel}}

\newcommand{\st}{{ \ \mathrm{st} \ }}

\newcommand{\ab} {\mathrm{ab}}

\DeclareMathOperator{\step}{step}

\newcommand{\Unif}{\mathrm{Unif}}
\newcommand{\iid}{\mathrm{iid}}

\newcommand{\Geom}{\mathrm{Geom}}

\newcommand{\tmix}{t_\mix}

\newcommand{\trel}{t_\rel}

\newcommand{\Ninn}{{N\in\mathbb{N}}}

\newcommand{\floor}[1]{\lfloor #1 \rfloor}

\newcommand{\ceil}[1]{\lceil #1 \rceil}

\newcommand{\midb}{\bigm\vert}

\newcommand{\one}  [1]{\bm1( #1 )}
\newcommand{\oneb} [1]{\bm1\bigl( #1 \bigr)}

\newcommand{\logk}[1][]{
	\ifthenelse{\equal{}{#1}}
	{\log k}
	{(\log k)^{#1}}
}

\newcommand{\logn}[1][]{
	\ifthenelse{\equal{}{#1}}
	{\log n}
	{(\log n)^{#1}}
}

\newcommand{\logm}[1][]{
	\ifthenelse{\equal{}{#1}}
	{\log m}
	{(\log m)^{#1}}
}

\newcommand{\loglogn}[1][]{
	\ifthenelse{\equal{}{#1}}
	{\log\log n}
	{(\log\log n)^{#1}}
}

\newcommand{\prt}[2][]{
	\ifthenelse{\equal{}{#1}}
	{\mathbb{P}(#2)}
	{\mathbb{P}_{#1}(#2)}
}
\newcommand{\pr}[2][]{
	\mathchoice
	{\ifthenelse{\isempty{#1}}
		{\mathbb{P}\bigl(#2\bigr)}
		{\mathbb{P}_{#1}\bigl(#2\bigr)}}
	{\ifthenelse{\isempty{#1}}
		{\mathbb{P}(#2)}
		{\mathbb{P}_{#1}(#2)}}
	{\ifthenelse{\isempty{#1}}
		{\mathbb{P}(#2)}
		{\mathbb{P}_{#1}(#2)}}
	{\ifthenelse{\isempty{#1}}
		{\mathbb{P}(#2)}
		{\mathbb{P}_{#1}(#2)}}
}
\newcommand{\prb}[2][]{
	\ifthenelse{\equal{}{#1}}
	{\mathbb{P}\bigl( #2 \bigr)}
	{\mathbb{P}_{#1}\bigl( #2 \bigr)}
}
\newcommand{\prB}[2][]{
	\ifthenelse{\equal{}{#1}}
	{\mathbb{P}\Bigl( #2 \Bigr)}
	{\mathbb{P}_{#1}\Bigl( #2 \Bigr)}
}
\newcommand{\prbb}[2][]{
	\ifthenelse{\equal{}{#1}}
	{\mathbb{P}\biggl( #2 \biggr)}
	{\mathbb{P}_{#1}\biggl( #2 \biggr)}
}
\newcommand{\prBB}[2][]{
	\ifthenelse{\equal{}{#1}}
	{\mathbb{P}\Biggl( #2 \Biggr)}
	{\mathbb{P}_{#1}\Biggl( #2 \Biggr)}
}
\newcommand{\prs}[2][]{
	\ifthenelse{\equal{}{#1}}
	{\mathbb{P}\left( #2 \right)}
	{\mathbb{P}_{#1}\left( #2 \right)}
}

\newcommand{\qr}[2][]{
	\mathchoice
	{\ifthenelse{\isempty{#1}}
		{\mathbb{Q}\bigl(#2\bigr)}
		{\mathbb{Q}_{#1}\bigl(#2\bigr)}}
	{\ifthenelse{\isempty{#1}}
		{\mathbb{Q}(#2)}
		{\mathbb{Q}_{#1}(#2)}}
	{\ifthenelse{\isempty{#1}}
		{\mathbb{Q}(#2)}
		{\mathbb{Q}_{#1}(#2)}}
	{\ifthenelse{\isempty{#1}}
		{\mathbb{Q}(#2)}
		{\mathbb{Q}_{#1}(#2)}}
}
\newcommand{\qrb}[2][]{
	\ifthenelse{\equal{}{#1}}
	{\mathbb{Q}\bigl( #2 \bigr)}
	{\mathbb{Q}_{#1}\bigl( #2 \bigr)}
}
\newcommand{\qrB}[2][]{
	\ifthenelse{\equal{}{#1}}
	{\mathbb{Q}\Bigl( #2 \Bigr)}
	{\mathbb{Q}_{#1}\Bigl( #2 \Bigr)}
}
\newcommand{\qrbb}[2][]{
	\ifthenelse{\equal{}{#1}}
	{\mathbb{Q}\biggl( #2 \biggr)}
	{\mathbb{Q}_{#1}\biggl( #2 \biggr)}
}
\newcommand{\qrBB}[2][]{
	\ifthenelse{\equal{}{#1}}
	{\mathbb{Q}\Biggl( #2 \Biggr)}
	{\mathbb{Q}_{#1}\Biggl( #2 \Biggr)}
}
\newcommand{\qrs}[2][]{
	\ifthenelse{\equal{}{#1}}
	{\mathbb{Q}\left( #2 \right)}
	{\mathbb{Q}_{#1}\left( #2 \right)}
}

\newcommand{\ext}[2][]{
\ifthenelse{\equal{}{#1}}
{\mathbb{E}(#2)}
{\mathbb{E}_{#1}(#2)}
}
\newcommand{\ex}[2][]{
	\mathchoice
	{\ifthenelse{\isempty{#1}}
		{\mathbb{E}\bigl(#2\bigr)}
		{\mathbb{E}_{#1}\bigl(#2\bigr)}}
	{\ifthenelse{\isempty{#1}}
		{\mathbb{E}(#2)}
		{\mathbb{E}_{#1}(#2)}}
	{\ifthenelse{\isempty{#1}}
		{\mathbb{E}(#2)}
		{\mathbb{E}_{#1}(#2)}}
	{\ifthenelse{\isempty{#1}}
		{\mathbb{E}(#2)}
		{\mathbb{E}_{#1}(#2)}}
}
\newcommand{\exb}[2][]{
	\ifthenelse{\equal{}{#1}}
	{\mathbb{E}\bigl( #2 \bigr)}
	{\mathbb{E}_{#1}\bigr( #2 \bigr)}
}
\newcommand{\exB}[2][]{
	\ifthenelse{\equal{}{#1}}
	{\mathbb{E}\Bigl( #2 \Bigr)}
	{\mathbb{E}_{#1}\Bigl( #2 \Bigr)}
}
\newcommand{\exbb}[2][]{
	\ifthenelse{\equal{}{#1}}
	{\mathbb{E}\biggl( #2 \biggr)}
	{\mathbb{E}_{#1}\biggl( #2 \biggr)}
}
\newcommand{\exBB}[2][]{
	\ifthenelse{\equal{}{#1}}
	{\mathbb{E}\Biggl( #2 \Biggr)}
	{\mathbb{E}_{#1}\Biggl( #2 \Biggr)}
}

\newcommand{\fx}[2][]{
	\ifthenelse{\equal{}{#1}}
	{\mathbb{F}(#2)}
	{\mathbb{F}_{#1}(#2)}
}
\newcommand{\fxb}[2][]{
	\ifthenelse{\equal{}{#1}}
	{\mathbb{F}\bigl( #2 \bigr)}
	{\mathbb{F}_{#1}\bigr( #2 \bigr)}
}
\newcommand{\fxB}[2][]{
	\ifthenelse{\equal{}{#1}}
	{\mathbb{F}\Bigl( #2 \Bigr)}
	{\mathbb{F}_{#1}\Bigl( #2 \Bigr)}
}
\newcommand{\fxbb}[2][]{
	\ifthenelse{\equal{}{#1}}
	{\mathbb{F}\biggl( #2 \biggr)}
	{\mathbb{F}_{#1}\biggl( #2 \biggr)}
}
\newcommand{\fxBB}[2][]{
	\ifthenelse{\equal{}{#1}}
	{\mathbb{F}\Biggl( #2 \Biggr)}
	{\mathbb{F}_{#1}\Biggl( #2 \Biggr)}
}

\newcommand{\Var}[1]{\mathbb{V}\mathrm{ar}(#1)}
\newcommand{\Varb}[2][]{
	\ifthenelse{\equal{}{#1}}
	{\mathbb{V}\mathrm{ar} \bigl(#2\bigr)}
	{\mathbb{V}\mathrm{ar}_{#1} \bigl(#2\bigr)}
}
\newcommand{\VAR}[2][]{
	\ifthenelse{\equal{}{#1}}
	{\mathrm{Var}(#2)}
	{\mathrm{Var}_{#1}(#2)}
}

\newcommand{\Oh}  [1]{\mathcal{O}( #1 )}

\newcommand{\oh}  [1]{o( #1 )}

\newcommand{\mba}{\mathbb{A}}
\newcommand{\mbb}{\mathbb{B}}

\newcommand{\mbn}{\mathbb{N}}

\newcommand{\mbp}{\mathbb{P}}

\newcommand{\mbr}{\mathbb{R}}

\newcommand{\mbz}{\mathbb{Z}}

\newcommand{\mcb}{\mathcal{B}}

\newcommand{\mcd}{\mathcal{D}}
\newcommand{\mce}{\mathcal{E}}

\newcommand{\mcg}{\mathcal{G}}

\newcommand{\mci}{\mathcal{I}}
\newcommand{\mcj}{\mathcal{J}}

\newcommand{\mcl}{\mathcal{L}}
\newcommand{\mcm}{\mathcal{M}}

\newcommand{\mcr}{\mathcal{R}}
\newcommand{\mcs}{\mathcal{S}}

\newcommand{\mcw}{\mathcal{W}}

\newcommand{\mfd}{\mathfrak{D}}

\newcommand{\mfgcd}{\mathfrak{g}}

\newcommand{\mfr}{\mathfrak{R}}

\newcommand{\toinf}[1]{\ensuremath{#1\to\infty}}
\newcommand{\asinf}[1]{\text{as \ensuremath{#1\to\infty}}}

\newcommand{\kinf}{{k\to\infty}}

\newcommand{\Ninf}{{N\to\infty}}

\newcommand{\tinf}{{t\to\infty}}


\newcounter{parentnumber}

\makeatletter
\newenvironment{subtheorem-num}[1]{%
	\def\subtheoremcounter{#1}%
	\refstepcounter{#1}%
	\protected@edef\theparentnumber{\csname the#1\endcsname}%
	\setcounter{parentnumber}{\value{#1}}%
	\setcounter{#1}{0}%
	\expandafter\def\csname the#1\endcsname{\theparentnumber.\arabic{#1}}%
	\expandafter\def\csname theH#1\endcsname{thm.\theparentnumber.\arabic{#1}}%
	\unskip\ignorespaces
}{%
	\setcounter{\subtheoremcounter}{\value{parentnumber}}%
	\ignorespacesafterend
}
\makeatother

\usepackage[colorlinks, pdfauthor = {Jonathan\ Hermon\ and\ Sam\ Olesker-Taylor}, pdfusetitle]{hyperref}
\hypersetup{
	colorlinks,
	linkcolor	= {blue},
	filecolor	= {red},
	urlcolor	= {blue},
	citecolor	= {ForestGreen}
}

\newcommand{\qedtriangle}{\renewcommand{\qedsymbol}{\ensuremath{\triangle}}}

\usepackage[capitalise]{cleveref}

\crefformat{equation}{\textup{(#2#1#3)}}
\crefmultiformat{equation}{\textup{(#2#1#3}}{\textup{, #2#1#3)}}{\textup{, #2#1#3}}{\textup{, #2#1#3)}}
\crefrangeformat{equation}{\textup{(#3#1#4--#5#2#6)}}

\newenvironment{Proof}[1][\proofname]{%
	\proof[\upshape\bfseries\sffamily\boldmath#1]
}{\endproof}

\usepackage{etoolbox}
\usepackage{needspace}

\newtheoremstyle{sfsl}
{1\baselineskip}		
{1\baselineskip}		
{\slshape}				
{}						
{\bfseries\sffamily}	
{.}						
{0.5em}					
{\thmname{#1}\thmnumber{ #2}\thmnote{ {\mdseries(#3)}}}

\newtheoremstyle{sfup}
{1\baselineskip}		
{1\baselineskip}		
{\upshape}				
{}						
{\bfseries\sffamily}	
{.}						
{0.5em}					
{\thmname{#1}\thmnumber{ #2}\thmnote{ {\mdseries(#3)}}}

\theoremstyle{sfsl}

\newtheorem*{thm*}{Theorem}
\newtheorem{thm} {Theorem}[section]
\crefname{thm}{Theorem}{Theorems}

\newtheorem*{introthm*}{Theorem}
\newtheorem{introthm}{Theorem}

\crefname{introthm}{Theorem}{Theorems}

\newtheorem*{cor*}{Corollary}
\newtheorem{cor} [thm]{Corollary}
\crefname{cor}{Corollary}{Corollaries}

\newtheorem*{introcor*}{Corollary}
\newtheorem{introcor}{Corollary}

\crefname{introcor}{Corollary}{Corollaries}

\newtheorem*{introconj*}{Conjecture}

\crefname{introconj}{Conjecture}{Conjectures}

\newtheorem*{introques*}{Question}

\crefname{introques}{Question}{Questions}

\newtheorem*{lem*}    {Lemma}
\newtheorem{lem} [thm]{Lemma}
\crefname{lem}{Lemma}{Lemmas}

\newtheorem*{introlem*}{Lemma}

\crefname{introlem}{Lemma}{Lemmas}

\newtheorem*{prop*}    {Proposition}
\newtheorem{prop} [thm]{Proposition}
\crefname{prop}{Proposition}{Propositions}

\newtheorem*{clm*}    {Claim}

\crefname{clm}{Claim}{Claims}

\newtheorem*{defn*}    {Definition}
\newtheorem{defn} [thm]{Definition}
\crefname{defn}{Definition}{Definitions}

\newtheorem*{introdefn*}{Definition}
\newtheorem{introdefn}{Definition}

\crefname{introdefn}{Definition}{Definitions}

\newtheorem{innercustomgenericsl}{\customgenericnamesl}
\providecommand{\customgenericnamesl}{}
\newcommand{\newcustomtheoremsl}[2]{%
	\newenvironment{#1}[1]
	{%
		\renewcommand\customgenericnamesl{#2}%
		\renewcommand\theinnercustomgenericsl{##1}%
		\innercustomgenericsl
	}
	{\endinnercustomgenericsl}
}

\newcustomtheoremsl{customthm}{Theorem}
\newcustomtheoremsl{customprop}{Proposition}
\newcustomtheoremsl{customlemma}{Lemma}
\newcustomtheoremsl{customrmk}{Remark}
\newcustomtheoremsl{customconj}{Conjecture}
\newcustomtheoremsl{customdefn}{Definition}
\newcustomtheoremsl{customquestion}{Question}
\newcustomtheoremsl{customopenproblem}{Open Problem}
\newcustomtheoremsl{customhyp}{Hypothesis}

\newtheorem*{conj*}   {Conjecture}

\crefname{conj}{Conjecture}{Conjectures}

\newenvironment{conj-ind*}
	{\begin{quote}\textsf{\textbf{Conjecture.}}\slshape}
	{\end{quote}}
\newenvironment{conj-ind}
	{\begin{quote}\vspace{-\glueexpr\baselineskip+\topsep}\begin{customconj}}
	{\end{customconj}\end{quote}}

\newenvironment{question-ind*}
	{\begin{quote}\textsf{\textbf{Question.}}\slshape}
	{\end{quote}}
\newenvironment{question-ind}
	{\begin{quote}\vspace{-\glueexpr\baselineskip+\topsep}\begin{customquestion}}
	{\end{customquestion}\end{quote}}

\newenvironment{openproblem-ind*}
	{\begin{quote}\textsf{\textbf{Open Problem.}}\slshape}
	{\end{quote}}
\newenvironment{openproblem-ind}
	{\begin{quote}\vspace{-\glueexpr\baselineskip+\topsep}\begin{customopenproblem}}
	{\end{customopenproblem}\end{quote}}

\newtheorem*{hypothesis*}{Hypothesis}

\newtheorem*{hyp*}{Hypothesis}
\newtheorem{hyp}{Hypothesis}
\renewcommand*{\thehyp}{\Alph{hyp}}
\crefname{hyp}{Hypothesis}{Hypotheses}

\newtheorem*{rmk*}{Remark}

\theoremstyle{sfup}

\newtheorem{innercustomgenericup}{\customgenericnameup}
\providecommand{\customgenericnameup}{}
\newcommand{\newcustomtheoremup}[2]{%
	\newenvironment{#1}[1]
	{%
		\renewcommand\customgenericnameup{#2}%
		\renewcommand\theinnercustomgenericup{##1}%
		\innercustomgenericup
	}
	{\endinnercustomgenericup}
}

\crefname{exm} {Example}{Examples}
\crefname{exmT}{Example}{Examples}

\newtheorem{rmkT}[thm]{Remark}
	\newenvironment{rmkt}
	{\pushQED{\qed}\renewcommand{\qedsymbol}{\ensuremath{\triangle}}\rmkT}
	{\popQED\endrmkT}
\crefname{rmk} {Remark}{Remarks}
\crefname{rmkT}{Remark}{Remarks}

\newtheorem*{rmkTT}{Remark}
\newenvironment{rmkt*}
	{\pushQED{\qed}\renewcommand{\qedsymbol}{\ensuremath{\triangle}}\rmkTT}
	{\popQED\endrmkTT}

\newcustomtheoremup{customrmkT}{Remark}

\crefname{rmks} {Remarks}{Remarks}
\crefname{rmksT}{Remarks}{Remarks}

\newtheorem*{rmks*} {Remarks}
\newtheorem*{rmksTT}{Remarks}
\newenvironment{rmkst*}
	{\pushQED{\qed}\renewcommand{\qedsymbol}{\ensuremath{\triangle}}\rmksTT}
	{\popQED\endrmksTT}

\newtheorem{intrormkT}{Remark}
	\newenvironment{intrormkt}
	{\pushQED{\qed}\renewcommand{\qedsymbol}{\ensuremath{\triangle}}\intrormkT}
	{\popQED\endintrormkT}

\crefname{intrormk} {Remark}{Remarks}
\crefname{intrormkT}{Remark}{Remarks}

\newtheorem*{intrormk*} {Remark}
\newtheorem*{intrormkTT}{Remark}
\newenvironment{intrormkt*}
	{\pushQED{\qed}\renewcommand{\qedsymbol}{\ensuremath{\triangle}}\intrormkTT}
	{\popQED\endintrormkTT}

\newtheorem*{exm*} {Example}
\newtheorem*{exmTT}{Lemma}
	\newenvironment{exmt*}
	{\pushQED{\qed}\renewcommand{\qedsymbol}{\ensuremath{\triangle}}\exmTT}
	{\popQED\endexmTT}

\newtheorem*{note*} {Note}
\newtheorem*{noteTT}{Note}
	\newenvironment{notet*}
	{\pushQED{\qed}\renewcommand{\qedsymbol}{\ensuremath{\triangle}}\noteTT}
	{\popQED\endnoteTT}

\catcode`\^^A=14

\usepackage{xparse}

\newcounter{mixedsubequations}

\ExplSyntaxOn

\int_new:N \g_mixedsubeq_int

\NewDocumentEnvironment{mixedsubequations}{o}
{
	\IfNoValueTF { #1 }
	{
		\addtocounter{equation}{-\g_mixedsubeq_int}
		\stepcounter{mixedsubequations}
	}
	{
		\int_gset:Nn \g_mixedsubeq_int { \clist_count:n { #1 } }
		\clist_map_inline:nn { #1 }
		{
			\refstepcounter{equation}\label{##1}
		}
		\addtocounter{equation}{-\g_mixedsubeq_int}
		\setcounter{mixedsubequations}{1}
	}
	\domixedsubequations
}
{\ignorespacesafterend}

\NewDocumentCommand{\domixedsubequations}{}
{
	\cs_set:Npx \theequation
	{
		\exp_not:o { \theequation }
		\exp_not:n { \alph{mixedsubequations} }
	}
	\ignorespaces
}
\ExplSyntaxOff

\makeatletter
\let\save@mathaccent\mathaccent
\newcommand*\if@single[3]{%
  \setbox0\hbox{${\mathaccent"0362{#1}}^H$}%
  \setbox2\hbox{${\mathaccent"0362{\kern0pt#1}}^H$}%
  \ifdim\ht0=\ht2 #3\else #2\fi
  }
\newcommand*\rel@kern[1]{\kern#1\dimexpr\macc@kerna}
\newcommand*\widebar[1]{\@ifnextchar^{{\wide@bar{#1}{0}}}{\wide@bar{#1}{1}}}
\newcommand*\wide@bar[2]{\if@single{#1}{\wide@bar@{#1}{#2}{1}}{\wide@bar@{#1}{#2}{2}}}
\newcommand*\wide@bar@[3]{%
  \begingroup
  \def\mathaccent##1##2{%
    \let\mathaccent\save@mathaccent
    \if#32 \let\macc@nucleus\first@char \fi
    \setbox\z@\hbox{$\macc@style{\macc@nucleus}_{}$}%
    \setbox\tw@\hbox{$\macc@style{\macc@nucleus}{}_{}$}%
    \dimen@\wd\tw@
    \advance\dimen@-\wd\z@
    \divide\dimen@ 3
    \@tempdima\wd\tw@
    \advance\@tempdima-\scriptspace
    \divide\@tempdima 10
    \advance\dimen@-\@tempdima
    \ifdim\dimen@>\z@ \dimen@0pt\fi
    \rel@kern{0.6}\kern-\dimen@
    \if#31
      \overline{\rel@kern{-0.6}\kern\dimen@\macc@nucleus\rel@kern{0.4}\kern\dimen@}%
      \advance\dimen@0.4\dimexpr\macc@kerna
      \let\final@kern#2%
      \ifdim\dimen@<\z@ \let\final@kern1\fi
      \if\final@kern1 \kern-\dimen@\fi
    \else
      \overline{\rel@kern{-0.6}\kern\dimen@#1}%
    \fi
  }%
  \macc@depth\@ne
  \let\math@bgroup\@empty \let\math@egroup\macc@set@skewchar
  \mathsurround\z@ \frozen@everymath{\mathgroup\macc@group\relax}%
  \macc@set@skewchar\relax
  \let\mathaccentV\macc@nested@a
  \if#31
    \macc@nested@a\relax111{#1}%
  \else
    \def\gobble@till@marker##1\endmarker{}%
    \futurelet\first@char\gobble@till@marker#1\endmarker
    \ifcat\noexpand\first@char A\else
      \def\first@char{}%
    \fi
    \macc@nested@a\relax111{\first@char}%
  \fi
  \endgroup
}
\makeatother

\DeclareMathOperator{\Cay}{Cay}

\numberwithin{equation}{section}

\title{\sffamily Geometry of Random Cayley Graphs of Abelian Groups}

\author{\sffamily Jonathan Hermon\quad Sam Olesker-Taylor}

\date{}

\usepackage{geometry}
\begin{document}

\maketitle

\acknofootnote

\vspace{-6ex}

\renewcommand{\abstractname}{\sffamily Abstract}
\begin{abstract}
Consider the random Cayley graph of a finite Abelian group $G$ with respect to $k$ generators chosen uniformly at random, with $1 \ll \log k \ll \log |G|$. Draw a vertex $U \sim \operatorname{Unif}(G)$.

We show that the graph distance $\operatorname{dist}(\mathsf{id},U)$ from the identity to $U$ concentrates at a particular value $M$, which is the minimal radius of a ball in $\mathbb Z^k$ of cardinality at least $|G|$, under mild conditions. In other words, the distance from the identity for all but $o(|G|)$ of the elements of $G$ lies in the interval $[M - o(M), M + o(M)]$. In the regime $k \gtrsim \log \abs G$, we show that the diameter of the graph is also asymptotically $M$. In the spirit of a conjecture of Aldous and Diaconis~\cite{AD:conjecture}, this $M$ depends only on $k$ and $|G|$, not on the algebraic structure of $G$.

Write $d(G)$ for the minimal size of a generating subset of $G$. We prove that the order of the spectral gap is $|G|^{-2/k}$ when $k - d(G) \asymp k$ and $|G|$ lies in a density-$1$ subset of $\mathbb N$ or when $k - 2 d(G) \asymp k$. This extends, for Abelian groups, a celebrated result of Alon and Roichman~\cite{AR:cayley-expanders}.

The aforementioned results all hold with high probability over the random Cayley graph.
\end{abstract}

%
%
%
%

\small
\begin{quote}
\begin{description}
	\item [Keywords:]
	typical distance, diameter, spectral gap, relaxation time, random Cayley graphs
	
	\item [MSC 2020 subject classifications:]
	05C12, 05C48, 05C80; 60B15, 60K37
\end{description}
\end{quote}
\normalsize


\vspace{2ex}
\numberingroman

\vfill
\printtoc{\value{tocdepth}}
\vspace*{2ex}

\newpage
\section{Introduction and Statement of Results}
\label{sec-p3:intro}

\subsection{Brief Overview of Results and Notation}

\subsubsection{Brief Overview of Results}

We analyse geometric properties of a \textit{Cayley graph} of a finite group; the focus is on Abelian groups.
The generators of this graph are chosen independently and uniformly at random.
Precise definitions are given in \S\ref{sec-p3:intro:cayley-def}. For now, let $G$ be a finite group, let $k$ be an integer (allowed to depend on $G$) and denote by $G_k$ the Cayley graph of $G$ with respect to $k$ independently and uniformly random generators.
We consider values of $k$ with $1 \ll \log k \ll \log \abs G$ for which $G_k$ is connected with high probability (abbreviated \textit{whp}), ie with probability tending to 1 as $\abs G$ grows.
For an Abelian group $G$,
write $d(G)$ for the minimal size of a generating subset of $G$.

\begin{itemize}[itemsep = 0pt, topsep = \smallskipamount, label = \bcdot]
	\item 
	\textit{Typical Distance.}
	Draw $U \sim \Unif(G)$.
	We show that the law of the graph distance between the identity and $U$ concentrates.
	The leading order term in this typical distance depends only on $k$ and $\abs G$ when
		$1 \ll k \ll \log \abs G / \log\log\log \abs G$ and $k - d(G) \asymp k$
	or
		$k \gg \log \abs G$.
	
	\item 
	\textit{Diameter.}
	For
		$k \asymp \log \abs G$
		under mild conditions on the group
	and
		$k \gg \log \abs G$
		for any Abelian group,
	we show that the diameter concentrates at the same value as the typical distance.
	
	\item 
	\textit{Spectral Gap.}
	For any $1 \ll k \lesssim \log \abs G$ with $k - d(G) \asymp k$,
	we determine the order of the spectral gap of the random walk on the random Cayley graph.
\end{itemize}

Introduced by \textcite{AD:conjecture}, there has been a great deal of research into these random Cayley graphs.
Motivation for this model and an overview of historical work is given~in~\S\ref{sec-p3:intro:previous-work}.

\subsubsection{Notation and Terminology}

Cayley graphs are either directed or undirected; we emphasise this by writing $G_k^+$ and $G_k^-$, respectively.
When we write $G_k$ or $G^\pm_k$, this means ``either $G^-_k$ or $G^+_k$'', corresponding to the undirected, respectively directed, graphs with generators chosen independently and uniformly at~random.

Conditional on being simple, $G^+_k$ is uniformly distributed over the set of all simple degree-$k$ Cayley graphs. Up to a slightly adjusted definition of \textit{simple} for undirected Cayley graphs, our results hold with $G_k$ replaced by a uniformly chosen simple Cayley graph of degree $k$; see \S\ref{sec-p3:intro:rmks:typ-simp}.

Our results are for sequences $(G_N)_\Ninn$ of finite groups with \toinf{\abs{G_N}} \asinf N.
For ease of presentation, we write statements like ``let $G$ be a group'' instead of ``let $(G_N)_\Ninn$ be a sequence of groups''.
Likewise, the quantities $d(G)$ and, of course, $k$ appearing in the statements all correspond to sequences, which need not be fixed (or bounded) unless we explicitly say otherwise.
In the same vein, an event holds \textit{with high probability} (abbreviated \textit{\whp}) if its probability tends~to~1.

We use standard asymptotic notation:
	``$\ll$'' or ``$\oh{\cdot}$'' means ``of smaller order'';
	``$\lesssim$'' or $\Oh{\cdot}$'' means ``of order at most'';
	``$\asymp$'' means ``of the same order'';
	``$\eqsim$'' means ``asymptotically equivalent''.

\subsection{Statements of Main Results}
\label{sec-p3:intro:res}

\subsubsection{Typical Distance for Abelian Groups}
\label{sec-p3:intro:typdist}

\nextresult

\begin{subtheorem-num}{intrormkT}

Our first result concerns typical distance in the random Cayley graph.

\begin{introdefn}
	For
		a group $G$,
		$k \in \mbn$
	and
		$\beta \in (0,1)$,
	define the \textit{$\beta$-typical distance} $\mcd_{G_k}(\beta)$ via
	\[
		\mcb_{G_k}(R) \cq \bra{ x \in G \midb \dist_{G_k}(\id,x) \le R }
	\Qand
		\mcd_{G_k}(\beta) \cq \min\brb{ R \ge 0 \midb \abs{ \mcb_{G_k}(R) } \ge \beta \abs G }.
	\]
\end{introdefn}

Informally,
we show that the mass (in terms of number of vertices) concentrates at a thin `slice', or `shell', consisting of vertices at a distance
$M \pm \oh M$ from the origin, with $M$ explicit.

Investigating this \emph{typical distance} for $G_k$ when $k$ diverges with $\abs G$ was suggested to us by \textcite{B:typdist:private}.
Previous work concentrated on fixed $k$, ie independent of $\abs G$; see \S\ref{sec-p3:intro:previous-work}.

\smallskip

For an Abelian group $G$, write $d(G)$ for the minimal size of a generating subset of $G$ and
\[
	m_*(G)
\cq
	\max\brb{ \mint{j \in [d]} m_j \midb \oplus_{j=1}^d \: \mbz_{m_j} \text{ is a decomposition of } G }.
\]
The condition $1 \ll \log k \ll \log \abs G$ is necessary for the type of concentration of measure we show in \cref{res-p3:intro:typdist}; see \cref{rmk-p3:intro:typdist:k-nec}.
Refinements of \cref{res-p3:intro:typdist} are given in \cref{res-p3:typdist1:res,res-p3:typdist2:res,res-p3:typdist3:res}. 

\begin{introthm}
\label{res-p3:intro:typdist}
	Let $G$ be an Abelian group.
	The following convergences are in probability \asinf{\abs G}.
	
	\begin{itemize}[itemsep = 0pt, topsep = \smallskipamount, label = \bcdot]
		\item 
		Consider $1 \ll k \ll \log \abs G$;
		suppose $k - d(G) \asymp k$ and $d(G) \ll \log \abs G / \log \log \abs G$.
		Write $\mfd^+ \cq \abs{G}^{1/k} / (2e)$ and $\mfd^- \cq \abs{G}^{1/k} / e$.
		For all $\beta \in (0,1)$,
		we have
		\(
			\mcd_{G_k^\pm}(\beta) / \mfd^\pm
		\to^\mbp
			1.
		\)
		
		\item 
		Consider $k \eqsim \lambda \log \abs G$ with $\lambda \in (0,\infty)$;
		suppose $d(G) \le \tfrac14 \log \abs G / \log\log \abs G$ and $m_*(G) \gg 1$.
		There exists a constant $\alpha^\pm_\lambda \in (0,\infty)$ so that,
		for all $\beta \in (0,1)$,
		we have
		\(
			\mcd_{G_k^\pm}(\beta) / \rbr{ \alpha^\pm_\lambda k }
		\to^\mbp
			1.
		\)
		
		\item 
		Consider $k \gg \log \abs G$ with $\log k \ll \log \abs G$;
		write $\rho \cq \log k / \log \log \abs G$ so that $k = (\log \abs G)^\rho$.
		(We allow $\rho \gg 1$.)
		For all $\beta \in (0,1)$, we have
		\(
			\mcd_{G_k^\pm}(\beta) / \rbb{ \tfrac{\rho}{\rho-1} \log_k \abs G }
		\to^\mbp
			1.
		\)
	\end{itemize}
	
	The implicit lower bound holds for all Abelian groups and all choices of $k$ generators.
\end{introthm}

\begin{intrormkt}
We establish the concentration of typical distance via three distinct approaches, in
	\S\ref{sec-p3:typdist1}, \S\ref{sec-p3:typdist2} and \S\ref{sec-p3:typdist3}.
Conceptually, all involve sizes of lattice balls and drawing elements uniformly from balls.
A precise statement for each approach is given, as is an outline of the proof.
In summary, \cref{res-p3:intro:typdist} is a direct consequence of
	\cref{res-p3:typdist1:res,res-p3:typdist2:res,res-p3:typdist3:res};
	see also
	\cref{hyp-p3:typdist1,hyp-p3:typdist2,hyp-p3:typdist3}.
\end{intrormkt}

\begin{intrormkt}
	For smaller $k$, namely $1 \ll k \ll \sqrt{\log \abs G / \log \log \log \abs G}$,
	we can relax $k - d(G) \asymp k$ to $k - d(G) \gg 1$.
	In order to generate the group, we certainly need $k \ge d(G)$, by definition.
	In many cases $k - d(G) \gg 1$ is necessary in order to generate the group whp, so this assumption can not be removed. For a characterisation of these cases and related discussion, see~%
	\cite[Lemma~8.1]{HOt:rcg:abe:extra}.
\end{intrormkt}

\begin{intrormkt}
\label{rmk-p3:intro:typdist:mixing-geom}
	Interesting is how we prove this theorem.
	It is common in mixing time proofs to use geometric properties of the graph, such as expansion or distance properties.
	We do the opposite: we use mixing techniques to prove this geometric result.
	This is in the same spirit as \cite{LP:ramanujan}; see~\S\ref{sec-p3:intro:previous-work}.
\end{intrormkt}


\begin{intrormkt}
\label{rmk-p3:intro:typdist:k-nec}
	We discuss briefly lack of concentration of measure if the condition $1 \ll \log k \ll \log \abs G$ fails.
	The method developed here can be applied to the regime $\log k \asymp \log \abs G$, for any
	$G$:
	it gives
	\(
		\abs{ \mcd_{G_k}(\beta) - m }
	\le
		1
	\)
	\whp
	with $m \cq \ceil{\log_k \abs G} \asymp 1$,
	for all $\beta \in (0,1)$.
	The values $m-1$, $m$ and $m+1$ differ by a constant factor.
	Hence, there is no concentration in the sense of \cref{res-p3:intro:typdist}.
	
	Results of \textcite{SZ:diam-cayley,MS:diameter-cayley-cycle} show that there exist Abelian $G$ and $1 \asymp k \ge d(G)$ such that there is no concentration; see \S\ref{sec-p3:intro:previous-work:typdist} for more details.
	\cite{SZ:diam-cayley,MS:diameter-cayley-cycle} analyse the asymptotic laws of the diameter, but the deduction holds for typical distance~too.
\end{intrormkt}

\end{subtheorem-num}

\subsubsection{Typical Distance for Nilpotent Groups}
\label{sec-p3:intro:nil}

We can extend \cref{res-p3:intro:typdist} to the set-up of nilpotent, rather than Abelian, groups $G$.

\nextresult

\begin{subtheorem-num}{introdefn}
	\label{def-p3:intro:nil}

\begin{introdefn}
	A group $G$ is \textit{nilpotent} if its \textit{lower central series}, ie the sequence $(G_\ell)_{\ell \ge 0}$ defined by $G_0 \cq G$ and $G_\ell \cq [G_{\ell-1}, G]$ for $\ell \ge 1$, stabilises at the trivial group.
	
	The \textit{step} $\ell(G)$ is the number of terms until stabilisation:
	\(
		\ell(G)
	\cq
		\inf\bra{ \ell \ge 0 \mid G_\ell = \bra{\id} }.
	\)
	The \textit{rank} $d(G)$ is the minimal size of a generating subset.
	It is a standard fact that
	a symmetric set $S \subseteq G$ generates $G$ if and only if $S^\ab \cq \bra{ [G, G] s \mid s \in S }$ generates $G^\ab$.
	Hence, $d(G) = d(G^\ab)$.
\end{introdefn}

\begin{introdefn}
	Let $G$ be a group and let $S$ be a symmetric mutlisubset of $G$.
	Let $\Cay(G, S)$ be the (right) Cayley graph of $G$ wrt $S$.
	Let $H \trianglerighteq G$.
	Let $\diam_S(H) \cq \max\bra{ d_S(\id, h) \mid h \in H }$ denote the diameter of $H$ wrt the graph distance $d_S(\cdot, \cdot\cdot)$ on $\Cay(G, S)$.
	
	Let $\diam_S(G/H)$ denote the diameter of $\Cay(G/H, HS)$ where $HS \cq \bra{H s \mid s \in S}$.
	Similarly, for $H \trianglerighteq H' \trianglerighteq G$, let $\diam_S(H' / H)$ denote the diameter of $\Cay(H'/H, HS)$; this is a slight abuse of notation, as $\diam_{\bra{ Hs \mid s \in S }}(H'/H)$ would be consistent with the definition of $\diam_S(H)$.
\end{introdefn}

\begin{introdefn}
	For $\beta \in (0,1]$, let $\mcd_S(G,\beta)$ be the minimal $r \in \mbn$ such that a ball of radius $r$ in $\Cay(G,S)$ contains at least $\beta |G|$ elements of $G$.
	Write $\mcd_S(G^\ab,\beta) \cq \mcd_{S^\ab}(G^\ab,\beta)$, ie for the minimal $r \in \mbn$ such that a ball of radius $r$ in $\Cay(G^\ab,S^\ab)$ contains at least $\beta \abs{G^\ab}$ elements~of~$G^\ab$.
\end{introdefn}

\end{subtheorem-num}

The following theorem bounds the typical distance $\mcd_S(G,\beta)$ for $\beta \in (0,1]$---and, in particular, the diameter by taking $\beta \cq 1$---for a nilpotent group $G$ and a symmetric set of generators $S$ in terms of $\mcd_S(G^\ab,\beta)$, ie the corresponding typical distance in the Abelianisation.
%

\begin{introthm}
\label{res-p3:intro:nil}
	Let $G$ be a finite nilpotent group of step $\ell$ and rank $d$.
	Let $k \in \mbn$ and let $s_1, ..., s_k \in G$.
	Let $S \cq [s_1, s_1^{-1}, ..., s_k, s_k^{-1}]$ be a symmetric multisubset of $G$.
	Then, for all $\beta \in (0, 1]$,
	we have
	\[
		0
	\le
		\mcd_S(G,\beta) - \mcd_S(G^\ab,\beta)
	\lesssim
		\diam_S(G^\ab)^{3/4}
	\le
		\rbb{ 3 \mcd_S(G^\ab,\beta)/\beta }^{3/4}.
	\]
	if
	\(
		d
	\le
		k
	\le
		\tfrac1{16} \ell^{-1} d^{-\ell} \log \abs G / \log \log \abs G.
	\)
	Further, for all $\beta \in (0,1/2)$, we have
	\[
		\mcd_S(G^\ab,1-\beta) - \mcd_S(G^\ab,\beta)
	\le
		2 \sqrt{\beta^{-1}\trel\rbb{ \Cay(G^\ab,S^\ab) }},
	\]
	where $\trel(\Cay(G^\ab,S^\ab))$ is the relaxation time of the simple random walk on $\Cay(G^\ab,S^\ab)$.
\end{introthm}


\cref{res-p3:intro:nil} is a corollary of a more general result to appear in an upcoming paper \cite{H+:nil-to-abe} by one of the authors.
We state it as \cref{res-p3:nil:thm} and explain how to deduce the first part of \cref{res-p3:intro:nil} from it in \S\ref{sec-p3:nil:thm}.
The second part of \cref{res-p3:intro:nil} is a standard Poincaré-type inequality using the fact that graph distance is trivially $1$-Lipschitz; this deduction was pointed out to us by \textcite{S:concentration:private}.

The argument for \cref{res-p3:nil:thm} builds on one due to \textcite{EbP:cayley-diam-nil} who, in turn, relied on ideas of \textcite{BT:moderate-growth}.
\citeauthor{EbP:cayley-diam-nil} considered
$\max\bra{\abs S, \ell, d} \asymp 1$.
We are primarily interested in $k \gg 1$, so we must keep track of the certain dependences on $|S|$ and $\ell$.

For now, we explore some simple consequences of \cref{res-p3:intro:nil}.

\begin{intrormkt}
\label{rmk-p3:intro:nil}
	%
Whenever $\trel\rbr{ \Cay(G^\ab, S^\ab) } \ll \rbr{ \diam_S(G^\ab) }^2$,
for all fixed $\beta \in (0,1/2)$,
we obtain
\[
	\mcd_S(G^\ab,1-\beta) - \mcd_S(G^\ab,\beta)
\ll
	\mcd_S(G^\ab,1/2).
\] 
If, in addition,
\(
	d
\le
	k
\le
	\tfrac1{16} \ell^{-1} d^{-\ell} \log \abs G / \log \log \abs G,
\)
in the notation of \cref{res-p3:intro:nil},
then
\[
	\mcd_S(G,1-\beta) - \mcd_S(G,\beta)
\ll
	\mcd_S(G^\ab,1/2).
\]

By \cref{res-p3:intro:typdist} above,
\(
	\diam_S(G^\ab)
\gtrsim
	k \abs{G^\ab}^{1/k}
\)%
---deterministically, in fact.
Hence,
\[
	\trel\rbb{ \Cay(G^\ab,S^\ab) }
\ll
	k^2 \abs{G^\ab}^{2/k}
\Quad{implies}
	\trel\rbb{ \Cay(G^\ab,S^\ab) }
\ll
	\rbb{ \diam_S(G^\ab) }^2.
\]
\cref{res-p3:intro:gap} below roughly gives $\trel(H) \asymp \abs H^{2/k}$ for Abelian $H$.
Precisely, it implies that
the former relation above
holds whp when $d(G) \asymp 1$, $k \gg 1$ and $S = Z \cup Z^{-1}$ with $Z_1, ..., Z_k \sim^\iid \Unif(G)$.
\end{intrormkt}

\cref{res-p3:intro:nil} and its consequences laid out in \cref{rmk-p3:intro:nil} describe, in a formal manner, some of the key differences between geometry in Abelian versus non-Abelian groups.
We lay out a more informal explanation in \S\ref{sec-p3:nil:intuition}.
We go through the Abelian proof, pointing out where it fails.

We discuss what adjustment to the approach is needed to establish a result analogous to \cref{res-p3:intro:typdist}.
In particular, we reference our companion paper \cite{HOt:rcg:matrix} in which we establish typical distance for a particular matrix group, describing the main changes in that set-up.

\subsubsection{Diameter}
\label{sec-p3:intro:diam}

\nextresult

We can extend our proof to consider the \textit{diameter}, ie the maximal distances between pairs of vertices in the graph, in the regime $k \gtrsim \log \abs G$.
For a graph $H$, denote by $\diam H$ its diameter.

Our first diameter result gives concentration for $k \gtrsim \log \abs G$ and $\log k \ll \log \abs G$.
Again, $\log k \ll \log \abs G$ is a necessary condition; see \cref{rmk-p3:intro:diam:conc:k-nec}.
A refinement of \cref{res-p3:intro:diam:conc} is given in \cref{res-p3:diam:conc:res}.

\begin{introthm}
\label{res-p3:intro:diam:conc}
	Let $G$ be an Abelian group.
	The following convergences are in probability \asinf{\abs G}.
	
	\begin{itemize}[itemsep = 0pt, topsep = \smallskipamount, label = \bcdot] 
		\item 
		Consider $k \eqsim \lambda \log \abs G$ with $\lambda \in (0,\infty)$;
		suppose $d(G) \le \tfrac14 \log \abs G / \log\log \abs G$ and $m_*(G) \gg 1$.
		Let $\alpha^\pm_\lambda \in (0,\infty)$ be the constant from \cref{res-p3:intro:typdist}.
		We have
		\(
			\diam G_k^\pm / \rbr{ \alpha^\pm_\lambda k }
		\to^\mbp
			1.
		\)
		
		\item 
		Consider $k \gg \log \abs G$ with $\log k \ll \log \abs G$;
		write $\rho \cq \log k / \log\log \abs G$ so that $k = (\log \abs G)^\rho$.
		(We allow $\rho \gg 1$.)
		We have
		\(
			\diam G_k^\pm / \rbr{ \tfrac{\rho}{\rho-1} \log_k \abs G }
		\to^\mbp
			1.
		\)
		The upper bound holds~for~all~$G$.
	\end{itemize}
	
	The implicit lower bound holds for all Abelian groups and all choices of $k$ generators.
\end{introthm}

\begin{subtheorem-num}{intrormkT}
	\label{rmk-p3:intro:diam:conc}

\begin{intrormkt}
	For any Cayley graph $H$ one has $\mcd_H(\tfrac12) \le \diam H \le 2 \, \mcd_H(\tfrac12) + 1$.
	Indeed, $(x_1, ..., x_\ell)$ is a path in the Cayley graph if and only if $(x_\ell, ..., x_1)$ is a path when all generators inverted.
	Hence, the typical distance and diameter are always equivalent up to constants.
	\cref{res-p3:intro:diam:conc} gives conditions under which they are asymptotically equivalent whp for random Cayley graphs.
	
	We establish cutoff for the simple random walk on $G_k$ for many Abelian groups when $1 \ll k \lesssim \log \abs G$ in
	\cite[Theorem~A]{HOt:rcg:abe:cutoff}.
	Combined with \cref{res-p3:intro:typdist}, it shows that $\tmix(G_k) \asymp (\diam G_k)^2/k$ whp when $k - d(G) \asymp k \gtrsim \log \abs G$.
	One can also consider non-Abelian groups; see
	\cite[Theorem~E]{HOt:rcg:matrix}.
\end{intrormkt}

\begin{intrormkt}
\label{rmk-p3:intro:diam:conc:k-nec}
	Analogously to \cref{rmk-p3:intro:typdist:k-nec},
	we have
	\(
		\abs{ \diam G_k^\pm - \ceil{\log_k \abs G} }
	\le
		1
	\)
	\whp
	when $\log k \asymp \log \abs G$.
	Again, \cite{MS:diameter-cayley-cycle,SZ:diam-cayley} imply that there is no concentration for some Abelian $G$ when $k \asymp 1$.
\end{intrormkt}

\end{subtheorem-num}

\nextresult

Our next diameter result shows, in a well-defined sense, that, amongst all groups, when $k - \log_2 \abs G \asymp k$ with $\log k \ll \log \abs G$, the group $\mbz_2^d$ gives rise to the largest typical diameter.

\begin{introdefn*}
	For two random sequences $\alpha \cq (\alpha_N)_\Ninn$ and $\beta \cq (\beta_N)_\Ninn$ of reals,
	we say that \textit{$\alpha \le \beta$ whp up to smaller order terms} if there exist non-random sequences $(\gamma_N)_\Ninn$ and $(\delta_N)_\Ninn$ of reals with $\delta_N \to 0$ \asinf N such that
	\(
		\rbr{ \bra{ \alpha_N \le (1 + \delta_N) \gamma_N } }_\Ninn
	\)
	and
	\(
		\rbr{ \bra{ (1 - \delta_N) \gamma_N \le \beta_N } }_\Ninn
	\)
	both~hold~whp.
	We say that \textit{$\alpha \eqsim \beta$ whp} if $\alpha \le \beta$ and $\beta \le \alpha$ whp up to smaller order~terms.
\end{introdefn*}

We now define the candidate radius which we show is an upper bound for $\diam G_k$ whp.

\begin{introdefn}
	Write $\mfr(k, n)$ for the minimal $R \in \mbn$ with $\binomt kR \ge n$.
\end{introdefn}

We now state our second diameter result.
A refinement of \cref{res-p3:intro:diam:univ} is given in \cref{res-p3:diam:univ:res}.

\begin{introthm}
\label{res-p3:intro:diam:univ}
	Let $G$ be an arbitrary group.
	Suppose that $k - \log_2 \abs G \asymp k$ and $1 \ll \log k \ll \log \abs G$.
	Then
	\(
		\diam G_k
	\le
		\mfr(k, \abs G)
	\)
	up to smaller order terms \whp;
	further,
	if $H \cq \mbz_2^d$, then
	\(
		\diam H_k
	\eqsim
		\mfr(k, 2^d)
	=
		\mfr(k, \abs H)
	\)
	\whp.
	(The limit is as the size of the group diverges.)
\end{introthm}

This gives a quantitative sense in which $\mbz_2^d$ is the group giving rise to the largest diameter.

\begin{introcor}
	For
		all diverging $d$ and $n$ with $n \le 2^d$
	and
		all groups $G$ of size $n$,
	if $k - \log_2 n \asymp k$ and $\log k \ll \log n$, then
	\(
		\diam G_k
	\le 
		\diam H_k
	\)
	where $H \cq \mbz_2^d$
	up to smaller order terms \whp.
\end{introcor}


\textcite[Conjecture~7]{W:rws-hypercube} conjectures an analogous statement for mixing times.
We prove an extension of this conjecture in
\cite[Theorems~C and~D]{HOt:rcg:abe:cutoff}
when restricted to nilpotent groups.

\subsubsection{Spectral Gap}
\label{sec-p3:intro:gap}

\nextresult

Our next result concerns the spectral gap and relaxation time of the random Cayley graph.

\begin{introdefn}
	Consider a reversible Markov chain with (real) eigenvalues $1 = \lambda_1 \ge \lambda_2 \ge \cdots \ge \lambda_n \ge -1$ of its transition matrix.
	The \textit{usual}, respectively \textit{absolute}, \textit{spectral gap} is defined as
	\[
		\gamma
	\cq
		\MIN{i \ne 1} \bra{1 - \lambda_i}
	=
		1 - \lambda_2,
	\Quad{respectively}
		\gamma_*
	\cq
		\MIN{i \ne 1} \bra{1 - \abs{\lambda_i}}
	=
		1 - \max\bra{\abs{\lambda_2}, \abs{\lambda_n}};
	\]
	the \textit{usual}, respectively \textit{absolute}, \textit{relaxation time} is defined as
	\(
		\trel \cq 1/\gamma,
	\Quad{respectively,}
		\trel^* \cq 1/\gamma_*.
	\)
	
	The \textit{spectral gap} or \textit{relaxation time of a graph}, is that of the simple random walk on the graph.
\end{introdefn}

It is classical that under reversibility in continuous-time the spectral gap asymptotically determines the exponential rate of convergence to equilibrium, whereas in discrete-time it is determined by the absolute spectral gap; see \cite[\S 12 and \S 20]{LPW:markov-mixing}.
%
For a multiset $z = [z_1, ..., z_k]$ with $z_1, ..., z_k \in G$, we write $G^-(z)$ for the undirected Cayley graph with respect to the generators $z_1, ..., z_k$.

A refinement of \cref{res-p3:intro:gap} is given in \cref{res-p3:gap:res}.
We do not require $\kinf$ as $\abs G \to \infty$.

\begin{introthm}
\label{res-p3:intro:gap}
	There exists a positive constant $c$ so that,
	for
		all Abelian groups $G$,
		all $k$
	and
		all multisets of generators $z$ of size $k$,
	we~have
	\[
		\trel^*\rbb{ G^-(z) }
	\ge
		\trel\rbb{ G^-(z) }
	\ge
		c \abs G^{2/k}.
	\]
	For all $\delta > 0$,
	there exists a constant $C_\delta > 0$ so that,
	for all Abelian groups $G$,
	if $k \ge (2 + \delta) d(G)$,~then
	\[
		\pr{
			\trel^*\rbr{ G^-_k }
		\le
			C_\delta \abs G^{2/k}
		}
	\ge
		1 - C_\delta 2^{-k/C_\delta}.
	\]
	Further, for all $\eps \in (0,1)$, there exists a density-$(1 - \eps)$ subset $\mba \subseteq \mbn$ so that if $\abs G \in \mba$ then the condition $k \ge (2 + \delta) d(G)$ can be relaxed to $k \ge (1 + \delta) d(G)$; the constants now also depend on $\eps$.
\end{introthm}


The inequality $\trel\rbr{ G^-(z) } \gtrsim \abs G^{2/k}$ was previously proved by \textcite[Theorem~1.1]{H:cutoff-cayley-<}.
That theorem is only stated for $G = \mbz_p$ with $p$ prime and $k \le \log p / \log \log p$, which is the set-up of the rest of the paper. It appears that this lower bound holds in greater generality, possibly with a different absolute constant.
The proof involves an elegant use of Minkowski's theorem.
We present our own proof of the lower bound; it is the easier, and arguably less interesting, part of \cref{res-p3:intro:gap}.

The method of proof for this result is rather different to our previous results and also somewhat different to those used by others to study the spectral gap of random Cayley graphs; see \S\ref{sec-p3:intro:previous-work:gap}.

\subsection{Historic Overview}
\label{sec-p3:intro:previous-work}

In this subsection, we give a fairly comprehensive account of previous work on distance metrics on and spectral gap of random Cayley graphs; we compare our results with existing ones.
We also mention, where relevant, other results which we have proved in companion papers;
see also \S\ref{sec-p3:intro:rmks:advert}.

\subsubsection{Motivation: Random Cayley Graphs and Cutoff for Random Walks}
\label{sec-p3:intro:previous-work:motivation}

In their seminal paper, \textcite{AD:conjecture,AD:shuff-stop} considered random walks on \emph{random} Cayley graphs.
\textcite{D:cayley:private} gave the following (paraphrased) motivation.
\begin{quote} \slshape
	Erd\H{o}s, when considering classes of mathematical objects, often combinatorial or graph theoretic, would often ask,
	``What does a typical object in this class `look like'?''
	If an object is chosen uniformly at random, are there natural properties which hold whp?
	
	It is then natural to ask,
	``How does a typical random walk on a group behave?''
\end{quote}
This lead \textcite{AD:conjecture,AD:shuff-stop} to consider the set of all Cayley graphs of a given group $G$ with $k$ generators.
Drawing such a Cayley graph uniformly at random corresponds to choosing generators $Z_1, ..., Z_k \sim^\iid \Unif(G)$, conditional on giving rise to a simple graph; see \S\ref{sec-p3:intro:rmks:typ-simp}.
We study random walks in \cite{HOt:rcg:abe:cutoff,HOt:rcg:matrix}, establishing cutoff and showing universal mixing bounds in different set-ups.

\subsubsection{Universality: The Aldous--Diaconis Conjecture}
\label{sec-p3:intro:previous-work:ad-conj}

\textcite{AD:conjecture,AD:shuff-stop} made the following (informal) conjecture:
	regardless of the particular group $G$,
	provided $k \gg \log \abs G$,
	the random walk on the random Cayley graph exhibits cutoff whp at a time which depends only on $k$ and $\abs G$.
This was established for Abelian groups by Dou and Hildebrand \cite{DH:enumeration-rws,H:cutoff-cayley->};
in \cite{HOt:rcg:matrix}, we provide a counterexample using unit upper triangular matrix groups.
For more details, see our companion articles \cite{HOt:rcg:abe:cutoff,HOt:rcg:matrix} where we study cutoff extensively.

The point of the Aldous--Diaconis conjecture is that certain statistics should be ``independent of the algebraic structure of the group'', ie only depend on $G$ through $\abs G$.
The current article shows how very related statements to those above hold when ``cutoff'' is replaced by ``typical distance''.
Namely, we give conditions under which the typical distances concentrates on a value that depends only on $k$ and $\abs G$; see \cref{res-p3:typdist1:res,res-p3:typdist2:res,res-p3:typdist3:res}. 

\subsubsection{Typical Distance and Diameter}
\label{sec-p3:intro:previous-work:typdist}

Previous work on distance metrics (detailed below) had concentrated on the case where the number of generators $k$ is a \emph{fixed} number.
The results establish (non-degenerate)
limiting laws.
This restricts the (sequences of) groups which can be studied; eg, in order for it to be even possible to generate the group---never mind having independent, uniform generators do so whp---one needs $d(G) \le k \asymp 1$.
We discuss generation of groups further in
\cite[\S 8]{HOt:rcg:abe:extra};
see in particular
\cite[\S 8.2]{HOt:rcg:abe:extra}
where we describe adaptations made in order to obtain connected graphs in the references given below.

Our results are in a different direction:
	for us, \toinf k \asinf{\abs G} and we establish concentration of the observables.
This allows us to consider a much wider range of groups, in particular with $d(G)$ diverging with $\abs G$.
This line of enquiry was suggested to us by \textcite{B:typdist:private}

\textcite{AGg:diam-cayley-Zq} studied the diameter of the random Cayley graph of cyclic groups of prime order.
They prove (for fixed $k$) that the diameter is order $\abs G^{1/k}$; see \cite[Theorems~1 and 2]{AGg:diam-cayley-Zq}.
They conjecture that the diameter divided by $\abs G^{1/k}$ converges in distribution to some non-degenerate distribution \asinf{\abs G}; see \cite[Conjecture~3]{AGg:diam-cayley-Zq}.

\textcite{MS:diameter-cayley-cycle} consider, as a consequence of a quite general framework, the diameter of the random Cayley graph of $\mbz_n$ with respect to a fixed number $k$ of random generators, for a random $n$, without any primality assumption.
They derive distributional limits for the diameter, the average distance (defined with respect to various $L_p$ metrics) and the girth.
They determine limit distributions for each of these, and in some cases derive explicit formulas.

\textcite{SZ:diam-cayley} build on the framework of \textcite{MS:diameter-cayley-cycle}, again for fixed $k$; they are able to consider non-random $n$, as well as Abelian groups of arbitrary (fixed) rank, instead of only cyclic groups.
In particular, they verify the conjecture of \textcite[Conjecture~3]{AGg:diam-cayley-Zq};
they additionally work with average distance and girth.

\smallskip

\textcite{LP:ramanujan} derive an analogous typical distance result for $n$-vertex, $d$-regular Ramanujan graphs:
	whp all but $\oh n$ of the vertices lie at a distance $\log_{d-1} n \pm \Oh{\log\log n}$;
	they establish this by proving cutoff for the non-backtracking random walk at time $\log_{d-1} n$.

Related work on the diameter of random Cayley graphs, including concentration of certain measures, can be found in \cite{LPSW:cayley-digraphs,S:diam-cayley}.

\smallskip

The Aldous--Diaconis conjecture for mixing can be transferred naturally to typical distance:
	the mass should concentrate at a distance $M$, where $M$ can be written as a function only of $k$ and $\abs G$;
	ie there is concentration of mass at a distance independent of the algebraic structure of~the~group.

In
\cite[Theorem~E]{HOt:rcg:matrix}
we consider typical distance analogously to this paper; there the underlying group is a non-Abelian matrix group.
In contrast with the Abelian groups in \cref{res-p3:intro:typdist}, the $M$ for these non-Abelian groups cannot be written as a function only of $k$ and $\abs G$.

\subsubsection{Spectral Gap}
\label{sec-p3:intro:previous-work:gap}

\textcite[Theorem~1.1]{H:cutoff-cayley-<} showed that, for any prime $p$, the relaxation time of the random walk on any Cayley graph of $\mbz_p$ with respect to an arbitrary set of $k$ generators is order at least $\abs{\mbz_p}^{2/k} = p^{2/k}$, provided that $k \le \log p/\log\log p$.
Using a different approach, we extend \citeauthor{H:cutoff-cayley-<}'s result, removing the restrictions on $p$ and $k$ and considering general Abelian groups; see \cref{res-p3:intro:gap}.

This extends, in the Abelian set-up, a celebrated result of \textcite[Corollary~1]{AR:cayley-expanders}, which asserts that, for any finite group $G$, the random Cayley graph with at least $C_\eps \log\abs G$ random generators is whp an $\eps$-expander, provided $C_\eps$ is sufficiently large (in terms of $\eps$).
(A graph is an \textit{$\eps$-expander} if its isoperimetric constant is bounded below by $\eps$; up to a reparametrisation, this is equivalent to the spectral gap of the graph being bounded below by $\eps$.)
There has been a considerable line of work building upon this general result of \citeauthor{AR:cayley-expanders}.
(\textcite{P:cayley-expanders} proves a similar result.)
Their proof was simplified and extended, independently, by \textcite{LS:cayley-expanders} and \textcite{LR:cayley-expanders}; both were able to replace $\log_2 \abs G$ by $\log_2 D(G)$, where $D(G)$ is the sum of the dimensions of the irreducible representations of the group $G$;
for Abelian groups $D(G) = \abs G$.
A `derandomised' argument for Alon--Roichman is given by \textcite{CMR:cayley-expanders}.
Both \cite{CMR:cayley-expanders,LR:cayley-expanders} use some Chernoff-type bounds on operator valued random variables.

\textcite{CM:cayley-expanders} improve these further by using matrix martingales and proving a Hoeffding-type bound on operator valued random variables.
They also improved the quantification for $C_\eps$, showing that one may take $C_\eps \cq 1 + c_\eps$ with $c_\eps \to 0$ as $\eps \to 0$; this means that, whp, the graph is an $\eps$-expander whenever $k \ge (1 + c_\eps) \log_2 D(G)$ and $c_\eps \to 0$ as $\eps \to 0$.
They also generalise Alon--Roichman to random coset graphs.
The proofs use tail bounds on the (random) eigenvalues.

\textcite[Theorem~2]{AR:cayley-expanders} also specifically consider Abelian groups.
There they do a calculation directly in terms of the eigenvalues, rather than using a probabilistic tail bound.

In
\cite[Theorem~E]{HOt:rcg:abe:cutoff},
we analyse the spectral gap of $G_k$ for a nilpotent group $G$ in the regime $k \asymp \log \abs G$: we show that $G_k$ is an expander whp under a certain natural condition on $k$.
In the special case of Abelian groups, this becomes $k - d(G) \asymp k$; the general condition is $k - d(\widebar G) \asymp k$ where $\widebar G$ is the direct product of the quotients in the lower central series of $G$.
Hence in this set-up it extends \cref{res-p3:intro:gap} by removing the restriction that $\abs G$ lies in a large-density subset of $\mbn$.


\medskip

There are some fairly standard ways in which one can get bounds on the (usual) spectral gap of a Markov chain.
The first is to look at the mixing time.
For $c > 0$ and $\eps \in (0, \pi_{\min}^c]$, we have
\[
	\tmix(\eps) \asymp \trel \log(1/\eps),
\]
where $n$ is the size of the state space of the (reversible) Markov chain, $\pi_{\min}$ is the minimal value of the invariant distribution of the Markov chain and $c$ is a constant;
see, eg, \cite[Theorem~20.6 and Lemma~20.11]{LPW:markov-mixing}.
Thus, if one can bound the mixing time at level $\pi_{\min}^c$ then one can bound the relaxation time.
This method is used by \textcite{AR:cayley-expanders,P:cayley-expanders}; we use it in \cite{HOt:rcg:abe:cutoff}.

Another method is to obtain a tail estimate on the value of a random eigenvalue; one can then use the union bound to say that all (non-unitary) eigenvalues are at most some fixed value, which in turn lower bounds the spectral gap (ie upper bounds the relaxation time).

All these references consider the regime $k \asymp \log \abs G$; our results also apply when $1 \ll k \ll \log \abs G$.
From a technical perspective, in order to obtain failure probability via a large deviation bound for a random eigenvector of $\Oh{1/\abs G}$, one needs $k \gtrsim \log \abs G$.
The purpose of this is to carry out a union bound over the $\abs G$ eigenvalues;
see, eg, \cite{CM:cayley-expanders}.
Likewise, arguments that bound the $1/\abs G^c$ mixing time, for some constant $c$, in terms of some generator getting picked once (cf \cite{R:random-random-walks}) cannot work unless $k \gtrsim \log \abs G$.
As such, to consider $1 \ll k \ll \log \abs G$, a different approach is needed.
We still use a union bound, but instead of asking for an error probability $\Oh{1/\abs G}$ for each eigenvalue, we group together eigenvalues according to a certain gcd and bound the error for each group.

\subsection{Additional Remarks}
\label{sec-p3:intro:rmks}

\subsubsection{Precise Definition of Cayley Graphs}
\label{sec-p3:intro:cayley-def}

Consider a finite group $G$.
Let $Z$ be a multisubset of $G$.
We consider geometric properties,
namely through distance metrics and the spectral gap,
of the \textit{Cayley graph} of $(G,Z)$;
we call $Z$ the \textit{generators}.
The
	\textit{undirected}, respectively \textit{directed},
\textit{Cayley graph of $G$ generated by $Z$}, denoted
	$G^-(Z)$, respectively $G^+(Z)$,
is the multigraph whose vertex set is $G$ and whose edge multiset is
\[
	\sbb{ \bra{ g,g \cdot z } \mid g \in G, \, z \in Z },
\Quad{respectively}
	\sbb{ \rbr{ g,g \cdot z } \mid g \in G, \, z \in Z }.
\]

We focus attention on the \emph{random} Cayley graph defined by choosing $Z_1, ..., Z_k \sim^\iid \Unif(G)$; when this is the case, denote $G^+_k \cq G^+(Z)$ and $G^-_k \cq G^-(Z)$.
While we do not assume that the Cayley graph is connected (ie, $Z$ may not generate $G$), in the Abelian set-up the random Cayley graph $G_k$ is connected whp whenever $k - d(G) \gg 1$; see
\cite[Lemma~8.1]{HOt:rcg:abe:extra}.


The graph depends on the choice of $Z$.
Sometimes it is convenient to emphasise this;
we use a subscript, writing $\pr[G(z)]{\cdot}$ if the graph is generated by the group $G$ and multiset~$z$.
Analogously, $\pr[G_k]{\cdot}$ stands for the \emph{random} law $\pr[G(Z)]{\cdot}$ where $Z = [Z_1, ..., Z_k]$ with $Z_1, ..., Z_k \sim^\iid \Unif(G)$.

\subsubsection{Typical and Simple Cayley Graphs}
\label{sec-p3:intro:rmks:typ-simp}

The directed Cayley graph $G^+(z)$ is simple if and only if no generator is picked twice, ie $z_i \ne z_j$ for all $i \ne j$.
The undirected Cayley graph $G^-(z)$ is simple if in addition no generator is the inverse of any other, ie $z_i \ne z_j^{-1}$ for all $i,j \in [k]$.
In particular, this means that no generator is of order 2, as any $s \in G$ of order 2 satisfies $s = s^{-1}$---this gives a multiedge between $g$ and $g s$ for each $g \in G$.
Abusing terminology, we relax the definition of simple Cayley graphs to allow order 2 generators, ie remove the condition $z_i \ne z_i^{-1}$ for all $i$.

Given a group $G$ and an integer $k$,
we are drawing the generators $Z_1, ..., Z_k$ independently and uniformly at random.
It is not difficult to see that the probability of drawing a given multiset depends only on the number of repetitions in that multiset.
Thus, conditional on being simple, $G_k$ is uniformly distributed on all simple degree-$k$ Cayley graphs.
Since $k \ll \sqrt{\abs G}$, the probability of simplicity tends to 1 \asinf{\abs G}.
So when we say that our results hold ``\whp (over $Z$)'', we could equivalently say that the result holds ``for almost all degree-$k$ simple Cayley graphs of $G$''.

Our asymptotic evaluation does not depend on the particular choice of $Z$, so the statistics in question depend very weakly on the particular choice of generators for almost all choices.
In many cases, the statistics depend only on $G$ via $\abs G$ and $d(G)$.
This is a strong sense of `universality'.

\subsubsection{Overview of Random Cayley Graphs Project}
\label{sec-p3:intro:rmks:advert}

This paper is one part of an extensive project on random Cayley graphs.
There are
	three main articles \cite{HOt:rcg:abe:cutoff,HOt:rcg:matrix,HOt:rcg:abe:geom}
	(including the current one \cite{HOt:rcg:abe:geom}),
	a technical report \cite{HOt:rcg:abe:extra}
and
	a supplementary document \cite{HOt:rcg:supp}
	containing deferred technical proofs.
\textit{Each main article is readable independently.}

The main objective of the project is to establish cutoff for the random walk and determining whether this can be written in a way that, up to subleading order terms, depends only on $k$ and $\abs G$; we also study universal mixing bounds, valid for all, or large classes of, groups.
Separately, we study the distance of a uniformly chosen element from the identity, ie typical distance, and the diameter; the main objective is to show that these distances concentrate and to determine whether the value at which these distances concentrate depends only on $k$ and $\abs G$.

\begin{itemize}[noitemsep, topsep = \smallskipamount, label = \bcdot]
	\item [\cite{HOt:rcg:abe:cutoff}]
	Cutoff phenomenon (and Aldous--Diaconis conjecture) for general Abelian groups; also, for nilpotent groups, expander graphs and comparison of mixing times with Abelian groups.
	
	\item [\cite{HOt:rcg:abe:geom}]
	Typical distance, diameter and spectral gap for general Abelian groups.
	
	\item [\cite{HOt:rcg:matrix}]
	Cutoff phenomenon and typical distance for upper triangular matrix groups.
	
	\item [\cite{HOt:rcg:abe:extra}]
	Additional results on cutoff and typical distance for general Abelian groups.
	
\end{itemize}

The proofs of a number of auxiliary lemmas are deferred to the supplementary document \cite{HOt:rcg:supp}.
These are primarily of a technical and computational nature.
We do this deferral in order to keep the current manuscript as focussed as possible on the conceptually important matters.

\subsubsection{Acknowledgements}
\label{sec-p3:intro:ackno}

This whole random Cayley graphs project has benefited greatly from advice, discussions and suggestions from many of our peers and colleagues.
We thank a few of them specifically here.

\begin{itemize}[itemsep = 0pt, topsep = \smallskipamount, label = $\bcdot$]
	\item 
	Justin Salez for reading this paper in detail and giving many helpful and insightful comments as well as stimulating discussions ranging across the entire random Cayley graphs project.
	
	\item 
	Itai Benjamini for discussions on typical distance.
	
	\item 
	Evita Nestoridi and Persi Diaconis for general discussions, consultation and advice.
\end{itemize}

\section{Typical Distance: $1 \ll k \ll \log \abs G$}
\label{sec-p3:typdist1}
\renewcommand{\mm}{\ensuremath{\gamma}}

This section focusses on concentration of distances from the identity in the random Cayley graph of an Abelian group when $1 \ll k \ll \log \abs G$.
(Subsequent sections deal with $k \gtrsim \log \abs G$.)
The main result of the section is \cref{res-p3:typdist1:res}. 

\smallskip

The outline of this section is as follows:
\begin{itemize}[noitemsep, topsep = 0pt, label = \bcdot]
	\item 
	\S\ref{sec-p3:typdist1:res} states precisely the main theorem of the section;
	
	\item 
	\S\ref{sec-p3:typdist1:outline} outlines the argument;
	
	\item 
	\S\ref{sec-p3:typdist1:balls} gives some crucial estimates on the size of lattice balls;
	
	\item 
	\S\ref{sec-p3:typdist1:lower} is devoted to the lower bound;
	
	\item 
	
	\S\ref{sec-p3:typdist1:upper} is devoted to the upper bound.
\end{itemize}

\subsection{Precise Statement and Remarks}
\label{sec-p3:typdist1:res}

To start the section, we recall the typical distance statistic.

\begin{defn}
\label{def-p3:typdist1:def}
Let $H$ be a graph and fix a vertex $0 \in H$.
For $r \in \mbn$,
write $\mcb_H(r)$ for the $r$-ball in the graph $H$, ie
\(
	\mcb_H(r)
\cq
	\bra{ h \in H \mid d_H(0, h) \le r },
\)
where $d_H$ is the graph distance in $H$.
Define
\[
	\mcd_H(\beta)
\cq
	\min\brb{ r \ge 0 \midb \abs{ \mcb_H(r) } \ge \beta \abs H }
\Qfor
	\beta \in (0,1).
\]

When considering sequences $(k_N, G_N)_\Ninn$ of integers and Abelian groups,
abbreviate
\[
	\mcd_N(\beta)
\cq
	\mcd_{G_N([Z_1, ..., Z_{k_N}])}(\beta)
\Qwhere
	Z_1, ..., Z_{k_N} \sim^\iid \Unif(G_N).
\]
Finally, considering such sequences, we define the candidate radius for the typical distance:
\[
	\mfd_N^+ \cq k_N \abs{G_N}^{1/k_N} / (2e)
\Qand
	\mfd_N^- \cq k_N \abs{G_N}^{1/k_N} / e
\Quad{for each}
	\Ninn.
\]
As always, if we write $\mcd_N$, then this is either $\mcd_N^+$ or $\mcd_N^-$ according to context.
\end{defn}

We show that, whp over the graph (ie choice of $Z$), this statistic concentrates.
The result will be valid for all Abelian groups, under some conditions on $k$ in terms of $G$.
Further, the value at which the typical distance concentrates, which will be $\mfd^\pm$ above, depends only on $k$ and $\abs G$.
This is in agreement with the spirit of the Aldous--Diaconis conjecture.

\begin{hyp}
\label{hyp-p3:typdist1}
The sequence $(k_N, G_N)_\Ninn$ satisfies \textit{\cref{hyp-p3:typdist1}} if
the following hold:
\begin{gather*}
	\LIMINF\Ninf \, \abs{G_N} = \infty,
\quad
	\LIMSUP\Ninf k_N / \log \abs{G_N} = 0,
\quad
	\LIMINF\Ninf \rbr{ k_N - d(G_N) } = \infty
\\
	\text{and}
\quad
	\frac{k_N - d(G_N) - 1}{k_N} \ge 5 \frac{k_N}{\log \abs{G_N}} + 2 \frac{d(G_N) \log\log k_N}{\log \abs{G_N}}
\text{ for all }
	\Ninn.
\end{gather*}
\end{hyp}

We study $1 \ll k \ll \log \abs G$ here.
In \cref{rmk-p3:typdist1:hyp} below, we give some sufficient conditions for \cref{hyp-p3:typdist1} to hold.
Throughout the proofs, we drop the $N$-subscript from the notation, eg writing $k$ or $n = \abs G$, considering sequences implicitly.
Write $\mcd_k(\beta)$ for the $\beta$-typical distance~of~$G_k$.

We now state the main theorem of this section.

\begin{thm}
\label{res-p3:typdist1:res}
	Let $(k_N)_\Ninn$ be a sequence of positive integers and $(G_N)_\Ninn$ a sequence of finite, Abelian groups;
	for each $\Ninn$, define $Z_{(N)} \cq [Z_1, ..., Z_{k_N}]$ by drawing $Z_1, ..., Z_{k_N} \sim^\iid \Unif(G_N)$.
	
	Suppose that $(k_N, G_N)_\Ninn$ satisfies \cref{hyp-p3:typdist1}.
	Then,
	for all $\beta \in (0,1)$,
	we~have
	\[
		\mcd^\pm_N(\beta) / \mfd^\pm_N
	\to^\mbp
		1
	\Quad{(in probability)}
		\asinf N.
	\]
	Moreover, the implicit lower bound holds deterministically, ie for all choices of generators, and for all Abelian groups, ie \cref{hyp-p3:typdist1} need not be satisfied---we just need $\limsup_N k_N / \log \abs{G_N} = 0$.
\end{thm}

%

\begin{rmkt}
\label{rmk-p3:typdist1:hyp}
Write $n \cq \abs G$.
Any of the following conditions imply \cref{hyp-p3:typdist1}:
\begin{alignat*}{2}
	1 \ll k &\lesssim \sqrt{\log n / \log\log\log n}
&
	\Qand
	k - d &\gg 1;
\\
	1 \ll k &\lesssim \sqrt{\log n}
&
	\Qand
	k - d &\gg \log\log k;
\\
	1 \ll k &\ll \log n/\log\log\log n
&
	\Qand
	k - d &\ge \delta k
	\Quad{for some suitable} \delta = \oh1;
\\
	d &\ll \log n/\log\log\log n
&
	\Qand
	k - d &\asymp k.
\tag*{\qedhere}
\end{alignat*}
\end{rmkt}

\subsection{Outline of Proof}
\label{sec-p3:typdist1:outline}

As remarked after the summarised statement (in \cref{rmk-p3:intro:typdist:mixing-geom}), when considering properties of the random walk on a graph, such as the mixing time, geometric properties of the graph are often derived and used.
In a reversal of this, we use knowledge about the mixing properties of a suitable random variable to derive a geometric result.
We explain this in a little more detail now.

\smallskip

For the lower bound,
for any Cayley graph $\mcg$ of an Abelian group of degree $k$,
(trivially) we have
\(
	\abs{\mcb_\mcg(R)} \le \abs{B_k(R)},
\)
where $B_k(R)$ is the $k$-dimensional lattice ball of radius $R$.
If
\(
	\abs{B_k(R)} \ll n,
\)
then immediately
\(
	\abs{\mcb_\mcg(R)} \ll n,
\)
and so $\mcd_\mcg(\beta) \ge R$ for all $\beta \in (0,1)$, asymptotically~in~$n$.

\smallskip

Consider now the upper bound.
We fix some target radius $kL$ and draw $W_1, ..., W_k \sim^\iid \Geom(1/L)$ in the directed case.
For the undirected case, we multiply each $W_i$ by a uniform sign.
It is well-known that the law of $W \cq (W_1, ..., W_k)$ given $\norm{W}_1 = R$ is uniform on the discrete $L_1$ sphere of radius $R$.
Since the $\norm{W}_1 = \sumt[k]{1} \abs{W_i}$ is an iid sum, it concentrates around its mean, ie $k L$, when $k L \gg 1$.
So this is roughly like drawing uniformly from a sphere of radius $k L$, except that we have the added benefit that the coordinates $W_1, ..., W_k$ are (unconditionally) independent.

We can then interpret $W_i$ as the number of times which generator $i$ is used in getting from the identity to $W \bcdot Z = W_1 Z_1 + \cdots + W_k Z_k$.
We show that $W \bcdot Z$ is well-mixed \whp when
$L$ takes a certain value given below.
Now, if the law of $W \bcdot Z$ is mixed in TV and $\norm{W}_1 \le k L (1 + \delta)$ \whp, then the law of $W \bcdot Z$ conditional on $\norm{W}_1 \le k L (1 + \delta)$ is also mixed in TV.
Thus, using the concentration of $\norm{W}_1$, we deduce that a proportion $1 - \oh1$ of vertices $x \in G$ can be written as $x = w \bcdot Z$ for some $w$ with $\norm{w}_1 \eqsim k L$; this gives a path of length approximately $k L$~from~the~identity~to~$x$.

We show this mixing estimate via a (modified) $L_2$ argument, where $W$ is conditioned to be `typical', namely we define a set $\mcw$ and condition that $W \in \mcw$.
The most important part is to bound the probability that two independent copies of $W$ are equal conditional on both being in $\mcw$; this must be $\oh{1/n}$.
Since $\norm{W}_1$ concentrates and $W$ is uniform on the sphere of this radius, we need to choose $L$ so that the sphere of radius $kL$ has volume slightly more than $n$.
In high dimensions---here we consider balls in $k \gg 1$ dimensions---(discrete) spheres and balls are of asymptotically the same volume.
Thus the desired radius coincides with that of the lower bound.

In an ideal world, we would directly sample $W$ uniformly from a ball of radius $k L$.
However, the lack of independence between the coordinate causes difficulties, in particular in \cref{res-p3:typdist1:gcd-ex} below.
We thus use this vector of geometrics as a proxy for the uniform distribution, but with the key property that the coordinates are independent.

\subsection{Estimates on Sizes of Balls in $\mbz^k$}
\label{sec-p3:typdist1:balls}

We desire an $\RR^\pm_0$ so that $\abs{ B^\pm_k(\RR^\pm_0) } \approx n$,
where $B^\pm_k(R)$ is the lattice ball of radius $R$, ie
\[
	B^-_k(R)
\cq
	\brb{ w \in \mbz^k \midb \norm{w}_1 \le R }
\Quad{and}
	B^+_k(R)
\cq
	\brb{ w \in \mbz_+^k \midb \norm{w}_1 \le R }.
\]

\begin{defn}
\label{def-p3:typdist1:radius}
	Set $\omega \cq \max\bra{ \logk[2], \: k/n^{1/(2k)} }$.
	Note that $1 \ll \omega \ll k \ll \log n$.
	Define
	\[
		\RR_0^\pm
	\cq
		\inf\brb{ R \in \mbn \midb \abs{ B_k(R) } \ge n e^\omega }.
	\]
\end{defn}

The following lemma controls the size of balls.
Its proof is given in
\cite[\S E]{HOt:rcg:supp};
see in particular
\cite[Lemmas~E.2a and~E.3a]{HOt:rcg:supp}
where the index $\qq$ corresponds to a type of $L_\qq$ lattice balls; take $\qq \cq 1$ to recover the usual $L_1$ lattice balls here.
Recall $\mfd^\pm$ from \cref{def-p3:typdist1:def}.

\begin{lem}
\label{res-p3:typdist1:balls}
	Assume that $1 \ll k \ll \log n$.
	For all $\xi \in (0, 1)$,
	we have
	\[
		\abs{ \RR_0 - \mfd } / \mfd \ll 1
	\Qand
		\absb{ B_k\rbb{ \mfd (1 - \xi) } } \ll n.
	\]
\end{lem}

\subsection{Lower Bound on Typical Distance}
\label{sec-p3:typdist1:lower}

From the results in \S\ref{sec-p3:typdist1:balls}, it is straightforward to deduce the lower bound in \cref{res-p3:typdist1:res}.

\begin{Proof}[Proof of Lower Bound in \cref{res-p3:typdist1:res}]
Let $\xi \in (0, 1)$ and set $R \cq \RR_0 (1 - \xi)$.
Since the underlying group is Abelian,
applying \cref{res-p3:typdist1:balls},
we have
\(
	\abs{ \mcb_k\rbr{ R } }
\le
	\abs{ B_k\rbr{ R } }
\ll
	n.
\)
Hence,
for all $\beta \in (0,1)$ and all $Z$,
we have $\mcd_k(\beta) \ge R = \RR_0 (1 - \xi)$, asymptotically in $n$.
\end{Proof}

\subsection{Upper Bound on Typical Distance}
\label{sec-p3:typdist1:upper}

The argument given here is in a similar vein to that of
\cite[\S 2.7]{HOt:rcg:abe:cutoff};
there we analysed the mixing time of the random walk on the (random) Cayley graph.
Let $\eps > 0 $ and set $\LL \cq (1 + 3\eps) \RR_0/k$.

Draw $\WW = (\WW_i)_1^k \sim \Geom(1/\LL)^{\otimes k}$; later, we condition on $\norm{\WW}_1 \le k \LL$.
Here the geometric random variables have support $\bra{1, 2, ...}$.
Define $\chi \cq (\chi_i)_1^k$ as follows:
	in the undirected case, $\chi_i \sim^\iid \Unif(\bra{\pm1})$;
	in the directed case,   $\chi_i \cq 1$ for all $i$.
Set $S \cq (\chi \WW) \bcdot Z$ where $\chi \WW \cq (\chi_i \WW_i)_1^k$.
Define $\WW'$ and $\chi'$ as independent copies of $\WW$ and $\chi$, respectively; set $S' \cq (\chi' \WW') \bcdot Z$.

In
\cite[\S 2.7]{HOt:rcg:abe:cutoff},
a key ingredient was conditioning that the auxiliary variable $W$ was `typical' in a precise sense.
There we were interested in the law of the random walk; the introduction of typicality was a tool to study this, for establishing mixing bounds for the random walk.
Here, somewhat in reverse, we can choose which random variable we study.

%
%
%
%

\begin{defn}
	Abbreviate $\LL_0 \cq \LL \rbr{ 1 - \log k / \sqrt k }$.
	Define
	\[
		\mcw
	\cq	
		\brb{ \ww \in \mbz_+^k \midb
			\LL_0 + 1 \le \norm{\ww}_1/k \le \LL, \:
			\maxt{i} \ww_i \le 3 \LL \log k
		}.
	\]
	When $\WW$ and $\WW'$ are independent copies, write
	\(
		\typ
	\cq
		\bra{ \WW,\WW' \in \mcw }.
	\)
\end{defn}

\begin{lem}[Typicality]
\label{res-p3:typdist1:typ}
	We have
	\(
		\pr{W \in \mcw} \asymp 1
	\)
	and hence
	\(
		\prt{\typ} \asymp 1.
	\)
\end{lem}

\begin{Proof}
Recall that $\ex{W_i} = L$ and $\Var{W_i} \le L^2$, so $\norm{W}_1 = \sumt[k]{1} W_i$ is approximately distributed as $N(kL, kL^2)$.
The three parts follow easily from this representation as an iid sum of geometrics.
\begin{itemize}[noitemsep, topsep = \smallskipamount, label = \bcdot]
	\item 
	The lower bound on $\norm{W}_1$ holds with probability $1 - \oh1$ by Chebyshev's inequality.
	
	\item 
	The upper bound on $\norm{W}_1$ holds with probability bounded away from 0 by Berry--Esseen.
	
	\item 
	The upper bound on $\max_i W_i$ holds with probability $1 - \oh1$ by the union bound.
\qedhere
\end{itemize}
\end{Proof}

We control the $L_2$ distance between $S$ conditional on $W \in \mcw$ and the uniform distribution.

\begin{prop}
\label{res-p3:typdist1:l2}
	Suppose that \cref{hyp-p3:typdist1} is satisfied.
	Then
	\[
		\ex{ \normb{ \pr[G_k]{ S \in \cdot \mid W \in \mcw } - \Unif(G) }_2^2 }
	=
		\oh1,
	\]
	where we recall that $\pr[G_k]{ \cdot }$ is the \emph{random} law corresponding to the random Cayley graph $G_k$.
	Here, the expectation $\ex{ \cdot }$ is over the random choice of generators, ie over $G_k$.
\end{prop}

The proof of this proposition uses a number of auxiliary lemmas.
Given the proposition,
we now have all the ingredients to prove the upper bound on typical distance; we show this immediately.

\begin{Proof}[Proof of Upper Bound in \cref{res-p3:typdist1:res}]
Let $\widebar W$ have the law of $W$ conditional on $W \in \mcw$.
By \cref{res-p3:typdist1:l2},
the $L_2$ distance between $\widebar S \cq \widebar W \bcdot Z$ and $\Unif(G)$ is $\oh1$ \whp.
Thus the support $\mcs$ of $\widebar S$ is a proportion $1 - \oh1$ of the vertices whp.
In particular, there is a path of length at most $\LL k$ from $\id$ to all vertices in $\mcs$ \whp, as $\norm{\widebar W}_1 \le \LL k$ by definition of typicality.
Hence $\mcd_k(\beta) \le \LL k = (1 + 3 \eps) \RR_0$ whp.
Applying \cref{res-p3:typdist1:balls} then gives $\rbr{ \mcd_k(\beta) - \mfd } / \mfd \le 4 \eps$ \whp.
\end{Proof}

%

The remainder of this subsection is devoted to proving \cref{res-p3:typdist1:l2}.
We have
\[
	\ex{ \normb{ \pr[G_k]{ S \in \cdot \mid W \in \mcw } - \Unif(G) }_2^2 }
=
	n \, \pr{ S = S' \mid \typ } - 1,
\]
recalling that $\chi'$ and $W'$ are independent copies of $\chi$ and $W$, respectively, and $S' \cq (\chi' W') \bcdot Z$.
We note here that the probability on the right-hand side is annealed over the choice $Z$ of generators.

First we control the probability that $\chi \WW = \chi' \WW'$; in this case we necessarily have $S = S'$.

\begin{lem}
\label{res-p3:typdist1:V=0}
	We have
	\(
		\pr{ \chi \WW = \chi' \WW' \mid \typ } = \oh{1/n}.
	\)
\end{lem}

\begin{Proof}
Recall that $\LL_0 \cq \LL \rbr{ 1 - \log k / \sqrt k }$.
Consider the directed case first, ie $\chi = 1 = \chi'$.
Then,
\[
&	\pr{ \WW = \WW', \: \typ }
\le
	\sumt{ \ww : \norm{\ww}_1 \ge k( \LL_0 + 1 )}
	\pr{ \WW = \ww = \WW'}
\\&\qquad
=
	\sumt{\ww : \norm{\ww}_1 \ge k( \LL_0 + 1 )}
	\pr{ \WW' = \ww }
	\prodt[k]{1} \pr{ \WW_i = \ww_i }
\\&\qquad
=
	\sumt{\ww : \norm{\ww}_1 \ge k( \LL_0 + 1 )}
	\pr{ \WW' = \ww }
	\prodt[k]{1} \LL^{-1} \rbr{ 1 - \LL^{-1} }^{\ww_i - 1}
\\&\qquad
=
	\sumt{\ww : \norm{\ww}_1 \ge k( \LL_0 + 1 )}
		\pr{ \WW' = \ww }
	\cdot
		\LL^{-k} \rbr{ 1 - \LL^{-1} }^{\norm{\ww}_1 - k}
\\&\qquad
\le
	\LL^{-k} \rbr{ 1 - \LL^{-1} }^{k \LL_0}
=
	\rbb{ \LL^{-1} \rbr{ 1 - \LL^{-1} }^{\LL(1 - \sqrt{\log k/k})} }^k
\\&\qquad
\le
	(e \LL)^{-k} \expb{ \sqrt{k \log k} }
\le
	n^{-1} e^{-\eps k},
\]
with the final inequality using the fact that $\LL \ge (1 + 2\eps) n^{1/k}/e$.
Indeed,
	$L = (1 + 3\eps) \RR_0 / k$ by definition
and
	$\RR_0 \ge \mcd (1 - \eps)$ by \cref{res-p3:typdist1:balls},
recalling that $\mcd = \mcd^- = k n^{1/k} / e$ from \cref{def-p3:typdist1:def}.

In the undirected case, we also need to impose $\chi = \chi'$, which happens with probability $2^{-k}$, and is independent of $(\WW, \WW')$.
Hence, the same inequality holds with the event $\bra{\WW = \WW'}$ replaced by $\bra{\chi \WW = \chi' \WW'}$, recalling that $2 \mcd^- = \mcd^+$ from \cref{def-p3:typdist1:def}.
Finally,
\(
	\pr{ \typ }
\asymp
	1.
\)

Thus, Bayes's rule combined with the above calculation gives
\[
	\pr{ \chi \WW = \chi' \WW' \mid \typ }
\le
	n^{-1} e^{-\eps k}
\ll
	1/n.
\qedhere
\]
\end{Proof}

The following lemma describing the distribution of $\vv \bcdot Z$ for a given $\vv \in \mbz^k$ is crucial.

\begin{lem}
\label{res-p3:typdist1:unif-gcd}
	For all $\vv \in \mbz^k$ with $\gcd(\vv_1, ..., \vv_k, n) = \gamma$,
	we have
	\[
		\vv \bcdot Z \sim \Unif(\gamma G).
	\]
\end{lem}

\begin{Proof}[Sketch of Proof]
We can decompose $G$ as $\oplus_1^d \: \mbz_{m_j}$.
Now, $(U_1, ..., U_d) \sim \Unif(\oplus_1^d \: \gamma \mbz_{m_j})$ if and only if $U_j \sim \Unif(\gamma \mbz_{m_j})$ independently over $j \in [d]$.
The different coordinates of $\vv \bcdot Z$ are independent and $\gamma \mbz_{m_j} = \gcd(\gamma, m_j) \mbz_{m_j} = \gcd(\vv_1, ..., \vv_k, m_j) \mbz_{m_j}$.
Thus, it suffices to consider $d = 1$, ie $G = \mbz_n$.

We use induction on $k$:
	write $\vv \bcdot Z = (\vv_1 Z_1 + \cdots + \vv_{k-1} Z_{k-1}) + \vv_k Z_k$.
By the induction hypothesis,
	$\vv_1 Z_1 + \cdots + \vv_{k-1} Z_{k-1} \sim \Unif(a \mbz_n)$
where
	$a \cq \gcd(\vv_1, ..., \vv_{k-1}, n) \wr n$
and
	$\vv_k Z_k \sim \Unif(b \mbz_n)$
where
	$b \cq \gcd(\vv_k, n) \wr n$.
It then suffices to show that $aX + bY \sim \Unif(c \mbz_n)$ where $X, Y \sim^\iid \Unif(\mbz_n)$ and $c \cq \gcd(a, b) \wr n$.
Let $a' \cq a c^{-1}$ and $b' \cq b c^{-1}$ so that $a X + b Y = c(a'X + b'Y)$ and $\gcd(a', b') = 1$.
The fact that $a' X + b' Y \sim \Unif(\mbz_n)$~follows, in essence, from Euclid's algorithm.

The full details, fleshing out all the rigorous details, are given in
\cite[Lemma~F.1]{HOt:rcg:supp}.
	%
\end{Proof}

We thus now need to control $\abs{\mm G}$,
as
\(
	\pr{ \mm U = x } = 1 / \abs{\mm G}
\)
if $U \sim \Unif(G)$ and $\gcd(\mm, n) \wr x$.

\begin{lem}
\label{res-p3:typdist1:G/gammaG}
	For all Abelian groups $G$ and all $\gamma \in \mbn$,
	we have
	\[
		\abs G/\abs{\gamma G}
	\le
		\gamma^{d(G)}.
	\]
\end{lem}

\begin{Proof}
	Decompose $G$ as $\oplus_1^d \: \mbz_{m_j}$ with $d = d(G)$ and some $m_1, ..., m_d \in \mbn$.
	Then $\mm G$ can be decomposed as $\oplus_1^d \: \gcd(\mm, m_j) \mbz_{m_j}$.
	Hence
	\(
		\abs{\mm G}
	=
		\prodt[d]{1} \rbr{ m_j/\gcd(\mm, m_j) }
	\ge
		\prodt[d]{1} \rbr{ m_j/\mm }
	=
		\abs G / \mm^d.
	\)
\end{Proof}

These two lemmas are used in
\cite[\S 2.7]{HOt:rcg:abe:cutoff}
in an analogous way to here.
Define
\[
	\VV \cq \chi \WW - \chi' \WW'
\Qand
	\mfgcd \cq \gcd(\VV_1, ..., \VV_k, n).
\]

\begin{cor}
\label{res-p3:typdist1:pr-gcd}
	We have
	\[
		n \, \pr{ \VV \bcdot Z = 0, \: \VV \ne 0 \mid \typ }
	\le
		\ex{ \mfgcd^{d(G)} \one{\VV \ne 0} \mid \typ }.
	\]
\end{cor}

\begin{Proof}
The conditioning does not affect $Z$.
The corollary follows from \cref{res-p3:typdist1:unif-gcd,res-p3:typdist1:G/gammaG}:
\[
	n \,
	\pr{ \VV \bcdot Z = 0, \: \VV \ne 0, \: \typ }
&
=
	n \,
	\ex{ \pr{ \VV \bcdot Z = 0 \mid V } \one{V \ne 0} \one{\typ} }
\\&
=
	\ex{ \rbb{ \abs{G} / \abs{\mfgcd G} } \one{V \ne 0} \one{\typ} }
\\&
\le
	\ex{ \mfgcd^{d(G)} \one{V \ne 0} \one{\typ} }.
\qedhere
\]
\end{Proof}

\begin{lem}
\label{res-p3:typdist1:gcd-ex}
	Given \cref{hyp-p3:typdist1},
	we have
	\(
		\ex{ \mfgcd^{d(G)} \one{\VV \ne 0} \mid \typ }
	=
		1 + \oh1.
	\)
\end{lem}

\begin{Proof}
Each coordinate of $\VV$ is unimodal and symmetric about 0.
This means that we can write
\[
	\abs{\VV_1}
\sim
	\Unif\bra{1, ..., Y}
\Quad{conditional on}
	\VV_1 \ne 0,
\]
where $Y$ is a certain $\mbn$-valued random variable.
This implies that
\[
	\pr{ \VV_1 \in \gamma \mbz \mid \VV_1 \ne 0 }
=
	\ex{ \floor{Y/\gamma} / T }
\le
	1/\gamma.
\]
The probability of $\VV_1 = 0$ is roughly $1/(2\LL) \asymp n^{-1/k}$;
in particular, it is at most $3 n^{-1/k}$.
The coordinates are independent.
Since
\(
	\pr{ \typ }
\asymp
	1,
\)
we thus have
\[
	\pr{ \mfgcd = \gamma \mid \typ }
\lesssim
	\rbb{ 1/\gamma + 3/n^{1/k} }^k.
\]

We have $\mfgcd = \gcd(V_1, ..., V_k, n) \le \min_i \VV_i$.
So, $\mfgcd \le 6 \LL \log k \le 3 n^{1/k} \log k$ under typicality.
Hence,
\[
	\ex{ \mfgcd^d \one{ V \ne 0 } \mid \typ }
\lesssim
	\sumt[3 n^{1/k} \log k]{\gamma=1}
	\gamma^d \rbb{ 1/\gamma + 3/n^{1/k} }^k.
\]
We handle almost exactly the same sum in
\cite[Corollary~2.15]{HOt:rcg:abe:cutoff}.
\cref{hyp-p3:typdist1} here is designed precisely to control this sum; it is identical to
\cite[Hypothesis~A]{HOt:rcg:abe:cutoff}.
There the $3/n^{1/k}$ part is replaced with $2/n^{1/k}$, but exactly the same arguments apply showing that the sum is $1 + \oh1$.
	%
\end{Proof}

\cref{res-p3:typdist1:l2} now follows immediately from \cref{res-p3:typdist1:V=0,res-p3:typdist1:pr-gcd,res-p3:typdist1:gcd-ex}.

\begin{Proof}[Proof of \cref{res-p3:typdist1:l2}]
By \cref{res-p3:typdist1:V=0,res-p3:typdist1:pr-gcd,res-p3:typdist1:gcd-ex}, we have
\[
	n \, \pr{ S = S' \mid \typ }
&
\le
	n \, \pr{ V = 0 \mid \typ }
+	n \, \pr{ V \bcdot Z = 0, \, V \ne 0 \mid \typ }
\\&
\le
	n \, \pr{ V = 0 \mid \typ }
+	\ex{ \mfgcd^d \, \one{V \ne 0} \mid \typ }
=
	1 + \oh1.
\qedhere
\]
\end{Proof}

\section{Typical Distance: $k \asymp \log \abs G$}
\label{sec-p3:typdist2}

This section focusses on concentration of distances from the identity in the random Cayley graph of an Abelian group when $k \asymp \log \abs G$.
(The previous section dealt with $1 \ll k \ll \log \abs G$ and the next deal with $k \gg \log \abs G$.)
The main result of the section is \cref{res-p3:typdist2:res}; see also~\cref{hyp-p3:typdist2}.

\smallskip

The outline of this section is as follows:
\begin{itemize}[noitemsep, topsep = 0pt, label = \bcdot]
	\item 
	\S\ref{sec-p3:typdist2:res} states precisely the main theorem of the section;
	
	\item 
	\S\ref{sec-p3:typdist2:outline} outlines the argument;
	
	\item 
	\S\ref{sec-p3:typdist2:balls} gives some crucial estimates on the size of lattice balls;
	
	\item 
	\S\ref{sec-p3:typdist2:lower} is devoted to the lower bound;
	
	\item 
	\S\ref{sec-p3:typdist2:upper} is devoted to the upper bound under additional constraints;
	
	\item 
	\S\ref{sec-p3:typdist2:relax-m*} describes how to relax these additional constraints;
	
	\item 
	\S\ref{sec-p3:typdistq} describes an extension for $L_1$-type graph distances to $L_\qq$-type.
\end{itemize}

\subsection{Precise Statement and Remarks}
\label{sec-p3:typdist2:res}

To start the section, we recall the typical distance statistic.

\begin{defn}
\label{def-p3:typdist2:def}
Let $H$ be a graph and fix a vertex $0 \in H$.
For $r \in \mbn$,
write $\mcb_H(r)$ for the $r$-ball in the graph $H$, ie
\(
	\mcb_H(r)
\cq
	\bra{ h \in H \mid d_H(0, h) \le r },
\)
where $d_H$ is the graph distance in $H$.
Define
\[
	\mcd_H(\beta)
\cq
	\min\brb{ r \ge 0 \midb \abs{ \mcb_H(r) } \ge \beta \abs H }
\Qfor
	\beta \in (0,1).
\]

When considering sequences $(k_N, G_N)_\Ninn$ of integers and Abelian groups,
abbreviate
\[
	\mcd_N(\beta)
\cq
	\mcd_{G_N([Z_1, ..., Z_{k_N}])}(\beta)
\Qwhere
	Z_1, ..., Z_{k_N} \sim^\iid \Unif(G_N).
\]
As always, if we write $\mcd_N$, then this is either $\mcd_N^+$ or $\mcd_N^-$ according to context.
\end{defn}

We show that, whp over the graph (ie choice of $Z$), this statistic concentrates.
Here we consider $k \eqsim \lambda \log \abs G$ for any $\lambda \in (0,\infty)$.
The result holds for a large class of Abelian groups.
Further, for these groups, the typical distance concentrates at $\alpha_\lambda k$ where $\alpha_\lambda \in (0,\infty)$ is a constant; so this depends only on $k$ and $\abs G$.
This is in agreement with the spirit of the Aldous--Diaconis conjecture.

Recall that any Abelian group can be decomposed as $\oplus_{j=1}^d \: \mbz_{m_j}$ for some $d, m_1, ..., m_d \in \mbn$.
For an Abelian group $G$, we define the \textit{dimension} and \textit{minimal side-length}, respectively, as follows:
\[
	d(G)
&\cq
	\min\brb{ d \in \mbn \midb \oplus_{j=1}^d \: \mbz_{m_j} \text{ is a decomposition of } G };
\\
	m_*(G)
&\cq
	\max\brb{ \mint{j \in [d]} m_j \midb \oplus_{j=1}^d \: \mbz_{m_j} \text{ is a decomposition of } G }.
\]
It can be shown that there is a decomposition which is optimal for both these statistics:
	there exist $d, m_1, ..., m_d \in \mbn$ so that $\oplus_{j=1}^d \: \mbz_{m_j}$ is a decomposition of $G$ with $d = d(G)$ and $\mint{j \in [d]} m_j = m_*(G)$.
From now on, we assume that we are always using such an optimal decomposition.

There are some conditions which the Abelian groups must satisfy.

\begin{hyp}
\label{hyp-p3:typdist2}
	The sequence $(k_N, G_N)_{N \in \mbn}$ satisfies \textit{\cref{hyp-p3:typdist2}} if
	\begin{gather*}
		\LIM{\Ninf} k_N = \infty,
	\quad
		\LIM{\Ninf} k_N / \log \abs{G_N} \in (0,\infty),
	\quad
		\LIMINF{\Ninf} m_*(G_N) = \infty
	\\
		\text{and}
	\quad
		d(G_N) \le \tfrac14 \log \abs{G_N} / \log\log \abs{G_N}
	\Qforall
		N \in \mbn.
	\end{gather*}
\end{hyp}

We are now ready to state the main theorem of this section.

\begin{thm}
\label{res-p3:typdist2:res}
	Let $(k_N)_\Ninn$ be a sequence of positive integers and $(G_N)_\Ninn$ a sequence of finite, Abelian groups;
	for each $\Ninn$, define $Z_{(N)} \cq [Z_1, ..., Z_{k_N}]$ by drawing $Z_1, ..., Z_{k_N} \sim^\iid \Unif(G_N)$.
	
	Suppose that $(k_N, G_N)_\Ninn$ satisfies \cref{hyp-p3:typdist2}.
	Let $\lambda \cq \limsup_N k_N / \log \abs{G_N}$.
	Then there exists a constant $\alpha^\pm_\lambda \in (0,\infty)$ so that,
	for all $\beta \in (0,1)$,
	we have
	\[
		\mcd^\pm_N(\beta) / (\alpha^\pm_\lambda k_N)
	\to^\mbp
		1
	\Quad{(in probability)}
		\asinf N.
	\]
	Moreover, the implicit lower bound holds deterministically, ie for all choices of generators, and for all Abelian groups, ie \cref{hyp-p3:typdist2} need not be satisfied---we just need $\lim_N k_N / \log \abs{G_N} \in (0,\infty)$.
\end{thm}

For ease of presentation, in the proof we drop the $N$-subscripts.

\begin{rmkt}
\label{rmk-p3:typdist2:Lq}
	In \S\ref{sec-p3:typdistq}, we describe an extension from the usual $L_1$-type graph distances to $L_\qq$-type.
	An analogous concentration of typical distance is given.
	See Hypothesis~\ref{hyp-p3:typdistq} and \cref{res-p3:typdistq:res}.
\end{rmkt}

\subsection{Outline of Proof}
\label{sec-p3:typdist2:outline}

The outline here is very similar to that from before; see \S\ref{sec-p3:typdist1:outline}.
In particular, the lower bound is exactly the same idea.
For the upper bound, we were trying to bound the expectation of a $d$-th power of a gcd.
Issues arose when $k$ became too large while $k - d$ is fairly small; see the proof of \cref{res-p3:typdist1:gcd-ex}.
Particularly, the factor $\gamma^d$ needs to be countered by $(1/\gamma + 3/n^{1/k})^k$ in a suitably strong sense.
This arose from the fact that we used the estimate
\[
	\pr{ V_1 \in \mm \mbz }
\le
	\pr{ V_1 \in \mm \mbz \mid V_1 \ne 0 } + \pr{ V_1 = 0 }
\le
	1/\mm + 3/n^{1/k}.
\]
Once this was raised to the power $k$, the second term became an issue.
We alleviate this by defining
\[
	\mci
\cq
	\brb{ i \in [k] \mid V_i \ne 0 }
\Quad{and studying}
	\pr{ V_i \in \mm \mbz \mid i \in \mci }.
\]
The problematic term $3/n^{1/k}$ then does not exist as we consider only non-zero coordinates of $V$.

If $G = \oplus_{j=1}^d \: \mbz_{m_j}$, then we are actually interested in $V_i \mod m_j$ for each $j$.
Recall that $m_* = \min_j m_j$.
`Typically', one has $\abs{\VV_i} \le m_*$.
Indeed, $k \asymp \log \abs G$, which means that $\ex{\abs{\VV_i}} \asymp 1$, but $m_* \gg 1$ by assumption.
We suppose initially that $m_*$ is large enough so that $\max_i \abs{V_i} < m_*$ whp.
Thus looking at $V_i = 0$ or $V_i \equiv 0 \mod m_j$ is no different.

For large $\abs \mci$, the gcd analysis goes through similarly to before.
When $\abs \mci$ is small, eg smaller than $d$, it is more difficult to control; in this case, we use a fairly naive bound on the gcd, but control carefully the probability of realising such an $\mci$.
The case $\mci = \emptyset$, which corresponds to $V = 0$, is handled by taking the lattice ball to be of large enough volume.

Previously we used a vector of geometrics as a proxy for a uniform distribution on a ball.
Here we are able to let $W$ be uniform on a ball.
The coordinates are no longer independent, which makes the gcd analysis slightly complicated.
However, since we only consider $i$ with $V_i \ne 0$, this can be handled; see \cref{res-p3:typdist2:gcd-ex}.
This uniformity simplifies some other calculations somewhat.

\subsection{Estimates on Sizes of Balls in $\mbz^k$}
\label{sec-p3:typdist2:balls}

We wish to determine the size of balls $B_k(R)$ when $k \asymp \log n$.
In particular, we are interested in the growth when the volume is around $n$.

\begin{defn}
\label{def-p3:typdist2:M*-def}
	Define $M^\pm_*(k, N)$ to be the minimal integer $M$ satisfying $\abs{ B^\pm_k(M) } \ge N$.
\end{defn}

\begin{lem}
\label{res-p3:typdist2:balls}
	For all $\lambda \in (0,\infty)$,
	there exists a function $\omega \gg 1$ and a constant $\alpha^\pm$ so that,
	for all $\eps \in (0,1)$,
	if $k \eqsim \lambda \log n$, then
	$\mcm^\pm_* \cq M^\pm_*(k, n e^\omega)$ satisfies
	\[
		\mcm^\pm_* \eqsim \alpha^\pm k \eqsim \alpha^\pm \lambda \log n
	\Qand
		\absb{ B^\pm_k\rbb{ \alpha^\pm k (1 - \eps) } } \ll n.
	\]
\end{lem}

This will follow easily from the following auxiliary lemma controlling the size of lattice balls.

%

\begin{lem}
	There exists a strictly increasing, continuous function $c^\pm : (0,\infty) \to (0,\infty)$ so that,
	for all $a \in (0,\infty)$,
	we have
	\[
		\absb{ B^\pm_k(ak) }
	=
		\expb{ k \rbb{ c^\pm(a) + \oh1 } }.
	\]
\end{lem}

\begin{Proof}
The directed case follows immediately from Stirling's approximation and the fact that
\[
	\absb{ B^+_k(ak) }
=
	\absb{ \brb{ b \in \mbz_+^k \midb \sumt[k]{1} b_i \le ak } }
=
	\binomt{\floor{ak}+k}{k}
=
	\binomt{\floor{(a+1)k}}{k}
\]

Consider now the undirected case.
Omit all floor and ceiling signs.
By considering the number of coordinates which equal 0, we obtain
\[
	\absb{ B^-_k(ak) }
=
	\sumt[k]{i=0} A_i
\Qwhere
	A_i
\cq
	A_i(k,a)
\cq
	\binomt ki 2^{k-i} \binomt{k-i+ak}{ak}.
\]
Choose $i_* \cq i_*(k,a)$ that maximises $A_i$.
Then
\(
	A_{i_*} \le \abs{ B^-_k(ak) } \le (k+1) A_{i^*}.
\)
Observe that
\[
	\frac{A_{i+1}}{A_i}
=
	\frac{(k-i)^2}{2(i+1)(k(1+a)-i)},
\]
and hence one can determine $i_*$ as a function of $k$ and $a$, conclude that $i_*(a,k)/k$ converges as $k \to \infty $ and thus determine $c^+(a)$ in terms of the last limit. We omit the details.
Knowing this limit allows us to plug this into the definition of $A_i$ and use Stirling's approximation to get
\[
	A_{i_*} = \expb{ k \rbb{ c^-(a) + \oh1 } },
\]
for some strictly increasing function $c^- : (0,\infty) \to (0,\infty)$.
Since $k+1 = e^{\oh{k}}$, the claim follows.
\end{Proof}

From this lemma, \cref{res-p3:typdist2:balls} follows easily.

\begin{Proof}[Proof of \cref{res-p3:typdist2:balls}]
Set $\alpha \cq c^{-1}(1/\lambda)$.
The upper bound is an immediate consequence of the continuity of $c$.
The lower bound follows from the exponential growth rate.
\end{Proof}

\subsection{Lower Bound on Typical Distance}
\label{sec-p3:typdist2:lower}

From the results in \S\ref{sec-p3:typdist2:balls}, it is straightforward to deduce the lower bound in \cref{res-p3:typdist2:res}.

\begin{Proof}[Proof of Lower Bound in \cref{res-p3:typdist2:res}]
Let $\xi \in (0,1)$ and set $R \cq \alpha^\pm_\lambda k (1 - \xi)$.
Since the underlying group is Abelian,
applying \cref{res-p3:typdist2:balls},
we have
\(
	\abs{ \mcb^\pm_k(R) } \le \abs{ B^\pm_k(R) } \ll n.
\)
Hence,
for all $\beta \in (0,1)$ and all $Z$,
we have $\mcd^\pm_k(\beta) \ge R = \alpha^\pm_\lambda k (1 - \xi)$, asymptotically in $n$.
\end{Proof}

\subsection{Upper Bound on Typical Distance Given $m_*(G) \gg k$}
\label{sec-p3:typdist2:upper}

Define $\mcm^\pm_*$, $\omega$ and $\alpha^\pm$ as in \cref{def-p3:typdist2:M*-def,res-p3:typdist2:balls}.
In this subsection we draw $W^\pm \sim \Unif\rbr{ B^\pm_k(\mcm^\pm_*) }$, ie uniform on a ball of radius $\mcm^\pm_*$.
We show that $W^\pm \bcdot Z$ is well-mixed on $G$, and hence its support contains almost all the vertices.

\begin{prop}
\label{res-p3:typdist2:l2}
	Suppose that \cref{hyp-p3:typdist2} is satisfied.
	Then
	\[
		\ex{ \normb{ \pr[G_k]{ W^\pm \bcdot Z \in \cdot } - \Unif(G) }_2^2 } = \oh1,
	\]
\end{prop}

Given this proposition, the upper bound in \cref{res-p3:typdist2:res} follows easily.

\begin{Proof}[Proof of Upper Bound in \cref{res-p3:typdist2:res} Given \cref{res-p3:typdist2:l2}]
The support $\mcs$ of $W^\pm \bcdot Z$ satisfies $\abs{\mcs^c}/n \le \eps$
if
\(
	\norm{ \pr[G_k]{ W^\pm \bcdot Z \in \cdot } - \Unif(G) }_2 \le \eps.
\)
Combined with \cref{res-p3:typdist2:balls,res-p3:typdist2:l2}, the upper bound in \cref{res-p3:typdist2:res} follows.
\end{Proof}

The remainder of this subsection is devoted to proving \cref{res-p3:typdist2:l2}.
We tend to drop the $\pm$-superscript from the notation, only writing $+$ or $-$ if there is ambiguity.
Let $W,W' \sim^\iid \Unif\rbr{ B_k(\mcm_*) }$ and let $V \cq W - W'$.
The standard $L_2$ calculation gives
\[
	\ex{ \normb{ \pr[G_k]{ W \bcdot Z \in \cdot } - \Unif(G) }_2^2 }
=
	\ex{ n \, \pr{ V \bcdot Z = 0 \mid Z } - 1 }
=
	n \, \pr{ V \bcdot Z = 0 } - 1.
\]

First, it is immediate that
\(
	\pr{ V = 0 }
=
	\pr{ W = W' }
=
	\abs{ B_k(\mcm_*) }^{-1}
\le
	n^{-1} e^{-\omega}
\ll
	n^{-1}.
\)
Now consider $V \ne 0$.
As in \S\ref{sec-p3:typdist1:upper}, it is key to analyse certain gcds.
In this section, we set
\[
	\mfgcd_j \cq \gcd\rbb{ V_1, ..., V_k, m_j }
\Quad{for each}
	j \in [d];
\Quad{set}
	\mfgcd \cq \gcd\rbb{ V_1, ..., V_k, n }.
\]
The following lemma is equivalent to \cref{res-p3:typdist1:unif-gcd}, rephrased slightly.

\begin{lem}
\label{res-p3:typdist2:unif-gcd}
	Conditional on $V$, we have
	\(
		V \bcdot Z \sim \Unif\rbr{ \oplus_{j=1}^d \: \mfgcd_j \mbz_{m_j} }.
	\)
\end{lem}


For the remainder of this subsection, we assume that the minimal side-length $m_* \cq m_*(G)$ satisfies $m_* \gg k \asymp \mcm_*$.
In the next subsection, we remove this assumption: we extend the proof to $m_* \gg 1$, as in \cref{hyp-p3:typdist2}.
Given $m_* \gg k$, we have $\max_{i \in [k]} \abs{V_i} < \min_{j \in [d]} m_j$.
Hence,
\[
	\mci
\cq
	\brb{ i \in [k] \midb V_i \not\equiv 0 \mod m_j \: \forall \, j \in [d] }
=
	\brb{ i \in [k] \midb W_i \ne W'_i }.
\]
To analyse the expected gcd, we breakdown according to the value of $\mci$.

\begin{lem}
\label{res-p3:typdist2:gcd-ex}
	There exists a constant $C$ so that,
	for all $I \subseteq [k]$ with $I \ne \emptyset$,
	we have
	\[
		n \, \pr{ V \bcdot Z = 0 \mid \mci = I }
	\le
		\ex{ \mfgcd^d \mid \mci = I }
	\le
	\begin{cases}
		C 2^d (2 \mcm_*)^{d - \abs I + 2}
	&\text{when}\quad
		\abs I \le d + 1,
	\\
		1 + 5 \cdot (\tfrac32)^{2d - \abs I}
	&\text{when}\quad
		\abs I \ge d + 2.
	\end{cases}
	\]
\end{lem}

%
%
%

\begin{lem}
\label{res-p3:typdist2:I}
	For all $I \subseteq [k]$ with $\abs I \ll k$,
	we have
	\(
		\pr{ \mci = I }
	\le
		e^{-\omega} n^{-1 + \oh1}.
	\)
	If $I = \emptyset$, then the $\oh1$ term may be taken to be 0.
\end{lem}

Given these two lemmas, we have all the ingredients required to prove \cref{res-p3:typdist2:l2}, from which we deduced the main theorem (\cref{res-p3:typdist2:res}).
We defer the proofs of \cref{res-p3:typdist2:gcd-ex,res-p3:typdist2:I} until after the proof of \cref{res-p3:typdist2:l2}, which we give now.

\begin{Proof}[Proof of \cref{res-p3:typdist2:l2}]
Here $k \eqsim \lambda \log n$, $M \cq \mcm_* \eqsim \alpha k \eqsim \alpha \lambda \log n$ and $d \le \tfrac14 \log n / \log\log n$.

As noted previously,
the standard $L_2$ calculation gives
\[
	\ex{ \normb{ \pr[G_k]{ W \bcdot Z \in \cdot } - \Unif(G) }_2^2 }
&
=
	\ex{ n \, \pr{ V \bcdot Z = 0 \mid Z } - 1 }
\\&
=
	n \, \pr{ V \bcdot Z = 0 } - 1
=
	n \sumt{I \subseteq [k]}
	\pr{ V \bcdot Z = 0, \: \mci = I }
-	1.
\]

Consider $I = \emptyset$.
Then $V \bcdot Z = 0$ (for all $Z$).
By \cref{res-p3:typdist2:I},
we have
\(
	\pr{ \mci = \emptyset }
\le
	n^{-1} e^{-\omega}.
\)
Thus
\[
	n \, \pr{ V \bcdot Z = 0, \: \mci = \emptyset }
\le
	e^{-\omega}
=
	\oh1.
\]

Consider $I \subseteq [k]$ with $1 \le \abs I \le d + 1$.
There are at most $(d+1) \binom{k}{d+1} \le k^{d+2}$ such sets $I$.
Since $\log k = \log\log n + \log \lambda + \oh1$,
we have
\(
	k^{d+2} \le n^{2/3}.
\)
Applying \cref{res-p3:typdist2:gcd-ex,res-p3:typdist2:I} gives
\[
	n \, \pr{ V \bcdot Z = 0, \: \mci = I }
\le
	C 2^d \rbr{ 3 \alpha \lambda \log n }^{d + 2 - \abs I} \cdot n^{-1+\oh1}
\le
	k^{-d-2} n^{-1/4},
\]
noting that $d \ll k \asymp \log n$ and so $2^d = n^{\oh1}$.
We now sum over all $I$ with $1 \le \abs I \le d + 1$:
\[
	n \sumt{1 \le \abs I \le d + 1}
	\pr{ V \bcdot Z = 0, \: \mci = I }
\le
	n^{-1/4}
=
	\oh1.
\]

Consider $I \subseteq [k]$ with $d + 2 \le \abs I \le L \cq \tfrac23 \log n / \log \log n$; then $L - 2d \gg 1$.
Similarly to above, there are at most $L \binom kL \le k^{L+1}$ such sets $I$.
Applying \cref{res-p3:typdist2:gcd-ex,res-p3:typdist2:I}~gives
\[
	n \, \pr{ V \bcdot Z = 0, \: \mci = I }
\le
	n^{-1 + \oh1}
\le
	k^{-L-1} n^{-1/4},
\]
noting that $k^L \le n^{2/3+\oh1}$.
We now sum over all $I$ with $d + 2 \le \abs I \le L$:
\[
	n \sumt{d + 2 \le \abs I \le L}
	\pr{ V \bcdot Z = 0, \: \mci = I }
\le
	n^{-1/4}
=
	\oh1.
\]

Finally consider $I \subseteq [k]$ with $\abs I \ge L$.
Sum over these using \cref{res-p3:typdist2:gcd-ex}:
\[
	n \sumt{L \le \abs I \le k}
	\pr{ V \bcdot Z = 0, \: \mci = I }
\le
	1 + 5 \cdot (\tfrac32)^{2d - L}
=
	1 + \oh1.
\]

Combining these four parts into a single sum, we deduce the result.
\end{Proof}

It remains to prove the auxiliary \cref{res-p3:typdist2:gcd-ex,res-p3:typdist2:I}.

\begin{Proof}[Proof of \cref{res-p3:typdist2:gcd-ex}]
The first inequality is an immediate consequence of \cref{res-p3:typdist2:unif-gcd}.

Note that $\mfgcd \le 2 \mcm_*$ since $\max_i \abs{V_i} \le 2 \mcm_*$.
For $\alpha, \beta \in \mbz$, write $\alpha \wr \beta$ if $\alpha$ divides $\beta$.
Thus
\[
	\ex{ \mfgcd^d \midb \mci = I }
\le
	\sumt[2M]{\mm = 1}
	\mm^d \, \pr{ \mm \wr V_i \: \forall \, i \in I \mid \mci = I }
\]


For a set $I \subseteq [k]$,
write
\(
	W_I \cq (W_i)_{i \in I}
\)
and
\(
	W_{\setminus I} \cq W_{[k] \setminus I}.
\)
Consider conditioning on $\mci = I$.
Let $W_{\setminus I}$ and $W'_{\setminus I}$ be given; since $\mci = I$, we have $W_{\setminus I} = W'_{\setminus I}$.
Let $U$ have the distribution of $W_I$ given $W_{\setminus I}$ and define $U'$ analogously.
Write
\(
	D_i \cq D_i(\mm) \cq \bra{ \mm \wr (U_i - U'_i) }.
\)
Then
\[
	\pr{ \mm \wr V_i \: \forall \, i \in I \midb \mci = I, \, \norm{W_{\setminus I}}_1 }
=
	\pr{ D_i \: \forall \, i \in I }.
\]
By exchangeability, it suffices to consider the case $I = \bra{1, ..., \ell}$.
We then have
\[
	\pr{ D_i \: \forall \, i \in I }
=
	\pr{ D_\ell } \pr{ D_{\ell-1} \midb D_\ell } \cdots \pr{ D_1 \midb D_2, ..., D_\ell }
=
	\prodt[\ell]{i=1} \pr{ D_i \midb D_{i+1}, ..., D_\ell }.
\]

For $i \in [k]$, define
\(
	M_i \cq \mcm_* - \norm{ W_{\setminus \bra{1, ..., i}} }_1
\)
and $M_i'$ analogously.
Let $i \in [\ell-1]$.
Let $(u_{i+1}, ..., u_\ell)$ and $(u'_{i+1}, ..., u'_\ell)$ be two vectors in the support of $(U_{i+1}, ..., U_\ell)$.
Then,
\begin{gather*}
	\text{conditional on}
\quad
	(U_{i+1}, ..., U_\ell) = (u_{i+1}, ..., u_\ell)
\Qand
	(U'_{i+1}, ..., U'_\ell) = (u'_{i+1}, ..., u'_\ell),
\\
	\text{we have}
\quad
	(U_1, ..., U_i) \sim \Unif\rbb{ B_i(R) }
\Qand
	(U'_1, ..., U'_i) \sim \Unif\rbb{ B_i(R') }
\Quad{for some}
	R,R' \in \mbr.
\end{gather*}
(Recall that the subscript in $B_k$ denotes the dimension of the ball.)

In the case of undirected balls, the law of $U_i - U'_i$ given this conditioning is symmetric and unimodal on $\mbz \setminus \bra{0}$; see \cite[Theorem~2.2]{P:unimodal}.
It follows,
as in the proof of \cref{res-p3:typdist1:gcd-ex},
that
\[
	\pr{D^-_i \midb D^-_{i+1}, ..., D^-_\ell} \le 1/\mm.
\]
Further, this holds not just conditional on $D^-_{i+1} \cap \cdots \cap D^-_\ell$, but conditional on any choice of $(U_{i+1}, ..., U_\ell)$ and $(U'_{i+1}, ..., U'_\ell)$ which satisfy $D^-_{i+1} \cap \cdots \cap D^-_\ell$.
By the same reasoning, $\pr{D^-_\ell} \le 1/\mm$.
Hence,
for undirected balls,
\[
	\pr{ D^-_i \: \forall \, i \in I }
=
	\pr{ \mm \wr V^-_i \: \forall \, i \in I \midb \mci = I }
\le
	\mm^{-\abs I}.
\]
(The $-$-superscript emphasises that this is for undirected balls.)

We now turn our attention to directed balls.
In this case, $U_i$ and $U'_i$ are both unimodal, but with potentially different modes, if $R \ne R'$.
Instead of direct computation, we compare with the undirected case.
Specifically, if $U_i$ and $U'_i$ have the same sign in the undirected case, then $\abs{V_i} = \abs{U_i - U'_i}$ has the same law as in the directed case.
The choice of sign is independent of everything else; the two have the same sign with probability $\tfrac12$.
Hence,
by conditioning on the specific values of $(U_{i+1}, ..., U_\ell)$ and $(U'_{i+1}, ..., U'_\ell)$,
we obtain
\[
	1/\gamma
\ge
	\pr{ D^-_i \midb D^-_{i+1}, ..., D^-_\ell }
\ge
	\tfrac12 \pr{ D^+_i \midb D^+_{i+1}, ..., D^+_\ell }.
\]

For $\gamma = 2$, note that the probabilities are actually the same: this is because $x - y$ is even if and only if $\abs x - \abs y$ is even, since $x$ and $-x$ have the same parity.

From this we deduce,
for both the undirected and directed cases,
that
\[
	\ex{ \mfgcd^d \mid \mci = I }
\le
	1 + 2^{d - \abs I}
+	\sumt[2M]{\mm=3} \mm^d (2/\mm)^{\abs I}
=
	1 + 2^{d - \abs I}
+	2^d \sumt[2M]{\mm=3} (\mm/2)^{d - \abs I}.
\]
A case-by-case analysis, according to $d - \abs I$, completes the proof.
\end{Proof}

\begin{Proof}[Proof of \cref{res-p3:typdist2:I}]
Recall from \cref{def-p3:typdist2:M*-def} that $\abs{B_k(\mcm_*)} \ge n e^\omega$.
Thus
\[
	\pr{ \mci = \emptyset }
=
	\pr{ W = W' }
=
	\absb{ B_k(\mcm_*) }^{-1}
\le
	n^{-1} e^{-\omega}.
\]

Using the law of $W_I$ given $W_{\setminus I}$ determined in the previous proof,
we have
\[
	\pr{ W_{\setminus I} = W'_{\setminus I} }
=
	\frac
		{ \pr{ W = W' } }
		{ \pr{ W = W' \mid W_{\setminus I} = W'_{\setminus I} } }
=
	\frac
		{ \abs{ B_k(\mcm_*) }^{-1} }
		{ \ex{ \abs{ B_{\abs I}(\mcm_* - \norm{W_{\setminus I}}_1) }^{-1} } }
\le
	\frac
		{ \abs{ B_{\abs I}(\mcm_*) } }
		{ \abs{ B_k(\mcm_*) } }.
\]

It is a standard balls-in-bins combinatorial identity that
\[
	\absb{ B^+_\ell(R) }
=
	\absb{ \brb{ b \in \mbz_+^\ell \midb \sumt[\ell]{1} b_i \le R } }
=
	\binomt{\floor{R}+\ell}{\ell}.
\]
For the undirected case, we can choose a sign for each coordinate. Hence we see that
\[
	\absb{ B^+_\ell(R) }
\le
	\absb{ B^-_\ell(R) }
=
	\absb{ \brb{ b \in \mbz^\ell \midb \sumt[\ell]{1} \abs{b_i} \le R } }
\le
	2^\ell \binomt{\floor{R}+\ell}{\ell}.
\]

Abbreviate $M \cq \mcm_*$ and $\ell \cq \abs I$.
It suffices to consider $I$ with $\ell \le c k$, for an arbitrarily small positive constant $c$.
From \cref{res-p3:typdist2:balls},
we have $M \le 2 \alpha k$.
So
\[
	\absb{ B^\pm_\ell(M) }
\le
	2^\ell \binomt{\floor{M}+\ell}{\ell}
\le
	\rbb{ 2e \rbr{2 \alpha k / \ell + 1} }^\ell
\le
	\rbr{ 8 e \alpha k / \ell }^\ell,
\]
with the last inequality requiring $2 \alpha k / \ell \ge 1$, which holds if $c$ is sufficiently small, as $\ell \le ck$.
Now, for $c$ sufficiently small, the map
\(
	\ell \mapsto \rbr{ 8 e \alpha k / \ell }^\ell
\)
is increasing on $[1, ck]$.
Hence
\[
	\absb{ B^\pm_\ell(M) }
\le
	\rbr{ 8 e \alpha k / \ell }^\ell
\le
	\rbr{ 8 e \alpha / c }^{c k}
\le
	\rbr{ 8 e \alpha / c }^{2 c \lambda \log n}
=
	n^{ 2 c \lambda \log(8 e \alpha / c) }.
\]
Taking $c \to 0$, this exponent tends to $0$.
Hence,
\(
	\abs{ B^\pm_\ell(M) }
=
	n^{\oh1}.
\)
This proves the lemma.
\end{Proof}

\subsection{Relaxing Condition on Minimal Side-Length to $m_*(G) \gg 1$}
\label{sec-p3:typdist2:relax-m*}

For the upper bound, we have been assuming that the minimal side length $m_*(G)$ satisfies $m_*(G) \gg \log \abs G$.
(Recall that the lower bound had no conditions on $m_*(G)$.)
We now describe how to relax this condition to $m_*(G) \gg 1$.
We could go even further, with statements like ``only a small number of $j$ in $G = \oplus_{j=1}^d \: \mbz_{m_j}$ have $m_j \asymp 1$''.
Since we have no reason to believe our other conditions are optimal, we settle for the simpler $m_*(G) \gg 1$.

\smallskip

In this proof we consider both $L_1$ and $L_\infty$ balls.
To distinguish these we use a superscript:
\begin{itemize}[noitemsep, topsep = \smallskipamount, label = \bcdot]
	\item 
	$B_{\ell, 1}(R)$ will be the $L_1$ ball in $\ell$ dimensions of radius $R$;
	
	\item 
	$B_{\ell, \infty}(R)$ will be the $L_\infty$ ball in $\ell$ dimensions of radius $R$.
\end{itemize}
For a set $I \subseteq [k]$, recall that we write $W_I \cq (W_i)_{i \in I}$ and $W_{\setminus I} \cq (W_i)_{i \notin I}$.

\smallskip

We describe the adaptations for undirected graphs.
The adaptations for directed graphs are completely analogous: simply replace appearances of $\mbz^k$ with $\mbz_+^k$ and $\abs{W_i}$ with $W_i$.

\begin{Proof}[Outline of Proof]
The idea behind the proof is intuitive.
Since $R \asymp k$, by symmetry we have $\ex{\abs{W_i}} \le R/k \asymp 1$ for all $i$.
Thus `almost all' the coordinates should be smaller than any diverging function (these coordinates are \textit{good}).
Further, the contribution to the radius $\norm{W}_1$ due to the \textit{bad} coordinates should be small, ie $\oh k$.
Roughly this allows us to replace $k$ with $\widebar k = k\rbr{ 1 - \oh1 }$ and $R$ with $\widebar R = R\rbr{ 1 - \oh1 }$.
Choosing $R \cq \alpha_{k/\log n} k \cdot (1 + 2\eps)$ for $\eps > 0$ then
gives
\[
	\widebar R \ge \alpha_{\widebar k / \log n} \widebar k \cdot (1 + \eps)
\Quad{and hence}
	\abs{B_{\widebar k, 1}(\widebar R)} \gg n.
\]
This was the key element in the proof previously; the remainder of the proof is as before.
\end{Proof}

We now proceed formally and rigorously.

\begin{Proof}[Relaxing Minimal Side-Length Condition]
Let $\eps > 0$ and $\lambda \cq \lim_N k_N / \log n_N$.
Set $R \cq \alpha_\lambda k (1 + 2 \eps)$ and draw $W \sim \Unif\rbr{ B_{k, 1}(R) }$.
Let $\nu$ satisfy $1 \ll \nu \ll m_*(G)$.
For $w \in \mbz^k$, define
\[
	\mcj(w)
\cq
	\brb{ i \in [k] \midb \abs{w_i} \le \nu }.
\]
Call these coordinates \textit{good}.
By Markov's inequality, $\abs{ [k] \setminus \mcj(W) } \lesssim k/\nu = \oh k$ whp as $\ex{\abs{W_i}} \asymp 1$.

As always, we look at two independent realisations $W$ and $W'$.
We then wish to look at coordinates $i \in [k]$ which are good for both $W$ and $W'$, ie in $\widebar \mcj \cq \mcj(W) \cap \mcj(W')$.
We need to make sure that the contribution to the radius from the (abnormally large) \textit{bad} coordinates is not too large.
For $\delta > 0$ and $w \in \mbz^k$,
write $\mcl_\delta(w)$ for the collection of the $\ceil{2 \delta k}$-largest (in absolute value) coordinates of $w$.
We then define typicality in the following way:
for $\delta, \delta' > 0$, set
\[
	\mcw
\cq
	\brb{
		w \in \mbz^k
	\midb
		\norm{w}_1 \le R, \:
		\absb{ [k] \setminus \mcj(w) } \le \delta k, \:
		\norm{w_{\mcl_\delta(w)}}_1 \le \delta ' k
	}.
\]
In particular now, if $w,w' \in \mcw$, then $\norm{w_{\mcj(w) \cap \mcj(w')}}_1 \ge k - 2 \delta' k$.
It is not difficult to see that we can choose $\delta,\delta' = \oh1$ with $\pr{W \in \mcw} = 1 - \oh1$; we give justification at the end of the proof.

\smallskip

Consider now $W, W' \sim^\iid \Unif(B_{k, 1}(R))$.
We have the following conditional law:
\begin{gather*}
	W_{\widebar J}, W'_{\widebar J}
\sim^\iid
	\Unif\rbb{ B_{\widebar k, 1}(\widebar R) \cap B_{\widebar k, \infty}(\nu) }
\Quad{conditional on}
	W_{\setminus \widebar J} = w_{\setminus \widebar J} = W'_{\setminus \widebar J}
\Quad{and}
	\widebar \mcj = \widebar J
\\
	\text{where}
\quad
	\widebar \mcj = \mcj(W) \cap \mcj(W'),
\quad
	\widebar k \cq \abs{\widebar \mcj}
\Quad{and}
	\widebar R \cq R - \norm{w_{\setminus \widebar J}}_1.
\end{gather*}
Write
\(
	\typ \cq \bra{ W, W' \in \mcw }.
\)
On the event $\typ$,
given $\widebar \mcj = \widebar J$ and $(W_{\widebar J}, W'_{\widebar J})$,
we have
\[
	\widebar k \ge k (1 - \delta) = k \rbb{ 1 - \oh1 }
\Quad{and}
	\widebar R \ge R (1 - \delta') = R \rbb{ 1 - \oh 1}.
\]
In particular, we may choose $\eta > 0$ sufficiently small but constant (depending on $\eps$) so that
\[
	\widebar R \ge \alpha_{\lambda(1-\eta)} k (1 - \eta) (1 + \eps)
\Quad{and}
	\widebar k \ge k (1 - \eta),
\Quad{and hence}
	\abs{B_{\widebar k, 1}(\widebar R)} \gg n.
\]
Since typicality holds with probability $1 - \oh1$,
we have
\[
	\absb{ B_{\widebar k, 1}(\widebar R) \cap B_{\widebar k, \infty}(\nu) } \gg n.
\]
The remainder of the proof follows similarly as before.
Formally, we define $\widebar W$ and $\widebar W{}'$ as follows:
\begin{alignat*}{5}
	\widebar W_i &\cq W_i&
&\Quad{and}&
	\widebar W{}'_i &\cq W'_i&
&\Qfor&
	i &\in \widebar \mcj;
\\
	\widebar W_i &\cq 0&
&\Quad{and}&
	\widebar W{}'_i &\cq 0&
&\Qfor&
	i &\notin \widebar \mcj.
\end{alignat*}
Since this is a projection,
\(
	\bra{ W_I = W'_I }
\subseteq
	\bra{ \widebar W_I = \widebar W{}_I }
\)
for any $I \subseteq [k]$.
Now instead of decomposing according to the value (or size) of
\(
	\mci \cq \bra{ i \in [k] \mid W_i \ne 0 },
\)
we use the set
\(
	\widebar \mci \cq \mci \cap \widebar \mcj.
\)
The fact that
\(
	\absb{ B_{\widebar k, 1}(\widebar R) \cap B_{\widebar k, \infty}(\nu) } \gg n
\)
allows all the previous estimates for $\mci$ to follow through for $\widebar \mci$ here.

The last change to mention is the gcd calculations of \cref{res-p3:typdist2:gcd-ex}.
The only property of the distribution of $(W, W')$ required was that each coordinate (while not independent) is unimodal and symmetric about 0, even conditional on $W_I = W'_I$ and $W'_I = w'_I$ for some $I \subseteq [k]$ and $w_I, w'_I \in \mbz^{\abs I}$.
For $(\widebar W, \widebar W{}')$, this property still holds.
Hence the identical argument applies here~too.

\smallskip

It remains to argue that $\pr{W \in \mcw} = 1 - \oh1$ for some $\delta, \delta' = \oh1$.
First,
as noted above,
\(
	\pr{ \abs{[k] \setminus \mcj(W)} > \delta k }
=
	\oh1
\)
by Markov's inequality
and the fact that $\ex{\abs{W_1}} \asymp 1$.
We now show that
\(
	\pr{ \norm{W_{\mcl_\delta(W)}}_1 > \delta' k }
=
	\oh1
\)
for an appropriate choice of $\delta' = \oh1$, to be determined later.
We do this via a union bound over all $\binom{k}{\ceil{2 \delta k}}$ possible values of the set $\mcl_\delta(W)$.

We first bound the above probability by a certain one involving independent random variables.
Consider $Y \cq (Y_1 \xi_1, ..., Y_k \xi_k)$ where $Y_1, ..., Y_k \sim^\iid \Geom_0(\beta)$,
where
	$\beta$ is picked so that $\ex{ \norm{Y}_1 } = k \rbr{ 1/\beta - 1 } = 2R$
and
	$\xi_1, ..., \xi_k \sim^\iid \Unif(\bra{-1, +1})$;
recall that $R = \alpha_\lambda k (1 + 2\eps)$.
Here, the Geometric distribution $\Geom_0$ is supported on $\bra{0, 1, 2, ...}$.
We have $\ex{ e^{Y_1} } \asymp 1$ since $k \asymp R$.
Given $\ell \in [0, R]$, the law of $Y$ conditioned on $\norm{Y}_1 = \ell$ is the same as that of $W$ conditioned on $\norm{W}_1 = \ell$.

It is straightforward that the law of $\abs Y \cq (Y_1, ..., Y_k)$ conditioned on $\norm{Y}_1 = \ell$ is stochastically increasing in $\ell$. Indeed, one can couple $\abs Y = (Y_1, ..., Y_k)$ conditioned on $\norm{Y}_1 = \ell$ with $\abs Y = (Y_1, ..., Y_k)$ conditioned on $\norm{Y}_1 = \ell+1$ by first sampling the former, then picking a random coordinate and increasing it by $+1$ to obtained the latter.
It follows that the law of $\abs Y$ given $\norm{Y}_1 > R$ stochastically dominates the law of $(\abs{W_1}, ..., \abs{W_k})$ which is a mixture of the laws of $\abs Y = (Y_1, ..., Y_k)$ conditioned on $\norm{Y}_1 = \ell$ for different $\ell \in [0, R$].

Since $\bra{ \norm{Y_{\mathcal{L}_{\delta}(Y)}}_1 > \delta' k}$ is a monotone increasing event wrt $\abs Y$, we have
\[
	\pr{ \norm{ W_{\mcl_\delta(W)} }_1 > \delta' k }
\le
	\pr{ \norm{ Y_{\mcl_\delta(Y)} }_1 > \delta' k \midb \norm{Y}_1 > R }.
\]
Next, note that $\pr{ \norm{Y}_1 > R } \asymp 1$, by our choice of $\beta$, and
\[
	\pr{ \sumt[\ceil{2 k \delta}]{i=1} Y_i \ge \delta' k }
\le
	e^{-\delta' k} \ex{ \expb{\sumt[\ceil{2 k \delta}]{i=1} Y_i} }
=
	e^{-\delta' k} \ex{ e^{Y_1} }^{\ceil{2 \delta k}}.
\]
These together with Bayes's rule and the union bound then gives
\[
	\pr{ \norm{ W_{\mcl_\delta(W)} }_1 > \delta' k }
\lesssim
	\pr{ \norm{ Y_{\mcl_\delta(Y)} }_1 > \delta' k }
\le
	\binomt{k}{\ceil{2 \delta k}}
	e^{-\delta' k}
	\ex{ e^{Y_1} }^{\ceil{2 \delta k}}.
\]
Next, there exists some universal constant $C > 0$ such that
\[
	\binomt{k}{\ceil{2 \delta k}}
\le
	\expb{ C k \delta \log(1/\delta) }.
\]
We then obtain
\(
	\pr{ \norm{ W_{\mcl_\delta(W)} }_1 > \delta' k }
=
	\oh1,
\)
as required,
by taking
\[
	\delta'
\cq
	\delta \rbb{ 2 C \log(1/\delta) + 4 \log \ext{ e^{Y_1} } }
\asymp
	\delta \log(1/\delta)
=
	\oh1.
\qedhere
\]
\end{Proof}

\begin{rmkt*}
	We believe that the typical distance should concentrate if $k \asymp \log \abs G$ and $k - d \gg 1$ without any condition like that on $m_*(G)$.
	However, without any such condition, we do have reason to believe that the \emph{value} at which this concentration happens should depend on more than just $k$ and $\abs G$---the algebraic structure of $G$ should be important.
	This exact phenomenon occurs when studying the mixing time of the random walk on the Cayley graph.
	See
	\cite[Theorem~A]{HOt:rcg:abe:cutoff},
	in particular contrasting the case
		$k \asymp \log \abs G \asymp d(G)$
	with
		$1 \ll k \lesssim \log \abs G$ and $d(G) \ll \log \abs G$.
\end{rmkt*}

\subsection{Typical Distances for $L_\qq$-Type Graph Distances}
\label{sec-p3:typdistq}

Graph distances in Cayley graphs have some special properties.
Consider a collection $z = [z_1, ..., z_k]$ of generators and distances in the Cayley graph $G(z)$.
For a path $\rho$ in $G(z)$,
for each $i \in [k]$,
write
	$\rho_{i,+}$ for the number of times $z_i$ is used,
	$\rho_{i,-}$ for the number of times $z_i^{-1}$ is used (if in the undirected case, otherwise $\rho_{i,-} \cq 0$)
and
	$\rho_i \cq \rho_{i,+} - \rho_{i,-}$.
The path connects the identity with $\rho \bcdot z$.
Then the length, in the usual graph distance, of $\rho$ is $\norm{\rho}_1 \cq \sumt[k]{1} \rbr{ \rho_{i,+} + \rho_{i,-} }$.

For any $\qq \in [1,\infty)$,
define the \textit{$L_\qq$ graph distance} of $\rho$ by
\(
	\norm{\rho}_\qq^\qq
\cq
	\sumt[k]{1} \rbr{ \rho_{i,+}^\qq + \rho_{i,-}^\qq }.
\)
For the \textit{$L_\infty$ graph distance},
define
\(
	\norm{\rho}_\infty
\cq
	\maxt{i} \bra{ \rho_{i,+} + \rho_{i,-}}.
\)
(The usual graph distance is given by $\qq = 1$.)

For Abelian groups, clearly for any $\qq \in [1,\infty)$ an \textit{$L_\qq$ geodesic}, ie a path of minimal $L_\qq$ weight, will only use either $z_i$ or $z_i^{-1}$, not both (since the terms in the product can be reordered), ie $\rho_{i,+} \rho_{i,-} = 0$ for all $i$.
Thus $\norm{\rho}_\qq^\qq = \sumt[k]{1} \abs{\rho_i}^\qq$.
Similarly, any $L_\infty$ geodesic $\rho$ can be adjusted into a new path $\rho'$ with $\rho \bcdot z = \rho' \bcdot z$ and $\norm{\rho}_\infty = \norm{\rho'}_\infty$ satisfying $\rho'_{i,+} \rho'_{i,-} = 0$ for all $i$.

\smallskip

We define the \textit{$L_\qq$ typical distance} $\mcd_{G(z),\qq}(\cdot)$ analogously to $\mcd_{G(z)}(\cdot)$, ie the $\qq = 1$ case.

\addtocounter{hyp}{-1} \refstepcounter{hyp}
\begin{customhyp}{\thehyp$'$}
\label{hyp-p3:typdistq}
	The sequence
		$(k_N, G_N)_\Ninn$
	and
		$\qq \in [1,\infty]$
	jointly satisfy \textit{Hypothesis~\ref{hyp-p3:typdistq}} if the following conditions hold (defining $k^{1/\infty} \cq 1$ for $k \in \mbn$):%
	\begin{gather*}
		\LIM{\Ninf} k_N = \infty,
	\quad
		\LIM{\Ninf} k_N / \log \abs{G_N} = 0
	\Qand
		\LIM{\Ninf} k_N^{1/\qq} \abs{G_N}^{1/k_N} / m_*(G_N) = 0;
	\\
		\text{if\quad $\qq \in (1,\infty)$\quad then additionally\quad $k_N \le \log \abs{G_N} / \log\log \abs{G_N}$ for all $\Ninn$};
	\\
		\LIMSUP{\Ninf} \frac{d(G_N)}{k_N}
	<
		\begin{cases}
			1			&\text{for undirected graphs}, \\
			\tfrac12	&\text{for directed graphs},
		\end{cases}
	\end{gather*}
	recalling that $d(H)$ is the minimal size of a generating set for a group $H$.
\end{customhyp}


Finally, we set up a little more notation.
Let $\Gamma(\cdot)$ denote the Gamma function.
Let
\begin{gather*}
	C^-_\qq \cq 2 \, \Gamma(1/\qq+1) (\qq e)^{1/\qq},
\quad
	C_\qq^+ \cq \tfrac12 C^-_\qq
\Qand
	\mfd^\pm_\qq(k,n) \cq k^{1/\qq} n^{1/k} / C^\pm_\qq,
\end{gather*}
where the case $\qq = \infty$ is to be interpreted as the limit $\qq \to \infty$;
eg,
\(
	C^-_\infty
=
	2
\)
and
\(
	\mfd^+_\infty(k,n)
=
	n^{1/k}.
\)
When these are sequences $(k_N, \abs{G_N})_\Ninn$,
for $\Ninn$ and $q \in [1,\infty]$,
write $\mfd^\pm_{N,\qq} \cq \mfd^\pm_\qq(k_N,\abs{G_N})$.

Similarly,
for a sequence $(G_N)_\Ninn$ of finite groups with corresponding multisubsets $(Z_{(N)})_\Ninn$ of sizes $(k_N)_\Ninn$,
for $\Ninn$, $\beta \in [0,1]$ and $\qq \in [1,\infty]$,
define
\(
	\mcd_{N,\qq}^\pm \cq \mcd_{G_N^\pm(Z_{(N)})}(\beta).
\)

Using an extension of the methodology from this section (\S\ref{sec-p3:typdist2}), along with analysis of $L_\qq$ lattice balls, we can prove the following theorem.
We have already considered $\qq = 1$ and $k \asymp \log \abs G$.

\begin{thm}
\label{res-p3:typdistq:res}
	Let $(k_N)_\Ninn$ be a sequence of positive integers and $(G_N)_\Ninn$ a sequence of finite, Abelian groups;
	for each $\Ninn$, define $Z_{(N)} \cq [Z_1, ..., Z_{k_N}]$ by drawing $Z_1, ..., Z_{k_N} \sim^\iid \Unif(G_N)$.
	
	Suppose that $(k_N, G_N)_\Ninn$ satisfies Hypothesis~\ref{hyp-p3:typdistq}.
	Then,
	for all $\beta \in (0,1)$,
	we~have
	\[
		\mcd^\pm_{N,\qq}(\beta) /\mfd^\pm_{N,\qq}
	\to^\mbp
		1
	\Quad{(in probability)}
		\asinf N.
	\]
	Moreover, the implicit lower bound holds for all choices of generators and for all Abelian groups, only requiring the conditions in Hypothesis~\ref{hyp-p3:typdistq} which depend only on $(k_N, \abs{G_N})_\Ninn$ and $q$.
\end{thm}

The arguments used to prove this theorem really are analogous to those used in this section (\S\ref{sec-p3:typdist2}).
The only real difference is that we have to look at lattice balls under an $L_\qq$ norm and in dimension $1 \ll k \ll \log n$, rather than $L_1$ and $k \asymp \log n$.
Other than this, the remainder of the analysis, in particular the reduction to a gcd and the consideration of the set $\mci$ of non-zero coordinates of $W$, is exactly the same. (Now $W$ is uniform on an $L_\qq$ ball of appropriate radius.)
We do not give the details here; they can be found in
\cite[\S 7]{HOt:rcg:abe:extra}.

We remark that $k \asymp \log \abs G$ is not covered when $q \in (1, \infty)$, ie $q \notin \bra{1, \infty}$.
This, in essence, is because we need estimates on the volume of $L_\qq$ balls.
These can be estimated very precisely when $\qq \in \bra{1, \infty}$, but our estimates are less precise otherwise.
See
\cite[Lemma~7.2b]{HOt:rcg:abe:extra}
for specifics.

\section{Typical Distance: $k \gg \log \abs G$}
\label{sec-p3:typdist3}

This section focusses on concentration of distances from the identity in the random Cayley graph of an Abelian group when $k \gg \log \abs G$.
(The previous sections dealt with $1 \ll k \lesssim \log \abs G$.)
The main result of the section is \cref{res-p3:typdist3:res}. 

\smallskip

The outline of this section is as follows:
\begin{itemize}[noitemsep, topsep = 0pt, label = \bcdot]
	\item 
	\S\ref{sec-p3:typdist3:res} states precisely the main theorem of the section;
	
	\item 
	\S\ref{sec-p3:typdist3:outline} outlines the argument;
	
	\item 
	\S\ref{sec-p3:typdist3:balls} gives some crucial estimates on the size of lattice balls;
	
	\item 
	\S\ref{sec-p3:typdist3:lower} is devoted to the lower bound;
	
	\item 
	\S\ref{sec-p3:typdist3:upper} is devoted to the upper bound.
\end{itemize}

\subsection{Precise Statement and Remarks}
\label{sec-p3:typdist3:res}

To start the section, we recall the typical distance statistic.

\begin{defn}
\label{def-p3:typdist3:def}
Let $H$ be a graph and fix a vertex $0 \in H$.
For $r \in \mbn$,
write $\mcb_H(r)$ for the $r$-ball in the graph $H$, ie
\(
	\mcb_H(r)
\cq
	\bra{ h \in H \mid d_H(0, h) \le r },
\)
where $d_H$ is the graph distance in $H$.
Define
\[
	\mcd_H(\beta)
\cq
	\min\brb{ r \ge 0 \midb \abs{ \mcb_H(r) } \ge \beta \abs H }
\Qfor
	\beta \in (0,1).
\]

When considering sequences $(k_N, G_N)_\Ninn$ of integers and Abelian groups,
abbreviate
\[
	\mcd_N(\beta)
\cq
	\mcd_{G_N([Z_1, ..., Z_{k_N}])}(\beta)
\Qwhere
	Z_1, ..., Z_{k_N} \sim^\iid \Unif(G_N).
\]
Finally, considering such sequences, we define the candidate radius for the typical distance:
\[
	\widebar \mfd_N
\cq
	\tfrac{\rho_N}{\rho_N-1} \log \abs{G_N} / \log k_N
\Qwhere
	\rho_N \cq \log k_N / \log \log \abs{G_N}
\Quad{for each}
	\Ninn.
\]
As always, if we write $\mcd_N$, then this is either $\mcd_N^+$ or $\mcd_N^-$ according to context.
Up to subleading order, the typical distance will be the same for the undirected graphs as for the directed graphs.
\end{defn}

We show that, whp over the graph (ie choice of $Z$), this statistic concentrates.
Here we consider $k \gg \log \abs G$.
The result holds for all Abelian groups; in fact, the implicit upper bound is valid for all groups.
Further, the typical distance concentrates at a distances which depends only on $k$ and $\abs G$.
This is in agreement with the spirit of the Aldous--Diaconis conjecture.

\begin{hyp}
\label{hyp-p3:typdist3}
	The sequence $(k_N, n_N)_\Ninn$ satisfies \textit{\cref{hyp-p3:typdist3}} if
	\[
		\LIMINF{\Ninf}
		\frac{ k_N }{ \log n_N }
	=
		\infty
	\Qand
		\LIMINF{\Ninf}
		\frac{ \log k_N }{ \log n_N }
	=
		0.
	\]
\end{hyp}

\begin{thm}
\label{res-p3:typdist3:res}
	Let $(k_N)_\Ninn$ be a sequence of positive integers and $(G_N)_\Ninn$ a sequence of finite, Abelian groups;
	for each $\Ninn$, define $Z_{(N)} \cq [Z_1, ..., Z_{k_N}]$ by drawing $Z_1, ..., Z_{k_N} \sim^\iid \Unif(G_N)$.
	
	Suppose that $(k_N, \abs{G_N})_\Ninn$ satisfies \cref{hyp-p3:typdist3}.
	Then,
	for all $\beta \in (0,1)$,
	we have
	\[
		\mcd^\pm_N(\beta) / \widebar \mfd_N
	\to^\mbp
		1
	\Quad{(in probability)}
		\asinf N.
	\]
	Moreover,
		the implicit lower bound holds deterministically, ie for all choices of generators,
	and
		the implicit upper bound holds for all groups, not just Abelian groups.
\end{thm}

As always, for ease of presentation, in the proof we drop the $N$-subscripts.

\subsection{Outline of Proof}
\label{sec-p3:typdist3:outline}

When $k \gg \log \abs G$,
one can see that
the typical distance statistic $\mcd$ must satisfy $\mcd \ll k$.
By symmetry, the expected number of times a generator is used when drawing from a ball is $\oh1$.

An approximation is then that each generator is either chosen once or not at all.
If $R$ are chosen, then there are precisely $\binom kR$ ways of doing this.
We choose $R$ with $\binom kR \approx \abs G$.

\subsection{Estimates on Sizes of Balls in $\mbz^k$}
\label{sec-p3:typdist3:balls}

We consider balls and spheres in the $L_1$ and $L_\infty$ senses:
	write $B_{k,1}(\cdot)$, respectively $S_{k,1}(\cdot)$, for the $L_1$ ball, respectively sphere, in $\mbz^k$;
	write $B_{k,\infty}(1)$ for the $L_\infty$ unit ball in $\mbz^k$.

\begin{lem}
\label{res-p3:typdist3:balls}
	For all $R \ge 0$,
	we have
	\[
		\abs{B_{k,1}^\pm(R)} \le 2^R \binomt{\floor{R}+k}{\floor{R}}
	\Quad{and}
		\absb{ S_{k,1}^\pm(R) \cap B_{k,\infty}^\pm(1) } \ge \binomt{k}{\floor{R}}.
	\]
	Furthermore, if $R \ll k$, then both
	\[
		2^R \binomt{\floor{R}+k}{\floor{R}}
	=
		\expb{ R \log(k/R) \cdot \rbr{ 1 + \oh1 } }
	\Qand
		\binomt{k}{\floor{R}}
	=
		\expb{ R \log(k/R) \cdot \rbr{ 1 + \oh1 } },
	\]
	with different $\oh1$ terms, naturally.
	In particular,
	if $k = (\log n)^\rho \gg \log n$ and $\eps > 0$ is constant,~then
	\[
		\absb{ S_{k,1}^\pm\rbb{ (1 + \eps) \tfrac{\rho}{\rho-1} \log_k n } \cap B_{k,\infty}(1) }
	\gg
		n.
	\]
\end{lem}

\begin{Proof}
In the first display,
	the upper bound is proved in
	\cite[Lemma~E.2a]{HOt:rcg:supp};
	the lower bound is the usual formula for the number of subsets of $[k]$ of size $R$.
The second display is a simple application of Stirling's approximation and asymptotics of the binary entropy function.
The final display follows by combining the previous two and performing a simple calculation.
Indeed,
take
\[
	R
\cq
	(1 + \eps) \tfrac{\rho}{\rho-1} \log_k n
=
	(1 + \eps) \log n / \log(k/\log n)
\ll
	k.
\]
Thus, applying the first lower bound followed by the second asymptotic equality,
\[
&	\absb{ S_{k,1}^\pm\rbr{ R } \cap B_{k,\infty}(1) }
\ge
	\binomt{k}{\floor{R}}
\\&\qquad
=
	\expb{ R \log(k/R) \cdot \rbr{ 1 + \oh1 } }
\\&\qquad
=
	\expb{ (1 + \eps) \log n \cdot  \rbr{ 1 + \oh1 } }
\ge
	n^{1 + \eps/2}
\gg
	n.
\qedhere
\]
\end{Proof}

\subsection{Lower Bound on Typical Distance}
\label{sec-p3:typdist3:lower}

From the results in \S\ref{sec-p3:typdist3:balls}, it is straightforward to deduce the lower bound in \cref{res-p3:typdist3:res}.

\begin{Proof}[Proof of Lower Bound in \cref{res-p3:typdist3:res}]
Let $\xi \in (0, 1)$ and set $R \cq \widebar \mfd (1 - \xi)$,
recalling that
\[
	\widebar \mfd = \tfrac{\rho}{\rho-1} \log n / \log k
\Qwhere
	\rho = \log k / \log \log n,
\Quad{ie}
	k = (\log n)^\rho.
\]
Since the underlying group is Abelian,
applying \cref{res-p3:typdist3:balls},
a simple calculation gives
\[
	\abs{ \mcb_k\rbr{ R } }
\le
	\abs{ B_{k,1}\rbr{ R } }
\le
	\expb{ \widebar \mfd \log(k / \widebar \mfd) \cdot (1 - \tfrac12 \xi) }
\ll
	n.
\]
Hence,
for all $\beta \in (0,1)$ and all $Z$,
the typical distance satisfies
$\mcd_k(\beta) \ge R =  \widebar \mfd (1 - \xi)$.
	%
\end{Proof}

\subsection{Upper Bound on Typical Distance}
\label{sec-p3:typdist3:upper}

\cref{res-p3:typdist3:balls} gives a quantitative sense in which
\(
	\abs{B_{k,1}(R)}
\approx
	\absb{ S_{k,1}(R) \cap B_{k,\infty}(1) } \ge \binomt{k}{\floor{R}};
\)
informally, this means that we do not really lose any volume by restricting to the sphere and requiring that each generator is used at most once.
We show the upper bound for arbitrary groups.

\begin{Proof}[Proof of Upper Bound in \cref{res-p3:typdist3:res}]
Let $\xi > 0$ and set $R \cq \widebar \mfd (1 + \xi)$.
Draw $\WW,\WW' \sim^\iid \Unif(S_{k,1}(R) \cap B_{k,\infty}(1))$.
Define $S \cq Z_1^{W_1} \cdots Z_k^{W_k}$ and $S'$ similarly.
We show that $S$ is well-mixed whp (this time in the $L_2$ sense) to deduce the upper bound.
By the standard $L_2$ calculation,
\[
	\ex{ \norm{ \pr[G_k]{ S \in \cdot } - \Unif(G) }_2^2 }
=
	n \, \pr{ S = S' } - 1.
\]
If $\WW \ne \WW'$, then there exists an $i \in [k]$ with $\WW_i = 1$ and $\WW'_i = 0$ or vice versa.
Then, by uniformity and independence of the generators,
$S' S^{-1} \sim \Unif(G)$
for \emph{all} (not just Abelian) groups.
Thus,
\[
	n \, \pr{ S = S' } - 1
\le
	n \, \pr{ \WW = \WW' }
=
	n \, \absb{ S_{k,1}(R) \cap B_{k,\infty}(1) }^{-1}
\ll
	1,
\]
using \cref{res-p3:typdist3:balls}.
This means that $\WW \bcdot Z$ is well-mixed in the $L_2$ sense, so has support $n - \oh n$.
Thus, $n - \oh n$ of the vertices are within distance $R = \widebar \mfd (1 + \xi)$ of the identity.
\end{Proof}

\begin{rmkt*}
	This upper bound, ie on typical distance with $k \gg \log \abs G$, can be easily deduced from mixing results proved in the '90s.
	Specifically, it was shown by \textcite[Theorem~1]{DH:enumeration-rws} that the mixing time for the usual random walk is upper bounded by $\tfrac{\rho}{\rho-1} \log_k \abs G$ for any group; \textcite[Theorems~1 and~2]{R:random-random-walks} subsequently gave a simpler proof, using an argument not that dissimilar from our proof above.
	The lower bound does not follow from mixing results,~though.
	
	There are a few reasons for including the proof above.
		Foremost is that we use the same argument in \S\ref{sec-p3:diam:univ} to obtain universal bounds for $k$ with $k - \log_2 \abs G \asymp k$, not just $k \gg \log \abs G$.
		Additionally, we need to do most of the work for the lower bound anyway, and it demonstrates how easily our method adapts to this new regime.
\end{rmkt*}

\section{Typical Distance for Nilpotent Groups}
\label{sec-p3:nil}

Recall the definition of \textit{nilpotent} and the corresponding distance-based definitions in Definition~\ref{def-p3:intro:nil}.

\subsection{Formal Statements for Dominance of Abelianisation}
\label{sec-p3:nil:thm}

	%
If $Z_1, ..., Z_k \sim^\iid \Unif(G)$, then $[G,G] Z_1, ..., [G,G] Z_k \sim^\iid \Unif(G^\ab)$.
In the case that $k \ge d(G)$ is constant, the asymptotic law of $\diam_{Z \cup Z^{-1}}(G^\ab)/\abs{G^\ab}^{1/k}$ as $|G| \to \infty$ was determined by \textcite{SZ:diam-cayley},
where $Z \cup Z^{-1} \cq [Z_1, Z_1^{-1}, ..., Z_k, Z_k^{-1}]$.
\textcite{EbP:cayley-diam-nil} proved that
the asymptotic laws of
\(
	\diam_{Z \cup Z^{-1}}(G) / \abs{G^\ab}^{1/k}
\)
and
\(
	\diam_{Z \cup Z^{-1}}(G^\ab) / \abs{G^\ab}^{1/k}
\)
are the same
when $k \ge d(G)$ and $\ell$ are fixed.
\cref{res-p3:nil:thm} below forms a quantitative version of the results of \cite{EbP:cayley-diam-nil}.

We are primarily concerned with the case $k \gg 1$. This means that we must keep track on the dependence of certain constants on $|S|$ and $\ell$.
\cref{eq-p3:nil:cor:dominance} below is stated for a general symmetric set of generators $S$, but we are particularly interested in the case that $S = Z \cup Z^{-1}$, where $Z_1, ..., Z_k \sim^\iid \Unif(G)$, as above.
The asymptotics of $\mcd_{Z \cup Z^{-1}}(G^\ab, \beta)$ can be determined using~\cref{res-p3:intro:typdist}.
	%

\medskip

The following theorem
is to appear in \textcite{H+:nil-to-abe}.
As stated before, its argument builds on arguments of \textcite{EbP:cayley-diam-nil}, who built on those of \textcite{BT:moderate-growth}.

\begin{thm}[\cite{H+:nil-to-abe}]
\label{res-p3:nil:thm}
	Let $G$ be a finite nilpotent group of step $\ell$ and rank $d$.
	Let $k \in \mbn$ and $s_1, ..., s_k \in G$.
	Write $S \cq [s_1, s_1^{-1}, ..., s_k, s_k^{-1}]$.
	Then, for all $\beta \in (0, 1]$,
	we have
	\begin{gather*}
		\mcd_S(G^\ab,\beta) \le \mcd_S(G,\beta) \le \mcd_S(G^\ab,\beta)+ \diam_S(G_2)
	\Qand
	\label{eq-p3:nil:thm:dominance}
	\nt
	\\
		\diam_S([G, G])
	\le
		\sumt[\ell]{i=2}
		\diam_S(G_i / G_{i+1})
	\le
		\sumt[\ell]{2}
		2^{i+3} k^i
		\rbb{ \ceil{ \diam_S(G^\ab) / k }^{1/i} + 2^4(i^2 + 4) }.
	\label{eq-p3:nil:thm:diam-Gcom}
	\nt
	\end{gather*}
\end{thm}


\cref{res-p3:intro:nil} follows as a corollary of this.
We repeat the statement below for convenience.

\begin{customthm}{\ref{res-p3:intro:nil}}
	Let $G$ be a finite nilpotent group of step $\ell$ and rank $d$.
	Let $k \in \mbn$ and let $s_1, ..., s_k \in G$.
	Let $S \cq [s_1, s_1^{-1}, ..., s_k, s_k^{-1}]$ be a symmetric multisubset of $G$.
	Then, for all $\beta \in (0, 1]$,
	we have
	\[
		0
	\le
		\mcd_S(G,\beta) - \mcd_S(G^\ab,\beta)
	\lesssim
		\diam_S(G^\ab)^{3/4}
	\le
		\rbb{ 3 \mcd_S(G^\ab,\beta)/\beta }^{3/4}.
	\label{eq-p3:nil:cor:dominance}
	\nt
	\]
	if
	\(
		d
	\le
		k
	\le
		\tfrac1{16} \ell^{-1} d^{-\ell} \log \abs G / \log \log \abs G.
	\)
	Further, for all $\beta \in (0,1/2)$, we have
	\[
		\mcd_S(G^\ab,1-\beta) - \mcd_S(G^\ab,\beta)
	\le
		2 \sqrt{\beta^{-1}\trel\rbb{ \Cay(G^\ab,S^\ab) }},
	\label{eq-p3:nil:cor:concentration}
	\nt
	\]
	where $\trel(\Cay(G^\ab,S^\ab))$ is the relaxation time of the simple random walk on $\Cay(G^\ab,S^\ab)$.
\end{customthm}


\begin{Proof}[Proof of \cref{res-p3:intro:nil} Given \cref{res-p3:nil:thm}]
The first inequality in \cref{eq-p3:nil:thm:dominance} implies the first in \cref{eq-p3:nil:cor:dominance}.

For the second inequality in \cref{eq-p3:nil:cor:dominance},
the second in \cref{eq-p3:nil:thm:dominance} implies that
it suffices to show that
\[
	\diam_S(G_2)
\lesssim
	\diam_S(G^\ab)^{3/4}
\Quad{whenever}
	d
\le
	k
\le
	\frac{\log |G|}{16\ell d^{\ell}\log \log |G|}.
\label{eq-p3:nil:cor:pf:1}
\nt
\]
One can show that
\(
	|G_i/G_{i+1}| \le \abs{G^\ab}^{d^{i-1}}
\)
using standard arguments à la \cite[Equation~(2.3)]{EbP:cayley-diam-nil}.
A precise proof will appear in \cite{H+:nil-to-abe}.
It follows that
\(
	|G|
=
	\prod_{i=1}^{\ell} |G_i/G_{i+1}|
\le
	\abs{G^\ab}^{2d^{\ell}}.
\)

Assume that
\(
	d
\le
	k
\le
	\tfrac1{16} \ell^{-1} d^{-\ell} \log \abs G / \log \log \abs G.
\)
The general lower bound on the diameter from \cref{res-p3:intro:typdist} for Cayley graphs of Abelian groups,
we have
\[
	\diam_S(G^\ab)
\gtrsim
	k \abs{G^\ab}^{1/k}
\ge
	|G|^{1/(2 d^\ell k)}  \ge (\log |G)|)^{8\ell}.
\label{e:lowerboundonDGab}
\]
Finally,
this
together with \cref{eq-p3:nil:thm:diam-Gcom} gives
\[
	\diam_S(G_2)
&
\le
	\sumt[\ell]{i=2}
	2^{i+3} k^i
	\rbb{ \lceil \diam_S(G^\ab)/k\rceil^{1/i}+2^7 i^2 } 
\\&
\lesssim
	(2k)^\ell
	\rbb{ \ell^2+  \diam_S(G^\ab)^{1/2} }
\le
	\diam_S(G^\ab)^{3/4},
\]
where in the last inequality we used the fact that $d \le k \le \frac{\log |G|}{16\ell d^{\ell}\log \log |G|}$.
This proves \cref{eq-p3:nil:cor:pf:1}.

The last inequality in \cref{eq-p3:nil:cor:dominance} follows from the fact that
\[
	\diam_S(G^\ab) \le 2 \mcd_S(G^\ab,1/2)+1
\Qfor
	\beta \ge 1/2
\]
and we can fit $\lfloor \diam_S(G^\ab)/(2r+1) \rfloor$ disjoint balls of radius $r = \mcd_S(G^\ab,\beta)$ in $\Cay(G^\ab,S^\ab)$.
	%

	%
We turn to \cref{eq-p3:nil:cor:concentration}.
Let $H$ be an Abelian group and $S'$ a symmetric set of generators.
We show
\[
	\mcd_{S'}(H,1-\beta)-\mcd_{S'}(H,\beta)
\le
	2 \sqrt{\beta^{-1}\trel\rbr{\Cay(H,S')}}
\Qforall
	\beta \in (0,1/2).
\label{eq-p3:nil:cor:pf:2}
\nt
\]
This implies \cref{eq-p3:nil:cor:concentration} by taking $H = G^\ab$ and $S' = S^\ab$.
It remains to prove \cref{eq-p3:nil:cor:pf:2}.

Let $f:H \to \mbz_+$ be given by $f(h) \cq \dist_{S'}(\id,h)$. Then, $f$ is 1-Lipschitz wrt $\dist_{S'}(\cdot,\cdot\cdot)$.
Hence,
\[
	\mce(f,f)
\cq
	\tfrac12
	\sumt{h \in H,s \in S'}
	\tfrac{1}{|H||S'|} \rbr{ f(h)-f(hs) }^2
\le
	1.
\]
Let $P$ be the transition matrix of simple random walk on $\Cay(H,S')$. Let $\pi$ be the uniform distribution on $H$.
For $g,g' \in \mbr^H$, let
\[
	\langle g,g' \rangle_\pi
\cq
	\ext[\pi]{gg'}
\Qwhere
	\ext[\pi]{u}
\cq
	\sumt{h \in H}
	\pi(h) u(h)
\Qfor
	u \in \mbr^G.
\]
We have $\mce(f,f) = \langle (I-P)f,f \rangle_{\pi}$ by \cite[Lemma~13.6]{LPW:markov-mixing}. By the Courant--Fischer characterization of the eigenvalues of $I-P$, explained in \cite[Remark~13.8]{LPW:markov-mixing}, we have
\[
	\VAR[\pi]{f}
\le
	\trel\rbb{ \Cay(H,S') } \mce(f,f)
\le
	\trel\rbb{ \Cay(H,S') }
\Qwhere
	\VAR[\pi]{f}
\cq
	\ext[\pi]{f^2} - \rbb{ \ext[\pi]{f} }^2.
\]
An application of Chebyshev's inequality concludes the proof of \cref{eq-p3:nil:cor:pf:1}, and hence of \cref{res-p3:intro:nil}.
\end{Proof}

\subsection{Intuition Behind Dominance of Abelianisation}
\label{sec-p3:nil:intuition}


Theorems~\ref{res-p3:intro:nil} and~\ref{res-p3:nil:thm} explain some key differences to the geometry of the graph in the non-Abelian versus Abelian set-ups in a formal manner.
The current subsection gives intuition for why this domination of the Abelianisation occurs, including where the Abelian proof breaks down.
The intuition here addresses non-Abelian groups more generally, in a wider sense than \cref{res-p3:nil:thm}.



Draw $W, W' \sim^\iid \Unif(B_k(R))$ for some $R$ as before.
Define $S$ and $S'$ according to $W$ and $W'$, respectively.
In the Abelian set-up, $W = W'$ implies $S = S'$, but this is not the case generally.~%
Then,
\[
	n \, \pr{ S = S' } - 1
&
=
	n \, \pr{ S = S' \mid W \ne W' } \pr{ W \ne W' }
-	1
\\&\qquad
+	n \, \pr{ S = S' \mid W = W' }   \pr{ W = W' }.
\]
We should really be doing this conditioned on typicality, but we omit this here for simplicity.
The handling of the first term is similar in the non-Abelian and Abelian cases:
	roughly, $S' S^{-1} \sim \Unif(G)$ given $W \ne W'$, under some conditions;
	this balances the $n$-factor and cancels with the $-1$.

It is the last term which depends on the Abelian property:
	$\pr[Z]{S = S' \mid W = W'} = 1$ if the group is Abelian;
	the probability of $W = W'$ must then balance the $n$-factor.
This is why we chose the balls precisely so that the volume was slightly larger than $n$.
If the group is non-Abelian, then $\pr{S = S' \mid W = W'}$ may be much smaller than $1$,
meaning that $\pr{W = W'}$ need not be so small.

\medskip

This exactly the situation in our companion paper \cite{HOt:rcg:matrix} in which we study certain non-Abelian matrix groups.
In the simplest case, we consider upper-triangular, $3 \times 3$ matrices with $1$s on the diagonal and all strictly-superdiagonal entries in $\mbz_p$, with $p$ prime.
One has
\[
\begin{pmatrix}
	1 & a & c \\
	0 & 1 & b \\
	0 & 0 & 1
\end{pmatrix}
\cdot
\begin{pmatrix}
	1 & a' & c' \\
	0 & 1  & b' \\
	0 & 0  & 1
\end{pmatrix}
=
\begin{pmatrix}
	1 & a+a' & c + c' + ab' + a'b \\
	0 & 1    & b + b' \\
	0 & 0    & 1
\end{pmatrix}
\Qforall
	a,a',b,b',c,c' \in \mbz.
\]
This group satisfies an appropriate version of $\pr{ S = S' \mid W \ne W' } \approx 1/n = 1/p^3$.
It is easy to see that $W = W'$ implies that the two immediately-superdiagonal terms in $S$ and $S'$ are equal---%
	indeed, the Abelianisation $G^\ab = G/[G,G]$ precisely corresponds to these.
One can prove that
\[
	\pr{ S = S' \mid W = W' }
\approx
	1/p
=
	p^2 / n
\]
when $1 \ll k \ll \log n \asymp \log p$,
under appropriate typicality conditions.
This means that we need only choose the radius $R$ so that the ball has volume $p^2$, not $n = p^3$ as in the Abelian case.
This is precisely the change from $\abs G = p^3$ to $\abs{G^\ab} = p^2$ predicted by \cref{res-p3:intro:nil}.
More details on this can be found in our companion paper \cite{HOt:rcg:matrix}; see, in particular,
\cite[\S 5 and Theorem~5.1]{HOt:rcg:matrix}.

\section{Diameter}
\label{sec-p3:diam}

In this section we consider the diameter of the random Cayley graph.
Our analysis is separated into two distinct sections.

\begin{itemize}[noitemsep, topsep = 0pt, label = \bcdot]
	\item [\S\ref{sec-p3:diam:conc}]
	We show that the diameter concentrates for $k \gtrsim \log \abs G$, and that the value at which it concentrates is the same as for typical distance.
	
	\item [\S\ref{sec-p3:diam:univ}]
	We show, for $k$ with $k - \log_2 \abs G \asymp k$, that the group giving rise to the largest diameter amongst all groups is $\mbz_2^d$.
\end{itemize}

\subsection{Concentration for $k \gtrsim \log \abs G$}
\label{sec-p3:diam:conc}

Recall that in \cref{res-p3:typdist2:res} we showed, in the regime $k \asymp \log n$ and under some assumptions, that, up to subleading order terms, the typical distance concentrates at $\alpha k$, for some constant $\alpha$.
The next theorem shows, in the same set-up, that the diameter does the same.
The argument uses the typical distance result as a `black box', then extending from this to diameter.

\begin{thm}
\label{res-p3:diam:conc:res}
	Let $(k_N)_\Ninn$ be a sequence of positive integers and $(G_N)_\Ninn$ a sequence of finite, Abelian groups;
	for each $\Ninn$, define $Z_{(N)} \cq [Z_1, ..., Z_{k_N}]$ by drawing $Z_1, ..., Z_{k_N} \sim^\iid \Unif(G_N)$.
	
	Suppose that $(k_N, G_N)_\Ninn$ satisfies either Hypotheses~\ref{hyp-p3:typdist2} or~\ref{hyp-p3:typdist3}.
	For $\lambda \in (0,\infty)$, let $\alpha^\pm_\lambda \in (0,\infty)$ be the constant from \cref{res-p3:typdist2:res};
	for each $\Ninn$, write $\rho_N \cq \log k_N / \log\log \abs{G_N}$, so that $k_N = (\log \abs{G_N})^{\rho_N}$.
	Then
	the following convergences in probability hold:
	\begin{alignat*}{2}
		\diam G_N(Z_{(N)}) / \rbb{ \alpha^\pm_\lambda k_N } &\to^\mbp 1
	&\Qwhen&
		\limt{N} k_N / \log \abs{G_N} = \lambda \in (0,\infty);
	\\
		\diam G_N(Z_{(N)}) / \rbb{ \tfrac{\rho_N}{\rho_N-1} \log_{k_N} \abs{G_N}} &\to^\mbp 1
	&\Qwhen&
		\limt{N} k_N / \log \abs{G_N} = \infty.
	\end{alignat*}
	Moreover,
		the implicit lower bound on the diameter holds deterministically, ie for all choices of generators, and for all Abelian groups,
	and,
	when $k \gg \log \abs G$,
		the implicit upper bound holds for all groups, not just Abelian groups.
\end{thm}

\begin{rmkt*}
	We only state and prove the result for $k \gtrsim \log \abs G$, but the argument can be extended to allow $k \ll \log \abs G$, provided $\log \abs G / k$ diverges sufficiently slowly.
	This requires a little more care; we do not explore the details here.
\end{rmkt*}

As always, we drop the $N$-subscripts in the proof, eg writing $\diam G_k$ or $\abs G$.

\begin{Proof}[Proof of \cref{res-p3:diam:conc:res}]
Clearly $\diam G_k = \mcd_k(1) \ge \mcd_k(\beta)$ for all $\beta \in [0,1]$.
Hence typical distance is trivially a lower bound on the diameter.
It remains to consider the upper bound.

\smallskip

Assume first \cref{hyp-p3:typdist2}, so $k \eqsim \lambda \log \abs G$ for some $\lambda \in (0,\infty)$.
Let $\eps \ll 1$, vanishing slowly and specified later.
Define $\alpha \cq \alpha^\pm_\lambda$ as in \cref{res-p3:typdist2:res}.
Let $A \cq [Z_1, ..., Z_{(1-\eps)k}]$ be the first $(1-\eps)k$ generators and $B \cq [Z_{(1-\eps)k+1}, ..., Z_k]$ be the remaining $\eps k$.
By transitivity, it suffices to consider distances from the identity.
The idea is to
	take $L$ steps using $A$ and then one more using $B$,
	where $L$ is the minimal radius of a ball in the $\abs A$-dimensional lattice of volume at least $n e^\omega$, for some slowly diverging $\omega$.
Write $M \cq \alpha k$.
By \cref{res-p3:typdist2:balls},
we have
\(
	L/M \eqsim 1 - \eps \eqsim 1.
\)
%
The key point is that when $k \asymp \log \abs G$ replacing $k$ with $(1 - \eps)k$ changes the typical distance by a factor $1 + o_{\eps\to0}(1)$.

\smallskip

By \cref{res-p3:typdist2:res}, \whp, $A$ is \textit{typical} in the sense that the proportion of elements of the group which can be reached via a word of length at most $L$, using only the generators from $A$, is at least $1 - e^{-\nu}$, for some $\nu \gg 1$, independent of $\eps$.

Condition on $A$, and that it is typical;
write $\widebar \mbp$ for the probability measure induced by this conditioning.
Denote by $H$ the set of elements which can be reached in the above sense.
(This is the vertex set of the ball of radius $L$ in $G(A)$.)
Fix $x \in G$.
Note that if $b \sim \Unif(G)$, then
\[
	\widebar \mbp\rbb{ x \in b H } \ge 1 - e^{-\nu}
\Qwhere
	b H \cq \bra{ b \cdot h \mid h \in H }.
\]
Furthermore, if $b,b' \sim \Unif(G)$ are independent then the events $\bra{ x \in b H }$ and $\bra{ x \in b' H }$ are $\widebar \mbp$-independent; this is because we have conditioned on $A$, and so $H$ is a deterministic set.

Using the $\eps k$ generators from $B$, informally we get $\eps k$ Bernoulli trials to get to $x$ using $b H$ for $b \in B$, and each trial has success probability $1 - \oh1$.
Formally, write $\mcr$ for the set of elements reachable from the identity via a word of length at most $L + 1$ (ie the `range'); let $b'$ be an arbitrary element of $B$, so $b' \sim \Unif(G)$.
(Recall that the conditioning makes $H$ non-random.)
Then
\[
	\widebar \mbp\rbb{ x \notin \mcr }
\le
	\widebar \mbp\rbb{ x \notin B H }
=
	\widebar \mbp\rbb{ x \notin b H \: \forall \, b \in B }
=
	\widebar \mbp\rbb{ x \notin b' H }^{\abs B}
\le
	e^{ - \nu \eps k }.
\]
Since $\nu \to \infty$, we may choose $\eps \to 0$ so that $\nu \eps \to \infty$.
Then, since $k \asymp \log n$, we have
\[
	\widebar \mbp\rbb{ \mcr \ne G }
=
	\widebar \mbp\rbb{ \exists \, x \in G \st x \notin \mcr }
\le
	n \, \widebar \mbp\rbb{ x \notin \mcr }
\le
	n e^{-\nu \eps k}
=
	\oh1.
\]

Averaging over $A$ establishes an upper bound of $\diam G_k \le L + 1$ whp, and $L \le M(1 + \eps)$.

\smallskip

Finally consider \cref{hyp-p3:typdist3}, so $k \gg \log \abs G$.
Exactly the same argument holds here, using the typical distance to first get to almost all the elements and then one more step.
Recall from \cref{res-p3:typdist3:res} that the upper bound is valid for arbitrary groups.
	%
\end{Proof}

\subsection{Universal Bounds for $k - \log_2 \abs G \asymp k$}
\label{sec-p3:diam:univ}

In this subsection we show that the group $\mbz_2^d$ gives rise to the random Cayley graph with the largest diameter when $k - \log_2 \abs G \asymp k$ whp, up to smaller order terms.

Recall that $\mfr(k, n)$ is the minimal $R \in \mbn$ with $\binom kR \ge n$.

\begin{thm}
\label{res-p3:diam:univ:res}
	Let $(k_N)_\Ninn$ be a sequence of positive integers and $(G_N)_\Ninn$ a sequence of finite groups;
	for each $\Ninn$, define $Z_{(N)} \cq [Z_1, ..., Z_{k_N}]$ by drawing $Z_1, ..., Z_{k_N} \sim^\iid \Unif(G_N)$.
	
	Suppose that $\liminf_N \rbr{ k_N - \log_2 \abs{G_N} } / k_N > 0$ and $\limsup_N \log k_N / \log \abs{G_N} = 0$.
	Then
	\[
		\limsupt{\Ninf}
		\diam G_N(Z_{(N)}) / \mfr(k_N, \abs{G_N})
	\le
		1
	\quad
		\text{in probability}.
	\]
	Further, if $G_N = \mbz_2^{d_N}$ for each $N$, then the diameter is given by $\mfr(k_N, \abs{G_N})$:
	\[
		\diam G_N(Z_{(N)}) / \mfr(k_N, \abs{G_N})
	\to^\mbp
		1.
	\]
\end{thm}

\begin{Proof}
From \cref{res-p3:typdist3:balls,res-p3:diam:conc:res},
when $k \gg \log \abs G$, the diameter concentrates at $\mfr(k, \abs G)$ when the underlying group is Abelian, and this is an upper bound for all~groups.

Thus it remains to consider $k$ with $k - \log_2 \abs G \asymp k$ and $k \asymp \log \abs G$.
All that was required for the upper bound on typical distance when $k \gg \log \abs G$ was that $\pr{ \WW = \WW' } \ll 1/\abs G$ where $\WW,\WW' \sim^\iid \Unif\rbr{ S_{k,1}(D) \cap B_{k,\infty}(1) }$ with $D \cq \widebar \mfd (1 + \xi)$, where $\widebar \mfd$ was the candidate typical distance radius and $\xi > 0$ was a constant.
We show that the analogous statement holds here.

Let $\xi > 0$ be fixed and set $R \cq \mfr(k, \abs G) (1 + \xi)$.
Before proceeding, let us determine some estimates on $\mfr$.
Let $h : (0,1) \to (0,1) : p \mapsto - p \log p - (1 - p) \log(1 - p)$ denote the binary entropy function (in nats).
It is standard that Stirling's approximation, like in \cref{res-p3:typdist3:balls}, gives
\[
	\binomt kr
=
	\expb{ k \, h(r/k) \cdot \rbr{ 1 + \oh1 } }.
\]
Thus, if $k - \log_2 \abs G \asymp k$, then we see that $\mfr(k, \abs G) \asymp k$.
Further, the fact that the derivative of $h$ is continuous and strictly positive on $(0,\tfrac12)$ gives $\binomt kR \gg \abs G$; hence $\pr{W = W'} \ll 1/\abs G$.

This shows that the typical distance $\mcd_k(\beta) \le \mfr(k, \abs G)$ whp up to smaller order terms for all constants $\beta \in (0,1)$.
This is then converted from a statement about typical distance to one about the diameter via the same method as used previously (in \S\ref{sec-p3:diam:conc}), noting that $\mfr(k, \abs G) \asymp k$.

Finally, we need to show a matching lower bound when $G = \mbz_2^d$, the hypercube.
A generator need never be applied more than once here, as all elements have order $2$.
Thus,
for $R \le \tfrac12 k$,
we~have%
\[
	\abs{ \mcb_k(R) }
\le
	\sumt[R]{r=0}
	\binomt kr
\le
	(R + 1) \binomt kR
\le
 	k \binomt kR
\quad
	\text{however the generators are chosen}.
\]
But, it is clear from the above asymptotic form for $\binom kr$ that
\(
	k \binomt kR
\ll
	\abs G
\)
if $R \cq \mfr(k, \abs G) (1 - \xi)$,
by continuity of the binary entropy function $h(\cdot)$.
Hence, the diameter is greater than $R$.
\end{Proof}

\section{Spectral Gap}
\label{sec-p3:gap}

In this section, we calculate the spectral gap; see \cref{res-p3:intro:gap}.
We first prove it for $k \ge 3 d(G)$.
In \S\ref{sec-p3:gap:ext}, we explain how to extend to $k - 2 d(G) \asymp k$ and then to $k - d(G) \asymp k$ for a density-$(1-\eps)$ subset of values for $\abs G$.
The lower bound holds deterministically, without any~conditions.

\subsection{Precise Statement}

For an Abelian group $G$, we write $d(G)$ for the minimal size of a generating set.
It is convenient to phrase the statement in terms of the \textit{relaxation time}, which is the inverse of the spectral gap.

\begin{thm}[Spectral Gap]
\label{res-p3:gap:res}
\begin{subequations}
\label{eq-p3:gap}
	First,
	there exists an absolute constant $c > 0$ so that,
	for	all Abelian groups $G$
	and
		all multisets $z$ of generators of size $k$,
	we have
	\[
		\trel^*\rbb{ G^-(z) }
	\ge
		\trel\rbb{ G^-(z) }
	\ge
		c \abs G^{2/k}.
	\label{eq-p3:gap:res:lower}
	\nt
	\]
	
	Second,
	for all $\delta > 0$,
	there exist constants $c_\delta, C_\delta > 0$ so that,
	for all Abelian groups $G$,
	if $k \ge (2 + \delta) d(G)$ and $Z_1, ..., Z_k \sim^\iid \Unif(G)$, then
	\[
		\pr{ \trel^*\rbr{ G^-_k } \le C_\delta \abs G^{2/k} }
	\ge
		1 - C_\delta 2^{-k/c_\delta}.
	\label{eq-p3:gap:res:upper}
	\nt
	\]
	Furthermore,
	for all $\eps \in (0,1)$,
	there exists a subset $\mba \subseteq \mbn$ of density at least $1 - \eps$ so that
	if $\abs G \in \mba$ then
	then condition $k \ge (2 + \delta) d(G)$ can be relaxed to $k \ge (1 + \delta) d(G)$ and \cref{eq-p3:gap:res:upper} still holds; the constant $C_\delta$ now also depends on $\eps$, ie becomes $C_{\delta, \eps}$, but $c_\delta$ need not be adjusted.
\end{subequations}
\end{thm}

We prove this for the non-absolute spectral gap, ie $\min_{\lambda \ne 1} \bra{1 - \lambda}$, where the minimum is over eigenvalues; the same proof also works for the absolute spectral gap, ie $\min_{\lambda \ne 1} \bra{1 - \abs \lambda}$.

\subsection{Lower Bound on Relaxation Time}
\label{sec-p3:gap:lower}

In this subsection, we establish the lower bound on the relaxation time in \cref{res-p3:gap:res}, ie \cref{eq-p3:gap:res:lower}.

\begin{Proof}[Proof of Lower Bound in \cref{res-p3:gap:res}]
Write $n \cq \abs G$.
Abbreviate \textit{simple random walk} by \textit{SRW}.
We may assume that $k \le \log_3(\tfrac12n)$, as otherwise \cref{eq-p3:gap:res:lower} indeed holds for some $c>0$.
Indeed, $\trel \gtrsim 1$ and $n^{2/k} \asymp 1$ if $k \gtrsim \log n$.
Let $L \cq \floor{ \tfrac12 \rbr{ (\tfrac12 n)^{1/k} - 1 } }$.
By our assumptions, $L \ge 1$.
Consider the set
\[
	A
\cq
	\brb{ w \cdot Z \mid w \in \mbz^k \text{ and } \abs{w_i} \le L \: \forall \, i = 1,...,k }
\subseteq
	G.
\label{eq-p3:gap:A-def}
\nt
\]
Clearly $\abs A \le (2L+1)^k \le \tfrac12 n$.
Let $t \ge 0$, and let $(Y_s)_{s\ge0}$ be a continuous-time rate-1 SRW on $\mbz$.
Write $\tau_{A^c}$ for the first exit time of a set $A$ by SRW on $G$ started from the identity.
Observe that
\[
	\pr[0]{ \tau_{A^c} > t }
\ge
	\pr[0]{ \maxt{s \in [0,t/k]} \abs{Y_s} \le L }^k,
\label{eq-p3:gap:escapeA}
\nt
\]
where $0 \in A$ is the identity of the group.
\cref{res-p3:gap:exitinterval} below provides a lower bound on the probability that the exit time from $\bra{ -L, ..., L }$ by $Y$ is at least $t/k$.
It gives
\[
	\pr[0]{ \maxt{s \in [0,t/k]} \abs{Y_s} \le L }
\ge
	\expb{ -\tfrac18 \pi^2 (t/k) / (L+1)^2 }.
\]
Substituting this into \cref{eq-p3:gap:escapeA} we get
\[
	\pr[0]{ \tau_{A^c} > t }
\ge
	\expb{ -\tfrac18 t \pi^2 / (L+1)^2 }.
\label{eq-p3:gap:escapeA2}
\nt
\]

The minimal Dirichlet eigenvalue of a set $A$ is defined to be the minimal eigenvalue of minus the generator of the walk killed upon exiting $A$; we denote it by $\lambda_A$.
For connected $A$,
\cref{res-p3:gap:diric} below states that, for all $a \in A$, we have
\[
	- \tfrac1t \log \pr[a]{ \tau_{A^c} > t }
\to	
	\lambda_A
\quad
	\text{as $\tinf$}.
\]
From this and \cref{eq-p3:gap:escapeA2}, it then follows that $\lambda_A \le \lambda$ where
\[
	\lambda
\cq
	\tfrac18 \pi^2 / (L+1)^2
\le
	\pi^2 / \rbb{ (\tfrac12n)^{1/k} + 1 }^2.
\]
\cite[Corollary~3.34]{AF:book} controls the relaxation time in terms of quasi-stationary hitting times.
It gives
\[
	\trel
\ge
	\rbr{ 1-\tfrac1n \abs A } / \lambda
\ge
	1/(2\lambda).
\]

This concludes the proof of the lower bound,
modulo \cref{res-p3:gap:exitinterval,res-p3:gap:diric}.
\end{Proof}

It remains to state and prove the quoted \cref{res-p3:gap:exitinterval,res-p3:gap:diric}.
The proofs are deferred~to~%
\mbox{\cite[\S D]{HOt:rcg:supp}}.

\begin{lem}
\label{res-p3:gap:exitinterval}
	Let $\ell \in \mbn$ and $\tau \cq \inf\bra{ s \ge 0 \mid \abs{Y_s} = \ell }$, where $(Y_s)_{s\ge0}$ is a continuous-time rate-1 simple random walk on $\mbz$.
	Let $\theta \cq \tfrac12 \pi / \ell$ and $\lambda \cq 1 - \cos\theta$.
	Then, for all $s \ge 0$, we have
	\[
		\pr[0]{ \tau > s }
	\ge
		e^{- \lambda s}
	\ge
		\expb{- \tfrac18 s (\pi/\ell)^2 }.
	\]
\end{lem}

The proof of the next lemma follows a standard quasi-stationarity argument.
For a transition matrix $P$ and a set $A$, let
	$\tau_{A^c}$ be the exit time of $A$
and
	$\lambda_A$ be the \textit{minimal Dirichlet eigenvalue}, defined to be the minimal eigenvalue of minus the generator of the chain killed upon exiting $A$,~ie~of
\[
	I_A - P_A
\Qwhere
	(I_A - P_A)(x,y) \cq \oneb{x,y \in A} \rbb{ \one{x = y} - P(x,y) }.
\]

\begin{lem}
\label{res-p3:gap:diric}
	Consider a rate-1, continuous-time, reversible Markov chain with transition matrix $P$.
	Let $A$ be a connected set, and let $\lambda_A$ and $\tau_{A^c}$ be as above.
	Then, for all $a \in A$,
	we have
	\[
		- \tfrac1t \log \pr[a]{ \tau_{A^c} > t }
	\to
		\lambda_A
	\quad
		\text{as $\tinf$}.
	\]
\end{lem}

\begin{rmkt*}
	%
Our proof gives an explicit form for $c$ in \cref{eq-p3:gap:res:lower}.
If $k \ll \log n$, then we get
\[
	\trel \ge 2 \pi^{-2} \abs G^{2/k} \cdot \rbb{1 + \oh1}.
\]
Indeed, in this case, in the definition of the set $A$ in \cref{eq-p3:gap:A-def},
we can take
\(
	L \cq \floor{ \tfrac12(\eps n)^{1/k} }
\)
for any $\eps > 0$, making $\abs A/\abs G$ arbitrary small.
One can improve the constant by replacing $A$ with
\[
	\brb{ w \cdot Z \midb w \in \mbz^k \text{ and } \sumt[k]{i=1} \abs{w_i}^2 \le L(k,n) },
\]
where $L(k,n)$ is the maximal integer satisfying
\(
	\abs{ \bra{ w \in \mbz^k \mid \sum_{i=1}^k \abs{w_i}^2 \le L(k,n) } } \le \tfrac12 n.
\)
\end{rmkt*}

\subsection{Upper Bound on Relaxation Time}

\newcommand{\cc}{c}
\newcommand{\CC}{C}

In this subsection, we establish the upper bound on the relaxation time in \cref{res-p3:gap:res}, ie \cref{eq-p3:gap:res:upper}.
We prove it for the usual relaxation time $\trel$; the same proof applies to bound the absolute relaxation time $\trel^*$.
In particular, we bound the probability that $1 - \lambda_2$ is small;
a completely analogous calculation bounds the probability that $1 + \lambda_n$ is small.
We only present~the~former~calculation.

For ease of presentation, we assume first that $k \ge 3 d(G)$.
In \S\ref{sec-p3:gap:ext}, we explain how to relax this condition, to prove the complete theorem.

\begin{Proof}[Eigenstatistic Preliminaries]
\qedtriangle
Decompose $G$ as $\oplus_1^d \: \mbz_{m_j}$.
An orthogonal basis of  eigenvectors for $P$, the transition matrix of the corresponding discrete-time walk, is given by
\[
	\rbb{ f_x \mid x \in G }
\Qwhere
	f_x(y)
\cq
	\cos\rbb{ 2 \pi \sumt[d]{j=1} x_j y_j / m_j },
\]
with corresponding eigenvalues given by
\begin{gather*}
	\rbb{ \lambda_x \mid x \in G }
\Qwhere
	\lambda_x
=
	\tfrac1k \sumt[k]{i=1} \cos\rbb{ 2 \pi (\bar x \cdot Z_i) },
\\
	\text{where}
\quad
	\bar x_j = x_j / m_j \text{ for all } j = 1,...,d
\Qand
	\bar x \cdot Z_i
=
	\sumt[d]{j=1}
	x_j Z_i^j / m_j
\end{gather*}
is the standard inner-product on $\mbr^d$, where $Z_i^j$ is the $j$-th coordinate of the $i$-th generator $Z_i$; here we identify $\bar x$ and $Z_i$ with elements of $\mbr^d$ in a natural manner.
This can be verified via an elementary calculation. Alternatively, it can also be derived from the fact that $(f_x \mid x \in G)$ are the real parts of the characters of $G$ in the representation-theoretic sense.

Observe that $\lambda_0 = 1$.
Our goal is to bound
	$\min\bra{ 1 - \lambda_x \mid x \in G \setminus \bra{0} }$
from below. 
For $\alpha \in \mbr$, let $\bra{\alpha} $ be the unique number in $(-\tfrac12,\tfrac12]$ so that $\alpha - \bra{\alpha} \in \mbz$.
A simple calculation shows that
\(
	2(\pi \varphi)^2
\ge
	1 - \cos(2 \pi \varphi)
\ge
	\tfrac23 (\pi \varphi)^2
\)
for all $\varphi \in [-\tfrac12, \tfrac12]$; see
\cite[Lemma~D.1]{HOt:rcg:supp}.
It follows from this that
\[
	1 - \lambda_x
\ge
	\tfrac{ 2\pi^2}{3k}
	\sumt[k]{i=1} \bra{ \bar x \cdot Z_i }^2.
\qedhere \:\:
\label{eq-p3:gap:lax}
\nt
\]
\end{Proof}

\begin{Proof}[Outline of Proof for $G = \mbz_n$]
\qedtriangle
It is instructive to consider the case $d=1$, as it serves as motivation for the definitions of $s_*(x)$ and $A(s)$ below. If $d = 1$ then $G = \mbz_n$. Let $s = s(x) \cq n/\gcd(n,x)$. Observe that $\{\bar x Z_i\} \sim \Unif \{\{1/s\},\{2/s\}, ...,\{1\} \}$.
Recall that, here, $\{m/s\} \in (-\tfrac12, \tfrac12]$, as above.

Consider first the case $s = s(x) \le C n^{1/k}$.
In this case, $1 - \lambda_x \ge c n^{-2/k}$ provided that at least $q \cq \lceil c k s^2 n^{-2/k}\rceil $ of the generators $Z_i$ do not satisfy that $\{\bar x Z_i\} = 0$.
Hence $\pr{ 1-\lambda_x < c n^{-2/k} } \le \frac kq s^{q-k} \lesssim s^{-9k/10}$, where the last inequality holds provided that $c > 0$ is sufficiently small.
There are at most $\ell$ different $x \in G$ with $s(x) = \ell$.
The union bound then says that there is no $x$ with $s(x) \le C n^{1/k}$ such that $1 - \lambda_x < c n^{-2/k}$ whp.

Now consider the case $s = s(x) > C n^{1/k}$.
Let $Y_i \cq \ell-1$ if $|\{\bar x Z_i\}| \in J_\ell$ where $J_1 \cq [0, \tfrac1{2M}]$ and $J_\ell \cq (\tfrac{\ell-1}{2M}, \tfrac{\ell}{2M}]$ for $\ell>1$, where $M \cq \lceil 4n^{1/k} \rceil$.
Then, $Y_i/M^2 \le Y_i^2/M^2 \lesssim \{\bar x Z_i\}^2$.
Finally, a simple combinatorial calculation, which we later present, gives $\pr{ \tfrac{1}{k} \sumt[k]{i=1} Y_i \le 1/10 } \le 2^{-k}/n$.
The proof can now be concluded by a union bound over all $x$ such that $s(x) > C n^{1/k}$.

In the above calculation, we obtain a better upper bound on $\pr{ 1-\lambda_x < c n^{-2/k} }$ when $s(x)$ is large.
When $d > 1$, loosely speaking, our argument allows us to reduce the analysis to the case that $G = \mbz_{m_j}$ for any $j \in [d]$.
This reduces the analysis to the one above, with the quantity $n$ above remaining the same, rather than taking the value $m_j$; above, the fact that $n = |G|$ only played a role in bounding $|\{x \in G \mid s(x) = \ell\}|$.
By the above analysis, we want to pick $j$ such that $m_j/\gcd(x_j,m_j) = \max_{i \in [d]} m_i/\gcd(x_i,m_i) \eqqcolon s_*(x)$.
We consider the two cases $s_*(x)> C n^{1/k}$ and $s_*(x) \le C n^{1/k}$ and derive the same estimates on $\pr{ 1-\lambda_x < c n^{-2/k} }$, with $s_*(x)$ playing here the role of $s(x)$ above, as in the case $d = 1$ outlined above.

To conclude by a union bound, we also require an upper bound on $r_s \cq |\{x \in G \mid s_*(x) = s\}|$ from which we can deduce that $\sum_{s > C n^{1/k}} r_s s^{-9k/10} = \oh1$.
This is where we require the assumption that $k \ge 3d$---or with a bit more care, that $k \ge (2+\eps)d$ for some constant $\eps >0$.
\end{Proof}


\begin{Proof}[Proof of Upper Bound in \cref{res-p3:gap:res}]
For each $x \in G$, we make the following definitions:
\begin{alignat*}{2}
	g_j \cq g_j(x) &\cq \gcd(x_j,m_j)
	\quad
	&&\text{for each $j \ge 1$};
\\
	s_* \cq s_*(x) &\cq \max\brb{m_j/g_j \mid j \in \{1,...,d\} };
\\
	A(s) &\cq \brb{x \in G \mid s_*(x)=s }
	\quad
	&&\text{for each $s \ge 1$};
\\
	\phi(j) &\cq \absb{ \brb{j' \in \{1,...,j\} \mid \gcd(j,j') = 1 } }
	\quad
	&&\text{for each $j \ge 1$}.
\end{alignat*}
From this, we claim that we are able to deduce, for $s \ge 2$, that
\[
	\abs{ A(s) }
\le
	\rbb{ \sumt[s]{j=1} \phi(j) }^d
\le
	\rbb{ 1 + \sumt[s]{j=2} (j-1) }^d
\le
	\rbb{ \tfrac12 s^2 }^d.
\label{eq-p3:gap:A(s):orig}
\nt
\]
Indeed, $\phi(j) \le j-1$ for $j \ge 2$, and observe that
\[
	\text{if $r$ divides $m$,}
\Quad{then}
	\absb{ \brb{ a \in \bra{1,...,m} \midb \gcd(a,m) = r } }
=
	\phi(m/r);
\]
hence, summing over the set of possible values for $m_j/g_j$, which by definition of $A(s)$ is $\{1, ..., s\}$, we have $\abs{A(s)}^{1/d} \le \sum_{j=1}^s \phi(j)$.
We are then able to deduce the upper bound, ie \cref{eq-p3:gap:res:upper}, from \cref{res-p3:gap:s_*}, which we state precisely below.
Indeed, first write
\[
	p(s)
\cq
	\MAX{x : s_*(x) = s}
	\pr{ 1 - \lambda_x \le c n^{-2/k} }.
\]
We control this probability using \cref{eq-p3:gap:lax,eq-p3:gap:A(s):orig},
along with \cref{res-p3:gap:s_*} below, which states that
\[
	\prb{ \tfrac1k \sumt[k]{i=1} \bra{ \bar x \cdot Z_i }^2 \le \cc n^{-2/k} }
\le
\begin{cases}
	s_*(x)^{-9k/10}
		&\text{where}\quad
	s_*(x) \le \CC n^{1/k},
\\
	2^{-k}/n
		&\text{where}\quad
	s_*(x) > \CC n^{1/k},
\end{cases}
\]
for some absolute constants $\cc$ and $\CC$.
Applying these and letting $\cc' \cq \cc \cdot \tfrac3{2\pi^2}$ gives
\[
	\sumt{x \in G \setminus \{0\} }
	\pr{ 1 - \lambda_x \le \cc' n^{-2/k} }
&\le
	n \maxt{s > \CC n^{1/k}} p(s)
+	\sumt{2 \le s \le \CC n^{1/k}} \abs{A(s)} \, p(s)
\\&
\le
	2^{-k} + 2^{-d} \sumt{s \ge 2} s^{2d} (2s)^{-9k/10}
\lesssim
	2^{-k},
\]
where we have used $k \ge 3d$ and the fact that $s_*(x) > 1$ for all $x \ne 0$.
Hence, by the union bound,
\[
	\pr{ \trel \ge n^{2/k} / \cc' }
=
	\pr{ \gamma \le \cc' n^{-2/k} }
=
	\pr{ \exists \, x \in G \st 1 - \lambda_x \le \cc' n^{-2/k} }
\lesssim
	2^{-k}.
\]

This concludes the proof of the upper bound when $k \ge 3 d(G)$,
modulo \cref{res-p3:gap:s_*}.
\end{Proof}

It remains to state and prove the quoted \cref{res-p3:gap:s_*},
then extend the range of $k$.

\begin{prop}
\label{res-p3:gap:s_*}
	There exist absolute constants $\cc \in (0,1)$ and $\CC$ such that 
	\begin{subequations}
		\label{eq-p3:gap:s_*}
	\begin{empheq}[ left = {%
		\prb{ \tfrac1k \sumt[k]{i=1} \bra{ \bar x \cdot Z_i }^2 \le \cc n^{-2/k} }
	\le
		\empheqlbrace}]
	{alignat=2}
		&s_*(x)^{-9k/10}
			&&\Qwhere
		s_*(x) \le \CC n^{1/k},
		\label{eq-p3:gap:s_*:<}
	\\
		&2^{-k}/n
			&&\Qwhere
		s_*(x) > \CC n^{1/k}.
		\label{eq-p3:gap:s_*:>}
	\end{empheq}
	\end{subequations}
\end{prop}

\begin{Proof}
Fix $x \in G$. First consider the case that $s \cq s_*(x) > \CC n^{1/k}$, ie \cref{eq-p3:gap:s_*:>}.
Let $j \cq j(x)$ be a coordinate satisfying $s = m_j/g_j$. Denote $m \cq m_{j(x)}$ and $g \cq g_{j(x)}$.
Observe that $x_j Z_i^j \sim^\iid \Unif\{g,2g,...,m\}$ for each $i$.
Hence, for each $i$, we have
\[
	U_i
\cq
	\bar x_j Z_i^j
\sim
	\Unif\bra{1/s, 2/s, ..., 1}.
\label{eq-p3:gap:Uiunif}
\nt
\]

By averaging over $(a_i)_{i=1}^k$, where
\(
	a_i
\cq
	\bra{ \sumt{\ell \in \bra{1,...,d} \setminus \bra{j} } x_\ell Z_i^\ell / m_\ell },
\)
recalling that $\bra{\alpha}$ is the unique number in $(-\tfrac12, \tfrac12]$ so that $\alpha - \bra{\alpha} \in \mbz$,
it suffices to show that
\[
	\MAX{b_1,...,b_k \in [-1/2,1/2]}
	\pr{ \tfrac1k \sumt[k]{i=1} \bra{ U_i + b_i }^2 \le \cc n^{-2/k} }
\le
	2^{-k}/n.
\label{eq-p3:gap:Ui}
\nt
\]
Replacing $\cc$ with $4\cc$ we may assume that $b_i \in \tfrac1s \mbz$ for all $i$.
Indeed, if
\[
	\abs{b_i - \ell/s} \le 1/(2s),
\Quad{ie}
	\abs{b_i - \ell/s}
=
	\min\brb{ \abs{b_i - \alpha} \mid \alpha \in \tfrac1s \mbz },
\]
then
\(
	\bra{ U_i + \ell/s }^2 \le 4 \bra{ U_i +b_i }^2.
\)
Hence
\[
	\text{if}
\quad
	\tfrac1k \sumt[k]{j=1} \bra{ U_i + b_i }^2 \le \cc n^{-2/k}
\quad
	\text{then}
\quad
	\tfrac1k \sumt[k]{j=1} \bra{ U_i + \ell/s }^2 \le 4 \cc n^{-2/k}.
\]
If $b_i \in \frac1s \mbz$, then $\bra{ U_i + b_i }$ has the same law as $\bra{ U_i }$ by \cref{eq-p3:gap:Uiunif}.
Hence, we may assume $b_i = 0$ for~all~$i$.

We now split $[0,\tfrac12]$ into $M \cq \ceil{4 n^{1/k}}$ consecutive intervals of equal length $J_1,...,J_M $,
where $J_1 \cq [0,\frac{1}{2M}]$ and $J_\ell \cq (\frac{\ell-1}{2M}, \frac{\ell}{2M}]$ for $\ell >1$.
Let $Y_i \cq \ell-1$ if $\abs{\{U_i\}} \in J_{\ell}$.
Clearly, $\tfrac14 Y_i/M^2 \le \tfrac14 Y_i^2/M^2 \le \bra{U_i}^2$.
It thus suffices to show that
\[
	\pr{ \tfrac1k \sumt[k]{i=1} Y_i \le \tfrac1{10} }
\le
	2^{-k}/n.
\]
This last claim follows by a simple counting argument:
	there are $M^k$ total assignments of the $Y_i$-s,
	but at most $L(k) \cq \binom{\ceil{11k/10}}{k-1} \le 2^k$ assignments satisfy
	$\tfrac1k \sum_{i=1}^k Y_i \le \tfrac1{10}$,
	since $L(k)/M^k \le 2^{-k}n^{-1}$.

\smallskip

We now prove the case $s \cq s_*(x) \le \CC n^{1/k}$, ie \cref{eq-p3:gap:s_*:<}.
By the same reasoning as for \cref{eq-p3:gap:Ui},
it suffices to show that
\[
	\MAX{b_1, ..., b_k \in [-1/2,1/2]}
	\pr{ \tfrac1k \sumt[k]{i=1} \bra{ U_i + b_i }^2 \le \cc n^{-2/k} }
\le
	s^{-9k/10}.
\label{eq-p3:gap:Ui2}
\nt
\]
Regardless of $b_i$, there is at most one $a \cq a(b_i) \in \{1/s, 2/s, ..., 1\}$ such that $\bra{ a + b_i }^2 < (2s)^{-2}$,
and hence by \cref{eq-p3:gap:Uiunif}, for all $i$, we have
\[
	\pr{ \bra{ U_i +b_i }^2 < (2s)^{-2} } \le 1/s.
\]
If there is no such value $a(b_i)$, then set $a(b_i) \cq -1$.

If $\bra{ U_i +b_i }^2 \ge (2s)^{-2}$ for at least $q \cq k \cdot 4 \cc s^2 n^{-2/k}$ of the $i$-s,
ie
\[
	\text{if}
\quad
	\absb{ \brb{ i \in \{1,...,k\} \mid U_i \ne a(b_i) } } \ge q,
\Quad{then}
	\tfrac1k \sumt[k]{i=1} \bra{ U_i + b_i }^{2} \ge \cc n^{-2/k}.
\]
as desired.
As $s \le \CC n^{1/k}$, by taking $\cc$ sufficiently small in terms of $\CC$, we can make $q/k$ sufficiently small so that the following holds:
\[
	\pr{ \absb{ \brb{ i \in \{1,...,k\}  \mid U_i \ne a(b_i) } } < q }
\lesssim
	\binomt{k}{q}s^{q-k}
\lesssim
	s^{-9k/10}.
\qedhere
\]
\end{Proof}


\subsection{Relaxing the Conditions on $k$}
\label{sec-p3:gap:ext}

In this subsection, we explain how to relax the conditions on $k$.
First we can relax from $k \ge 3 d(G)$ to $k - 2 d(G) \asymp k$, valid for every group size $n = \abs G$.

\smallskip

We now give conditions under which this can be relaxed to $k - d(G) \asymp k$.
If $G = \mbz_p^d$ for a prime $p$, then one can relax this further to $k - d \gtrsim d$, and even allow $k - d(G) \ll d(G)$, provided $p$ diverges.
	(In this case, the term  $2^{-k}$ has to be replaced by another term which tends to zero at a slower rate as $k \to \infty$.)
This follows from the fact that now we only need to consider \cref{eq-p3:gap:A(s):orig} above with $s \cq p$ and we can replace \cref{eq-p3:gap:A(s):orig} with $\abs{A(p)} = p^d - 1$.
So the condition $k - d(G) \asymp k$ is sufficient when $G = \mbz_p^d$ with $p$ prime.

We now show that if $\abs G$ is `typical' (in a precise sense), then the same condition is sufficient.
In the proof above, in \cref{eq-p3:gap:A(s):orig}, we used the crude bound
\[
	\abs{A(s)}
\le
	\rbb{ \sumt{i \in [s]} \phi(i) }^d
\le
	\rbb{\tfrac12 s^2}^d.
\]
Instead, recalling that we write $i \wr n$ to mean that $i$ divides $n$, we can use the improved bound
\[
	\abs{A(s)}
\le
	\rbb{ \sumt{i \in [s]} i \, \one{i \wr n} }^d.
\]
In
\cite[Lemma~F.7]{HOt:rcg:supp},
we show that,
for all $\eps > 0$,
there exists a constant $C'_\eps$ and a density-$(1-\eps)$ set $\mbb_\eps \subseteq \mbn$ such that,
for
	all $n \in \mbb_\eps$
and
	all $2 \le s \le n$,
we have
\[
	\sumt{i \in [s]}
	i \, \one{i \wr n}
\le
	C'_\eps s (\log s)^2.
\]
Using this to derive an improved bound on $\abs{A(s)}$, and adjusting some of the constants in the proof in an appropriate manner, an inspection of the proof reveals that, for all $n \in \mbb_\eps$ and all $\delta > 0$, there exists a positive constant $C_{\eps,\delta}$ so that, for all Abelian groups of size $n$, if $k \ge (1 + \delta) d$, then
\[
	\pr{ \trel(G_k) \ge C_{\eps,\delta} n^{2/k} } \le e^{-k/C_{\eps,\delta}}.
\]

\section{Open Questions and Conjectures}
\label{sec-p3:open-conj}

We close the paper with some questions which are left open.

\newcounter{oq}

\refstepcounter{oq}
\label{oq-p3:typdist}
\subsubsection*{\theoq: Typical Distance and Diameter for All Abelian Groups}

In our typical distance theorem, there were some conditions on the group.
We allowed any group with $d(G) \ll \log \abs G / \log\log k$ if $1 \ll k \ll \log \abs G$, but once $d(G)$ became larger than this or $k$ became order $\log \abs G$, we had to impose conditions.
We conjecture that these are artefacts~of~the~proof.

\begin{conj-ind}\theoq
	Let $G$ be an Abelian group.
	Suppose that $1 \ll k \lesssim \log \abs G$ and $k - d(G) \gg 1$.
	Then the typical distance statistic concentrates
	at a value which depends only on $k$ and $G$, not the particular realisation of the generators.
	Further, if $k \ll \log \abs G$ and $k - d(G) \asymp k$, then it concentrates at a value which depends only on $k$ and $\abs G$.
\end{conj-ind}

The claim when $1 \ll k \ll \log \abs G$ and $k - d(G) \asymp k$ is a natural extension of \cref{res-p3:typdist1:res}. Further, if $k \ll \sqrt{\log \abs G / \log\log\log \abs G}$, then $k - d(G) \gg 1$ is sufficient, by \cref{hyp-p3:typdist1}.
Once we relax to $k - d(G) \gg 1$, for larger $k$, we still expect concentration of typical distance for all Abelian groups, but now the value will likely depend on the specific group.
Compare this with the occurrence of cutoff for the random walk on the random Cayley graph established in \cite{HOt:rcg:abe:cutoff}.

\smallskip

There are two levels on which concentration occurs:
	first, for a fixed graph $G(z)$, one draws a $U \sim \Unif(G)$ and looks for concentration of $\dist(\id,U)$ at some value, say $f(z)$;
	second, one draws $Z$ uniformly and looks for concentration of $f(Z)$.
The second is the meat of Conjecture~\ref{oq-p3:typdist}.
Indeed, our lower bound on typical distance holds for all Abelian groups and all Cayley graphs with $k$ generators, thus necessarily $\pr[G(z)]{\dist(\id,U) \gtrsim k \abs G^{1/k}} = 1 - \oh1$ for all such Cayley graphs~$G(z)$.
Additionally, our spectral gap estimate (\cref{res-p3:intro:gap}) says that the gap is order $\abs G^{-2/k}$ if $k - 2 d(G) \asymp k$ (or when $k - d(G) \asymp k$ and $\abs G$ is `typical') whp over uniform $Z$.

Since $u \mapsto \dist(\id,u)$ is a 1-Lipschitz function, by Poincaré's inequality $\VAR[G(z)]{\dist(\id,U)} \le \trel(G(z))$.
For all multisets $z$ of size $k$ satisfying the aforementioned spectral gap estimate from \cref{res-p3:intro:gap} (which holds whp for $G_k$),
using our deterministic lower bound on the typical distance,
we see that $\dist_{G(z)}(\id,U)$ concentrates at some value $f(z)$, which may depend on $z$, by Chebyshev's inequality.
We conjecture that in fact $f(Z)$ concentrates at some value~$\mcd$.

\medskip

It is easy to see that the typical distance and diameter are always the same up to constants.
We conjecture that the diameter of $G_k$ concentrates whp whenever $1 \ll k \lesssim \log \abs G$ and $k - d(G) \gg 1$.
We leave open the question of finding conditions under which the diameter and typical distance are asymptotically equivalent whp.

%
%
%

\refstepcounter{oq}
\label{oq-p3:cheeger}
\subsubsection*{\theoq: Isoperimetry for Random Cayley Graphs}
 
The \textit{isoperimetric}, or \textit{Cheeger}, \textit{constant} of a finite $d$-regular graph $G = (V,E)$ is defined as
\[
	\Phi_*
\cq
	\tfrac1d
	\MIN{1 \le \abs S \le \frac12 \abs V}
	\Phi(S)
\Qwhere
	\Phi(S)
\cq
	\tfrac1{\abs S} \absb{ \brb{ \bra{ a, b } \in E \midb a \in S, \, b \in S^c } }.
\]
More generally, the isoperimetric constant is defined for Markov chains; see \cite[\S 7.2]{LPW:markov-mixing}. For a given stochastic matrix $P$, it is easy to see that the original chain $P$, the time-reversal $P^*$ and the additive symmetrisation $\tfrac12(P + P^*)$ all have the same isoperimetric profile.
Thus the isoperimetric constant for a directed Cayley graphs is the same as that for the undirected version.

The following conjecture asserts that the Cheeger constant is, up to a constant factor, the same as that of the standard Cayley graph of $\mbz_L^k$ where $L$ is such that $n \asymp L^k$.

\begin{conj-ind}\theoq
	There exists a constant $c$ so that,
	for all $\eps \in (0,1)$,
	there exist constants $n_\eps$ and $M_\eps$ so that,
	for every finite group $G$ of size at least $n_\eps$,
	when $k \ge M_\eps$,
	we have
	\[
		\pr{ \Phi_*(G_k) \le c \abs G^{-1/k} }
	\le
		\eps,
	\]
	where $\Phi_*(G_k)$ is the Cheeger constant of a random Cayley graph with $k$ generators.
\end{conj-ind}

By \cite[Theorem~6.29]{LP:prob-on-trees-networks}, which regards expansion of general Cayley graphs, along with out upper bound on typical distance (and hence on diameter), we can prove this conjecture up to a factor $k$.

By the well-known discrete analogue of Cheeger's inequality, discovered independently by multiple authors---see, for example, \cite[Theorem~13.10]{LPW:markov-mixing}---we have
$\tfrac12 \gamma \le \Phi_* \le \sqrt{ 2 \gamma }$.
Determining the correct order of $\Phi_*$ in our model remains an open problem. We conjecture that the correct order of $\Phi_*$ is given by $\sqrt{ \gamma }$, ie order $\abs G^{-1/k}$, using \cref{res-p3:intro:gap} for the order of the spectral gap.

\smallskip

The celebrated Alon--Roichman theorem states that the Cayley graph of any finite group $G$ is a $(1-\eps)$-expander (ie $\Phi_* \ge 1 - \eps$) whp when $k \ge C_\eps\log\abs G$, for some constant $C_\eps$; the best known upper bound on $C_\eps$ is $\Oh{1/\eps^2}$.
\textcite[Theorem~1.2]{N:cheeger-small-sets} refines this for Abelian groups:
	he showed that one can in fact bound $\abs{ \Phi(S) - 1 } \le \eps \sqrt{\log\abs{S}/\log\abs G}$ for all $S$ with $1 \le \abs S \le \tfrac12 \abs V$ simultaneously, when $k/\log n \ge C/\eps^2$, for a constant $C$.
In recent work, \textcite{SSZ:cayley-eigenbasis} extended Naor's result to all groups.

\renewcommand{\bibfont}{\sffamily}
\renewcommand{\bibfont}{\sffamily\small}
\printbibliography[heading=bibintoc]

\end{document}